%% file: newdual.tex
\documentclass{tran-l}

\input xy
\xyoption{all}
\xyoption{2cell}
\UseAllTwocells

\usepackage{amsmath}
\usepackage{amscd}

\def\dlim{\underset{\rightarrow}{\mbox{lim}}}

\def\HomO{{\mathcal{H}}{\it om}_{{\mathcal{O}}_{X}}}
\def\HomO2{{\mathcal{H}}{\it om}_{{\mathcal{O}}_{X^{2}}}}
\def\HomA{{\mathcal{H}}{\it om}_{{\sf Gr }\mathcal{A}}}
\def\HU{{\underline{{\mathcal{H}}{\it om}}_{{\sf Gr} \mathcal{A}}}}
\def\Ten{{\underline{\otimes}}_{\mathcal{A}}}
\def\Hom{{\operatorname{Hom}}}
\def\Homa{{\operatorname{Hom}}_{{\sf Gr}\mathcal{A}}}
\def\Homo{{\operatorname{Hom}}_{{\mathcal{O}}_{X}}}
\def\Ext{{\mathcal{E}}{\it xt}_{{\sf Gr}\mathcal{A}}}
\def\ExtU{\underline{{\mathcal{E}}{\it xt}}_{{\sf Gr}\mathcal{A}}}

\def\A{\mathcal{A}}
\def\C{\mathcal{C}}
\def\M{\mathcal{M}}
\def\N{\mathcal{N}}
\def\F{\mathcal{F}}
\def\D{\mathcal{D}}
\def\quota{\A/A_{\geq n}}
\def\R{\mathcal{R}}
\def\Q{\mathcal{Q}}
\def\B{\mathcal{B}}
\def\R{\mathcal{R}}
\def\E{\mathcal{E}}
\def\G{\mathcal{G}}

\copyrightinfo{2003}{American Mathematical Society}

\newtheorem{theorem}{Theorem}[section]
\newtheorem{lemma}[theorem]{Lemma}
\newtheorem{proposition}[theorem]{Proposition}
\newtheorem{corollary}[theorem]{Corollary}

\theoremstyle{definition}
\newtheorem{definition}[theorem]{Definition}

\theoremstyle{remark}
\newtheorem{remark}[theorem]{Remark}

\numberwithin{equation}{section}

\begin{document}

\title{Serre Duality For Non-commutative ${\mathbb{P}}^{1}$-bundles}

\author{Adam Nyman}
\address{Department of Mathematical Sciences, Mathematics Building, University of Montana, Missoula, MT 59812-0864}
\email{nymana@mso.umt.edu}


\date{}

\begin{abstract}
Let $X$ be a smooth scheme of finite type over a field $K$, let $\E$ be a locally free $\mathcal{O}_{X}$-bimodule of rank $n$, and let $\A$ be the non-commutative symmetric algebra generated by $\E$.  We construct an internal $\operatorname{Hom}$ functor, $\HU(-,-)$, on the category of graded right $\A$-modules.  When $\E$ has rank 2, we prove that $\A$ is Gorenstein by computing the right derived functors of $\HU(\mathcal{O}_{X},-)$.  When $X$ is a smooth projective variety, we prove a version of Serre Duality for ${\sf{Proj }} \A$ using the right derived functors of $\underset{n \to \infty}{\lim} \HU (\A/\A_{\geq n}, -)$.
\end{abstract}

\maketitle

\tableofcontents

\include{introduction}

\include{ncsym}

\include{hom}

\include{cohomology}

\include{serre}

\include{compats}

\include{proofs}

\include{bib}

\end{document}

%% file: introduction.tex
\section{Introduction}
Although the classification of non-commutative surfaces is nowhere in sight, some classes of non-commutative surfaces are relatively well understood.  For example, non-commutative deformations of the projective plane and the quadric in $\mathbb{P}^{3}$ have been classified in \cite{atv} and \cite{quad}, respectively.  

Non-commutative $\mathbb{P}^{1}$-bundles over curves have not been studied as extensively as non-commutative $\mathbb{P}^{2}$'s or non-commutative quadrics, but certainly play a prominent role in the theory of non-commutative surfaces.  For example, certain non-commutative quadrics are isomorphic to non-commutative $\mathbb{P}^{1}$-bundles over curves \cite{quad}.  In addition, every non-commutative deformation of a Hirzebruch surface is given by a non-commutative $\mathbb{P}^{1}$-bundle over $\mathbb{P}^{1}$ \cite[Theorem 7.4.1, p. 29]{p1bundles}.  Quotients of four-dimensional Sklyanin algebras by central homogeneous elements of degree two provide important examples of coordinate rings of non-commutative $\mathbb{P}^{1}$-bundles over $\mathbb{P}^{1}$ \cite[Theorem 7.2.1]{translations}.  

The purpose of this paper is to prove a version of Serre duality for non-commutative ${\mathbb{P}}^{1}$-bundles over smooth projective varieties of dimension $d$ over a field $K$.  Before describing this result in more detail, we review some important notions from non-commutative algebraic geometry.

If $X$ is a quasi-compact and quasi-separated scheme, then ${\sf Mod }X$, the category of quasi-coherent sheaves on $X$, is a Grothendieck category.  This leads to the following generalization of the notion of scheme, introduced by Van den Bergh in order to define a notion of blowing-up in the non-commutative setting.

\begin{definition} \cite{blowup}
A {\bf quasi-scheme} is a Grothendieck category ${\sf Mod }X$, which we denote by $X$.  $X$ is called a {\bf noetherian} quasi-scheme if the category ${\sf Mod }X$ is locally noetherian.  $X$ is called a {\bf quasi-scheme over $\mathbf{K}$} if the category ${\sf Mod }X$ is $K$-linear. 
\end{definition}
If $R$ is a ring and ${\sf Mod }R$ is the category of right $R$-modules, ${\sf Mod }R$ is a quasi-scheme, called the non-commutative affine scheme associated to $R$.  If $A$ is a graded ring, ${\sf Gr }A$ is the category of graded right $A$-modules, ${\sf Tors }A$ is the full subcategory of ${\sf Gr }A$ consisting of direct limits of right bounded modules, and ${\sf Proj }A$ is the quotient category ${\sf Gr }A/{\sf Tors }A$, then ${\sf Proj }A$ is a quasi-scheme called the non-commutative projective scheme associated to $A$.  If $A$ is an Artin-Schelter regular algebra of dimension 3 with the same hilbert series as a polynomial ring in 3 variables, ${\sf Proj }A$ is called a non-commutative ${\mathbb{P}}^{2}$.  These definitions motivate the following     

\begin{definition} \label{def.pbund} \cite{p1bundles}
Suppose $X$ is a smooth scheme of finite type over $K$, $\E$ is a locally free $\mathcal{O}_{X}$-bimodule of rank 2 and $\A$ is the non-commutative symmetric algebra generated by $\E$ \cite[Section 5.1]{p1bundles}.  Let ${\sf Gr} \A$ denote the category of graded right $\A$-modules, let ${\sf Tors }\A$ denote the full subcategory of ${\sf Gr }\A$ consisting of direct limits of right-bounded modules, and let ${\sf Proj }\A$ denote the quotient of ${\sf Gr }\A$ by ${\sf Tors }\A$.  The category ${\sf Proj }\A$ is a {\bf non-commutative $\mathbf{\mathbb{P}^{1}}$-bundle over $\mathbf{X}$.}  
\end{definition}
Although no generally accepted definition of ``smooth non-commutative surface" exists yet, it seems reasonable to insist that such an object should be a noetherian quasi-scheme of cohomological dimension $2$.  While an intersection theory for non-commutative surfaces exists (\cite{intersect}, \cite{intersect2}), it has yet to be applied to the study of non-commutative $\mathbb{P}^{1}$-bundles over curves.  Mori shows \cite[Theorem 3.5]{intersect2} that if $Y$ is a noetherian quasi-scheme over a field $K$ such that
\begin{enumerate}
\item{} 
$Y$ is $\operatorname{Ext}$-finite,

\item{}
the cohomological dimension of $Y$ is $2$, and 

\item{}
$Y$ satisfies Serre duality
\end{enumerate}   
then there is an intersection theory on $Y$ such that the Riemann-Roch theorem and the adjunction formula hold.  

We prove a version of Serre duality for non-commutative ${\mathbb{P}}^{1}$-bundles over smooth projective varieties of dimension $d$ over $K$.  In the course of the proof, we show that such a $\mathbb{P}^{1}$-bundle has cohomological dimension $d+1$.  Thus, the second and third items on the list above are verified for non-commutative $\mathbb{P}^{1}$-bundles over smooth projective curves.  The author is currently using the techniques of this paper to prove the first item on the list above holds for non-commutative $\mathbb{P}^{1}$-bundles.      

We now describe the main results of this paper.  Suppose $\A$ is as in Definition \ref{def.pbund}, $\pi: {\sf Gr }\A \rightarrow {\sf Proj }\A$ is the quotient functor, $\omega$ is right adjoint to $\pi$, and $\tau:{\sf Gr }\A \rightarrow {\sf Gr }\A$ is the torsion functor.  For any collection $\{\C_{ij}\}_{i,j \in \mathbb{Z}}$ of $\mathcal{O}_{X}$-bimodules, let $e_{i}\C$ denote the $\mathcal{O}_{X}$-bimodule $\underset{j}{\oplus}\C_{ij}$.

We prove the following version of Serre duality (Theorem \ref{theorem.serre}):
\begin{theorem} 
If $X$ is a smooth projective variety of dimension $d$ over $K$, there exists an object $\omega_{\A}$ in $\sf{Proj }\A$ such that for $0 \leq i \leq d+1$ there is an isomorphism
$$
\operatorname{Ext}_{{\sf Proj} \A}^{i}(\pi(\mathcal{O}_{X} \otimes e_{0}\A),\M)' \cong \operatorname{Ext}_{{\sf Proj} \A}^{d+1-i}(\M,\omega_{\A})
$$
natural in $\M$.  The prime denotes dualization with respect to $K$.
\end{theorem}
To prove this theorem, we apply the Brown representability theorem following J$\ddot{\operatorname{o}}$rgenson \cite{duality}.  In order to apply the Brown representability theorem to prove our version of Serre duality, we need to prove technical results (Theorem \ref{theorem.cohom} and Theorem \ref{theorem.cd}) regarding the cohomology of the functor
\begin{equation} \label{eqn.gamma0}
\operatorname{Hom}_{{\sf Proj }\A}(\pi(\mathcal{O}(i) \otimes e_{m}\A), -).
\end{equation}

Our strategy for studying the cohomology of (\ref{eqn.gamma0}) is indirect.  We construct a bifunctor, $\HU(-,-)$ (Definition \ref{def.hom}) whose left input is a locally free $\A-\A$ bimodule and whose right input and output are graded right $\A$-modules.  We prove the functor $\HU(-,-)$ enjoys the expected properties:
\begin{itemize}
\item{}
$\HU(\A,-) \cong \operatorname{id}_{{\sf Gr }\A}$,

\item{}
$\Homo(\mathcal{L}, \HU(\C,\M)_{m}) \cong \Homa(\mathcal{L} \otimes e_{m}\C, \M)$ for an quasi-coherent $\mathcal{O}_{X}$-module $\mathcal{L}$, and

\item{}
$\tau(-) \cong \underset{n \to \infty}{\lim} \HU (\A/\A_{\geq n}, -)$
\end{itemize}
(Propositions \ref{prop.homisom}, \ref{prop.global} and \ref{prop.tau}, respectively).  Since $\A$ is not a sheaf of algebras, $\HU(-,-)$ cannot be defined locally.  Instead, we use the natural categorical definition. 

We can use the first two properties above to find an alternative description for (\ref{eqn.gamma0}):

\begin{align*}
\operatorname{Hom}_{{\sf Proj }\A}(\pi(\mathcal{O}(i) \otimes e_{m}\A), -) & \cong \Homa(\mathcal{O}(i) \otimes e_{m}\A, \omega(-)) \\
& \cong \Homo(\mathcal{O}(i), \HU(\A,\omega(-))_{m}) \\
& \cong \Homo(\mathcal{O}(i), (\omega(-))_{m}) \\
& \cong \Gamma \circ \mathcal{O}(-i) \otimes (\omega(-))_{m}
\end{align*}
where $\Gamma:{\sf Mod }\mathcal{O}_{X} \rightarrow {\sf Mod }\Gamma(X,\mathcal{O}_{X})$ is the global sections functor.  Since $\operatorname{R}^{i}\omega \cong (\operatorname{R}^{i+1}\tau)\pi$ for $i \geq 1$ (Theorem \ref{theorem.exacts}), the cohomology of (\ref{eqn.gamma0}) is thus related to the cohomology of $\tau$.  By the third property of $\HU(-,-)$ above, the cohomology of (\ref{eqn.gamma0}) is thus related to the cohomology of $\underset{n \to \infty}{\lim} \HU (\A/\A_{\geq n}, -)$.  In order to compute the cohomology of $\underset{n \to \infty}{\lim} \HU (\A/\A_{\geq n}, -)$, we use the following theorem of Van den Bergh:

\begin{theorem} \label{theorem.sequence} \cite[Theorem 7.1.2]{p1bundles}
Let $\mathcal{O}_{X}$ denote the trivial right $\A$-module.  The sequence of $\mathcal{O}_{X}$-$\A$ bimodules
\begin{equation} \label{eqn.exact}
0 \rightarrow \Q \otimes e_{m+2}\A \rightarrow \E \otimes e_{m+1}\A \rightarrow e_{m}\A \rightarrow \mathcal{O}_{X} \rightarrow 0
\end{equation}
whose maps come from the structure of the relations in $\A$, is exact.
\end{theorem}
In order to use this theorem, we prove, using a variant of the first property of $\HU(-,-)$ above, that each term in (\ref{eqn.exact}) to the left of $\mathcal{O}_{X}$ is $\HU(-,\M)_{m}$-acyclic.  We may thus use (\ref{eqn.exact}) to compute the cohomology of $\HU(\mathcal{O}_{X},-)$.  In fact, we prove that $\A$ is Gorenstein (Theorem \ref{theorem.gor}):
\begin{theorem} 
Let $\mathcal{L}$ be a coherent, locally free $\mathcal{O}_{X}$-module.  Then 
$$
\ExtU^{i}(\mathcal{O}_{X},\mathcal{L} \otimes e_{l}\A)=0 \mbox{ for $i \neq 2$}
$$
and
$$
\ExtU^{2}(\mathcal{O}_{X},\mathcal{L} \otimes e_{l}\A)_{j}  \cong 
\begin{cases}
\mathcal{L} \otimes \Q_{l-2}^{*} \text{ if $j = l -2$}, \\
0& \text{otherwise}
\end{cases}
$$
where $Q_{l-2}$ is the image of the unit map $\eta:  \mathcal{O}_{X} \rightarrow \A_{l-2,l-1} \otimes \A_{l-2,l-1}^{*}$  (See Section 2.1 for the definition of $\eta$).
\end{theorem}
Using this result, we compute the cohomology of the limit 
$$
\underset{n \to \infty}{\lim} \HU (\A/\A_{\geq n}, -)
$$
which, in turn, allows us to prove the technical results regarding (\ref{eqn.gamma0}) we need to prove Serre duality.

\subsection{Notation}
Suppose we are given a diagram of categories, functors, and natural transformations between functors:
$$
\xymatrix{
\sf{X} \rtwocell^F_{F'}{\Delta} & \sf{Y} \rtwocell^G_{G'}{\Theta} & \sf{Z}.
}
$$
The {\bf horizontal composition of $\Delta$ and $\Theta$}, denoted $\Theta * \Delta$ is defined, for every object $X$ of $\sf{X}$, by the formula   

$$(\Theta * \Delta)_{X} = \Theta_{F'X} \circ G(\Delta_{X})=G'(\Delta_{X}) \circ \Theta_{FX}.$$

Suppose $\sf{A}$ is an abelian category.  If $\sf{C}$ is a localizing subcategory, we let $\pi: \sf{A} \rightarrow \sf{A}/\sf{C}$ denote the quotient functor, $\omega: \sf{A}/\sf{C} \rightarrow \sf{A}$ denote the section functor and $\tau:\sf{A} \rightarrow \sf{C}$ denote the torsion functor.  

We let ${\sf Ch(A)}$, ${\sf K(A)}$ and ${\sf D(A)}$ denote the category of cochain complexes of objects of $\sf{A}$, cochain complexes of objects of $\sf{A}$ with morphisms the chain homotopy equivalence classes of maps between complexes, and the derived category of $\sf{A}$, respectively.  We will sometimes just write $\sf{Ch}$, $\sf{K}$ and $\sf{D}$.  We will write $Q:\sf{K} \rightarrow \sf{D}$ for the localization functor.  If $F: \sf{A} \rightarrow \sf{B}$ is a left exact functor between abelian categories whose derived functor exists, we will write ${\bf R}F$ for the derived functor of $F$.

Throughout, we let $X$ denote a smooth scheme of finite type over a field $K$, and we let ${\sf Mod }X$ denote the category of quasi-coherent $\mathcal{O}_{X}$-modules.  
\vfill
\eject

%% file: ncsym.tex
\section{Non-commutative Symmetric Algebras and ${\mathbb{P}}^{1}$-bundles}

\subsection{Preliminaries}
We remind the reader of some definitions from \cite{p1bundles}.  A {\bf coherent} $\boldsymbol{\mathcal{O}_{X}}${\bf -bimodule}, $\E$, is a coherent $\mathcal{O}_{X^{2}}$-module whose support over both copies of $X$ is finite.  An $\boldsymbol{\mathcal{O}_{X}}$-{\bf bimodule} is a direct limit of coherent $\mathcal{O}_{X}$-bimodules.

For $i=1,2$ and $1 \leq j < l \leq 3$, let $pr_{i}:X^{2} \rightarrow X$ and $pr_{jl}:X^{3} \rightarrow X^{2}$ denote the standard projections.  Let $d:X \rightarrow X \times X$ be the diagonal map and let $\mathcal{O}_{\Delta}=d_{*}\mathcal{O}_{X}$.

If $\E$ and $\F$ are $\mathcal{O}_{X}$-bimodules, define the tensor product by
$$
\E \otimes_{\mathcal{O}_{X}} \F := pr_{13*}(pr_{12}^{*}\E \otimes_{\mathcal{O}_{X^{3}}} pr_{23}^{*}\F).
$$
If $\M$ is an $\mathcal{O}_{X}$-module, define $\M \otimes_{\mathcal{O}_{X}} \E$ by $pr_{2*}(pr_{1}^{*}\M \otimes \E)$.  In what follows, we will drop the subscript on the tensor product.

Let $\mu_{\mathcal{O}}:\E \otimes \mathcal{O}_{\Delta} \rightarrow \E$ and $_{\mathcal{O}}\mu:\mathcal{O}_{\Delta} \otimes \E \rightarrow \E$ denote the left and right scalar multiplication morphisms (see \cite[Proposition 3.7, p. 35]{me} for the definitions).

A coherent $\mathcal{O}_{X}$-bimodule $\E$ is said to be {\bf locally free of rank} $\boldsymbol{n}$ if $pr_{i*}\E$ is locally free of rank $n$ for $i=1,2$.

We shall use the following result without comment in the sequel.
\begin{proposition} \cite[p. 6]{p1bundles}
If $\E$ is a locally free $\mathcal{O}_{X}$-bimodule of rank $n$, there exists a locally free $\mathcal{O}_{X}$-bimodule $\E^{*}$, the dual of $\E$, of rank $n$ and natural transformations
$$
\eta_{\E}:\operatorname{id}_{{\sf Mod}X} \rightarrow (- \otimes \E)\otimes \E^{*}
$$
and
$$
\epsilon_{\E}:(- \otimes \E^{*}) \otimes \E \rightarrow \operatorname{id}_{{\sf{Mod }X}}
$$
such that $(- \otimes \E,- \otimes \E^{*}, \eta_{\E}, \epsilon_{\E})$ is an adjunction.
\end{proposition} \label{prop.adjointfunct}
Consider the composition of functors
\begin{equation} \label{eqn.init}
- \otimes \mathcal{O}_{\Delta} \rightarrow \operatorname{id}_{{\sf{Mod }X}} \overset{\eta_{\E}}{\rightarrow} (- \otimes \E) (-\otimes \E^{*}) \rightarrow - \otimes (\E \otimes \E^{*})
\end{equation}
whose left arrow is induced by scalar multiplication and whose right arrow is induced by the associativity isomorphism (\cite[Propostion 2.5, p. 442]{translations}).  By \cite[Lemma 3.1.1, p.4]{p1bundles}, this morphism of functors corresponds to a morphism $\mathcal{O}_{\Delta} \rightarrow \E \otimes \E^{*}$.  We will often abuse notation by referring to this morphism as $\eta$.  Similarly, there is a morphism $\epsilon: \E^{*} \otimes \E \rightarrow \mathcal{O}_{\Delta}$ induced by $\epsilon_{\E}$.  the maps $\eta$ and $\epsilon$ are monomorphisms and epimorphisms, respectively \cite[Proposition 3.1.1, p. 7]{p1bundles}

Let $E=-\otimes \E$ and let $E^{*}=- \otimes \E^{*}$.  By Proposition \ref{prop.adjointfunct}, the compositions
$$
\begin{CD}
E\operatorname{id}_{{\sf Mod }X}  @>{E*\eta_{\E}}>> EE^{*}E @>{\epsilon_{\E}*E}>> \operatorname{id}_{{\sf Mod }X} E
\end{CD}
$$
and
$$
\begin{CD}
\operatorname{id}_{{\sf Mod }X} E^{*} @>{\eta_{\E}* E^{*}}>> E^{*}EE^{*} @>{E^{*}*\epsilon_{\E}}>> E^{*} \operatorname{id}_{{\sf Mod }X}
\end{CD}
$$
equal $E$ and $E^{*}$ respectively.  As a consequence, one can show that the compositions
$$
\begin{CD}
\E @>{_{\mathcal{O}}\mu^{-1}}>> \mathcal{O}_{\Delta} \otimes \E @>{\eta * \E}>> (\E \otimes \E^{*}) \otimes \E @>{\cong}>>
\end{CD}
$$
$$
\begin{CD}
\E \otimes (\E^{*} \otimes \E) @>{\E * \epsilon}>> \E \otimes \mathcal{O}_{\Delta} @>{\mu_{\mathcal{O}}}>> \E
\end{CD}
$$
and
$$
\begin{CD}
\E^{*} @>{{\mu_{\mathcal{O}}}^{-1}}>> \E^{*} \otimes \mathcal{O}_{\Delta} @>{\E^{*}*\eta}>> \E^{*} \otimes (\E \otimes \E^{*}) @>{\cong}>>
\end{CD}
$$
$$
\begin{CD}
(\E^{*} \otimes \E) \otimes \E^{*} @>{\epsilon * \E^{*}}>> \mathcal{O}_{\Delta} \otimes \E^{*} @>{_{\mathcal{O}}\mu}>> \E^{*}
\end{CD}
$$
whose central isomorphisms are associativity isomorphisms, equal $\E$ and $\E^{*}$ respectively.
\newline
\noindent{\it Remark:}  In what follows, we will abuse the notation we have introduced above.  Suppose $\E$ and $\F$ are locally free finite rank $\mathcal{O}_{X}$-bimodules.  As an example of the abuse which follows, we will write
$$
\E \overset{\eta}{\rightarrow} \E \otimes \E^{*} \otimes \E \overset{\eta}{\rightarrow} (\E \otimes \F \otimes \F^{*}) \otimes \E^{*} \otimes \E
$$
instead of
$$
\begin{CD}
\E @>{\cong}>> \mathcal{O}_{\Delta} \otimes \E @>{\eta \otimes \E}>> (\E \otimes \E^{*}) \otimes \E @>{((\E \otimes \eta) \otimes \E^{*}) \otimes \E}>> ((\E \otimes (\F \otimes \F^{*})) \otimes \E^{*}) \otimes \E.
\end{CD}
$$
We will also use, without further comment, the fact that the category of $\mathcal{O}_{X}$-bimodules together with the tensor product forms a monoidal category.  Hence, the coherence theorem for monoidal categories holds, so that any diagram whose arrows are associativity isomorphisms commutes.  As another example of the notational abuse we will be guilty of, if $\phi:\E \rightarrow \F$, we will write $\E \otimes \E \overset{\phi}{\rightarrow} \F \otimes \E$ instead of $\E \otimes \E \overset{\phi \otimes \E}{\rightarrow} \F \otimes \E$.

In addition, we will sometimes omit tensor product symbols, and will write $\mathcal{O}_{X}$ instead of $\mathcal{O}_{\Delta}$ where no confusion arises.

If $\E$ and $\F$ are locally free, finite rank $\mathcal{O}_{X}$-bimodules, there is an isomorphism $- \otimes (\E \otimes \F) \rightarrow (- \otimes \E) \otimes \F$ of functors from ${\sf Mod }X$ to ${\sf Mod }X$ (\cite[Propostion 2.5, p. 442]{translations}).  Thus, there exists a canonical choice of unit, $\eta$ and counit, $\epsilon$, making $(- \otimes (\E \otimes \F), (- \otimes \E^{*}) \otimes \F^{*}, \eta, \epsilon)$ an adjunction, and there is a canonical isomorphism $\phi:  (- \otimes \F^{*}) \otimes \E^{*} \rightarrow - \otimes (\E \otimes \F)^{*}$ (see Lemma \ref{lemma.commute} for more details).  The composition
$$
-\otimes (\F^{*} \otimes \E^{*}) \rightarrow (- \otimes \F^{*}) \otimes \E^{*} \overset{\phi}{\rightarrow} - \otimes (\E \otimes \F)^{*}
$$
whose left arrow is the associativity isomorphism, induces an isomorphism
\begin{equation} \label{eqn.canonicalisom0}
\F^{*} \otimes \E^{*} \rightarrow (\E \otimes \F)^{*}
\end{equation}
by \cite[Lemma 3.1.1, p.4]{p1bundles}.

\begin{definition} \label{definition.dual}
If $g:\C \rightarrow \D$ is a morphism of coherent, locally free $\mathcal{O}_{X}$-bimodules, the {\bf dual of }$\mathbf{g}$, $g^{*}:\D^{*} \rightarrow \C^{*}$, is the composition
$$
\D^{*} \overset{\eta}{\rightarrow} \D^{*}\otimes \C \otimes \C^{*} \overset{D^{*} \otimes g \otimes \C^{*}}{\rightarrow} D^{*} \otimes \D \otimes \C^{*} \rightarrow \C^{*}.
$$
\end{definition}

We will need the following result to prove Proposition \ref{lem.inclusion}.
\begin{lemma} \label{lemma.added}
Let
$$
\begin{CD}
\mathcal{O}_{X} @>{i}>> \Q @>{\alpha}>> \E \otimes \E^{*}
\end{CD}
$$
be a factorization of the unit $\eta:\mathcal{O}_{X} \rightarrow \E \otimes \E^{*}$, where $\Q$ is the image of $\eta$ in $\E \otimes \E^{*}$.  Let $\delta:\Q^{*}\otimes \E \rightarrow \E$ denotes the isomorphism
$$
\begin{CD}
\Q^{*}\E @>{_{\mathcal{O}}\mu^{-1}}>> \mathcal{O}\Q^{*}\E @>{i}>> \Q\Q^{*}\E @>{\eta^{-1}}>> \E.
\end{CD}
$$
Then the diagram
$$
\begin{CD}
\E @>{\eta}>> \Q\Q^{*}\E \\
@V{\eta}VV @VV{\alpha}V \\
\E \E^{*} \E @<<{\delta}< \E\E^{*}\Q^{*}\E
\end{CD}
$$
commutes.
\end{lemma}

\begin{proof}
It suffices to show the outer circuit in the diagram
$$
\begin{CD}
\mathcal{O}_{X} @>{\eta}>> \Q\Q^{*} @>{\alpha}>> \E\E^{*}\Q^{*} \\
@V{\eta}VV @VV{\eta}V @VV{_{\mathcal{O}}\mu^{-1}}V \\
\E\E^{*} @>>{\eta}> \E\E^{*}\Q^{*} @>>{i^{-1}}> \E\E^{*}\mathcal{O}\Q^{*}
\end{CD}
$$
commutes.  The left square commutes by functoriality of the tensor product.  To prove the right square commutes, it suffices to prove
$$
\begin{CD}
\Q\Q^{*} @>{i^{-1}}>> \mathcal{O}\Q^{*} \\
@V{_{\mathcal{O}}\mu^{-1}}VV @VV{\eta}V \\
\mathcal{O}\Q\Q^{*} & & \E\E^{*}\Q^{*} \\
\E\E^{*}\Q\Q^{*} @>>{i^{-1}}> \E\E^{*}\mathcal{O}\Q^{*}
\end{CD}
$$
commutes, where we have used our convention that
$$
\eta:\Q\Q^{*} \rightarrow \E\E^{*}\Q\Q^{*}
$$
denotes the composition
$$
\begin{CD}
\Q\Q^{*} @>{_{\mathcal{O}}\mu^{-1}}>> \mathcal{O}\Q\Q^{*} @>{\eta \Q \Q^{*}}>> \E\E^{*} \Q\Q^{*}.
\end{CD}
$$
Thus, it suffices to show
$$
\begin{CD}
\Q @>{i^{-1}}>> \mathcal{O} @>{\eta}>> \E\E^{*} \\
@V{_{\mathcal{O}}\mu^{-1}}VV @VV{_{\mathcal{O}}\mu^{-1}}V @VV{\mu_{\mathcal{O}}^{-1}}V \\
\mathcal{O}\Q @>>{i^{-1}}> \mathcal{O}\mathcal{O} @>>{\eta}> \E\E^{*}\mathcal{O}
\end{CD}
$$
commutes.  Since $_{\mathcal{O}}\mu:\mathcal{O} \mathcal{O} \rightarrow \mathcal{O}$ equals $\mu_{\mathcal{O}}:\mathcal{O}\mathcal{O} \rightarrow \mathcal{O}$, the diagram commutes by the functoriality of $\otimes$.
\end{proof}

\subsection{The algebra $\A$}
We now review the definition of {\it non-commutative symmetric algebra}, and we study some of its basic properties.
\begin{definition} \cite[Section 4.1]{p1bundles}
Let $\E$ be a locally free $\mathcal{O}_{X}$-bimodule.  A {\bf non-commutative symmetric algebra in standard form generated by} $\boldsymbol{\E}$ is the sheaf $\mathbb{Z}$-algebra $\underset{i,j \in \mathbb{Z}}{\oplus}\A_{ij}$ with components
\begin{itemize}
\item{}
$\A_{ii}= \mathcal{O}_{\Delta}$

\item{}
$\A_{i,i+1}=\E^{i*}$,

\item{}
$\A_{ij}= \A_{i,i+1}\otimes \cdots \otimes \A_{j-1,j}/\R_{ij}$ for $j>i+1$, where $\R_{ij} \subset \A_{i,i+1} \otimes \cdots \otimes \A_{j-1,j}$ is the $\mathcal{O}_{X}$-bimodule
$$
\overset{j-2}{\underset{k=i}{\Sigma}}\A_{i,i+1} \otimes \cdots \otimes \A_{k-1,k} \otimes \Q_{k} \otimes \A_{k+2,k+3} \otimes \cdots \otimes \A_{j-1,j},
$$
and $\Q_{i}$ is the image of the unit map $\mathcal{O}_{\Delta} \rightarrow \A_{i,i+1}\otimes \A_{i+1,i+2}$, and
\item{}
$\A_{ij}= 0$ if $i>j$
\end{itemize}
and with multiplication defined as follows: for $i<j<k$,
\begin{align*}
\A_{ij} \otimes \A_{jk} & = \frac{\A_{i,i+1} \otimes \cdots \otimes \A_{j-1,j}}{\R_{ij}} \otimes \frac{\A_{j,j+1} \otimes \cdots \otimes \A_{k-1,k}}{\R_{jk}} \\
& \cong \frac{\A_{i,i+1} \otimes \cdots \otimes \A_{k-1,k}}{\R_{ij} \otimes \A_{j,j+1} \otimes \cdots \otimes \A_{k-1,k}+ \A_{i,i+1} \otimes \cdots \otimes \A_{j-1,j} \otimes \R_{jk}}
\end{align*}
by \cite[Corollary 3.18, p.38]{me}.  On the other hand,
$$
\R_{ik} \cong \R_{ij} \otimes \A_{j,j+1} \otimes \cdots \otimes \A_{k-1,k}+\A_{i,i+1} \otimes \cdots \otimes \A_{j-1,j} \otimes \R_{jk}+
$$
$$
\A_{i,i+1} \otimes \cdots \otimes \A_{j-2,j-1} \otimes \Q_{j-1} \otimes \A_{j+1,j+2} \otimes \cdots \otimes \A_{k-1,k}.
$$
Thus there is an epi $\mu_{ijk}:\A_{ij} \otimes \A_{jk} \rightarrow \A_{ik}$.

If $i=j$, let $\mu_{ijk}:\A_{ii} \otimes \A_{ik} \rightarrow \A_{ik}$ be the scalar multiplication map $_{\mathcal{O}}\mu:\mathcal{O}_{\Delta} \otimes \A_{ik} \rightarrow \A_{ik}$.  Similarly, if $j=k$, let $\mu_{ijk}:\A_{ij} \otimes \A_{jj} \rightarrow \A_{ij}$ be the scalar multiplication map $\mu_{\mathcal{O}}$.  Using the fact that the tensor product of bimodules is associative, one can check that multiplication is associative.
\end{definition}
We define $e_{k}\A_{\geq k+n}$ to be the sum of $\mathcal{O}_{X^{2}}$-modules $\underset{i}{\oplus}e_{k}\A_{k+n+i}$ and $\A_{\geq n}= \underset{k}{\oplus}(e_{k}\A)_{\geq k+n}$.  We define $\A_{\leq n}$ similarly.  If $\A_{ij}=0$ for $i>j$, we write $\A_{n}$ instead of $\A_{\leq n}$.

{\it Remark:}  Instead of writing $\mu_{ijk}:\A_{ij} \otimes \A_{jk} \rightarrow \A_{ik}$, we will write $\mu$.  Furthermore, if $\M$ is a right $\A$-module, we will sometimes write $\mu$ for every component of its multiplication.

We will use the fact that $\A_{ij}$ is locally free of rank $j-i+1$ \cite[Theorem 7.1.2]{p1bundles} without further comment.

The proof of the following result is similar to the proof of \cite[Proposition 7, p.18]{bourb}, so we omit it.
\begin{lemma} \label{lemma.intersect}
Let $\B$, $\B'$, $\C$ and $\C'$ be $\mathcal{O}_{X}$-bimodules such that $\B \subset \B'$ and $\C \subset \C'$.  Suppose, further, that either $\B/\B'$ is locally free or $\C/\C'$ is locally free.  If $(\B \otimes \C') \cap (\B' \otimes \C)$ denotes the appropriate pullback, the natural morphism
$$
\B \otimes \C \rightarrow (\B \otimes \C') \cap (\B' \otimes \C)
$$
is an isomorphism.
\end{lemma}
We will need the following lemma to prove Proposition \ref{lem.inclusion} and Theorem \ref{theorem.gor}.
\begin{lemma} \label{lemma.monomult}
If $\alpha:\Q_{k} \rightarrow \A_{k,k+1} \otimes \A_{k+1,k+2}$ denotes the inclusion map, the composition
$$
\A_{lk}\otimes \Q_{k} \overset{\alpha}{\rightarrow} \A_{lk} \otimes \A_{k,k+1} \otimes \A_{k+1,k+2} \overset{\mu \otimes \A_{k+1,k+2}}{\rightarrow} \A_{l,k+1} \otimes \A_{k+1,k+2}
$$
is monic.
\end{lemma}

\begin{proof}
The proof is almost identical to the proof of \cite[Theorem 7.1.2(2)]{p1bundles} and is omitted.
\end{proof}
We now present a definition and Lemma which will be useful in the proof of Lemma \ref{lemma.tor2}. Recall that the objects of ${\sf{Gr }\A}$ are sums $\underset{i \in \mathbb{Z}}{\oplus} \M_{i}$ of $\mathcal{O}_{X}$-modules with a graded right $\A$-module structure.
\begin{definition}
Let $\M$ be an object of $\sf{Gr }\A$.  The {\bf submodule of $\mathbf{\M}$ generated by $\mathbf{\M_{i}}$}, $\M^{i}$, is the graded $\A$-module with $\M^{i}_{k}=\operatorname{im}\mu_{ik}$ and with multiplication defined as follows.  On the diagram
$$
\begin{CD}
\operatorname{im }\mu_{ik} \otimes \A_{kj} @<<< \M_{i} \otimes \A_{ik} \otimes \A_{kj} @>>> \M_{i} \otimes \A_{ij} \\
@VVV @VVV @VVV \\
\M_{k}\otimes \A_{kj} @>>{=}> \M_{k} \otimes \A_{kj} @>>> \M_{j} @>>> \operatorname{cok }\mu_{ij}
\end{CD}
$$
whose left vertical is induced by inclusion, whose bottom right-most horizontal is the cokernel of $\mu_{ij}$ and whose other unlabeled maps are multiplication maps, consider the path starting at the top center and continuing left.  Since the right square commutes, this path is $0$.  Since the upper left horizontal is an epimorphism, the path starting at the upper left is $0$.  By the universal property of the cokernel, there is a morphism
$$
\M_{k}^{i} \otimes \A_{kj} \overset{=}{\rightarrow} \operatorname{im }\mu_{ik} \otimes \A_{kj} \rightarrow \operatorname{im }\mu_{ij}
$$
denoted $\mu_{jk}^{i}$.  The collection $\{\mu_{kj}^{i}\}_{k \leq j}$ gives $\M^{i}$ a right $\A$-module structure.
\end{definition}

\begin{lemma} \label{lemma.torneeds}
Let $\M$ be an object of $\sf{Gr }\A$.  For $i \leq j$ and for $k \in \mathbb{Z}$, Let $(\phi_{ij})_{k}:\operatorname{im }\mu_{ik} \rightarrow \operatorname{im }\mu_{jk}$ be the canonical morphism induced by the associativity diagram
$$
\begin{CD}
\M_{i} \otimes \A_{ij} \otimes \A_{jk} @>>> \M_{i} \otimes \A_{ik} \\
@VVV @VVV \\
\M_{j} \otimes \A_{jk} @>>> \M_{k}
\end{CD}
$$
whose arrows are multiplication.  Then $\phi_{ij}:\M^{i} \rightarrow \M^{j}$ is a morphism of right $\A$-modules making $\{\M^{i},\phi_{ij}\}$ a direct system.  Furthermore, the direct limit of this system is $\M$.
\end{lemma}

\subsection{Generators for ${\sf D}({\sf Proj }\A)$}
The purpose of this section is to construct a set of generators for the derived category of a non-commutative $\mathbb{P}^{1}$-bundle ${\sf Proj}\A$, $\sf{D}({\sf Proj}\A)$.  We will later prove (Theorem \ref{theorem.cd}, Proposition \ref{prop.main}) these generators are compact.

In order to construct a set of generators for the category ${\sf D}({\sf Proj }\A)$, we must first study generators of the category ${\sf{Gr }\A}$.

\begin{definition}
A {\bf locally noetherian category} is a Grothendieck category with a set of noetherian generators.
\end{definition}

\begin{lemma} \label{lem.fgnoeth} \cite[Exercise 5.4, p.56]{smith2}
In a locally noetherian category, every finitely generated object is noetherian.
\end{lemma}

\begin{lemma} \label{lem.noetherian}
If $\sf{C}$ is locally noetherian, any object $\M$ of $\sf{C}$ is the direct limit of its noetherian subobjects.
\end{lemma}
We remind the reader that an object in ${\sf{Gr }\A}$ is {\bf torsion} if it is a direct limit of right-bounded objects.

\begin{lemma} \label{lemma.smithuse}
If ${\sf Gr} \A$ is locally noetherian, an essential extension of a torsion module in $\sf{Gr }\A$ is torsion.
\end{lemma}

\begin{proof}
Suppose $\M$ is torsion, and let $\F$ be an essential extension of $\M$.  Since $\sf{Gr }\A$ is locally noetherian, $\F$ is a direct limit of its noetherian subobjects by Lemma \ref{lem.noetherian}.  Let $\N$ be a noetherian subobject of $\F$.  Then $\N \cap \M$ is torsion, hence right-bounded since $N$ is noetherian.  Thus $\N_{\geq n} \cap \M = 0$ for $n>>0$, whence $\N$ is right-bounded.
\end{proof}

\begin{lemma} \label{lem.not}
If ${\sf Gr}\A$ is locally noetherian and $\M$ is an object in ${\sf Gr }\A$, $\M$ is torsion if and only if every noetherian subobject is right bounded.
\end{lemma}

\begin{proof}
Since ${\sf Gr }\A$ is locally noetherian, any object $\M$ in ${\sf Gr }\A$ is the direct limit of its noetherian subobjects by Lemma \ref{lem.noetherian}.  Thus, if every noetherian subobject is right bounded, $\M$ is torsion.

Conversely, suppose $\M$ is torsion.  Since ${\sf Tors }\A$ is closed under subquotients, every noetherian subobject $\N$ of $\M$ is torsion.  Since $\N$ is torsion, it is a direct limit of its right bounded submodules.  On the other hand, since $\N$ is noetherian, it is a direct limit of finitely many of its right bounded submodules.  Thus, $\N$ is right bounded.
\end{proof}

\begin{corollary} \label{cor.not}
If ${\sf Gr }\A$ is locally noetherian, $\M$ is an object of ${\sf Gr }\A$ which is not torsion and $M \in \mathbb{Z}$ is given, there exists an $m \geq M$ and a noetherian subobject $\N$ of $\M$ generated in degree $m$ which is not torsion.
\end{corollary}

\begin{proof}
Since $\M$ is not torsion, Lemma \ref{lem.not} implies that $\M$ has a noetherian subobject $\N$ which is not right bounded.  If, for all $m \geq M$, $\N_{\geq m}$ were torsion, then since $\N_{\geq m}$ is noetherian, $\N_{\geq m}$ would be right bounded.  This contradicts the fact that $\N$ is not right bounded.
\end{proof}
The following lemma is a variant of \cite[Lemma 3.4, p. 450]{translations}.
\begin{lemma} \label{lem.change}
Let $\M$ be an $\mathcal{O}_{X}$-module, let $\N$ be a graded $\A$-module and let $u:\A_{kk} \rightarrow e_{k}\A$ be the inclusion morphism.  Then there is an isomorphism
$$
\Homa(\M \otimes e_{k}\A,\N) \rightarrow \Homo(\M, \N_{k})
$$
given by sending $f \in \Homa(\M \otimes e_{k}\A,\N)$ to the composition
$$
\M \overset{\cong}{\rightarrow} \M \otimes \A_{kk} \overset{\operatorname{id} \otimes u}{\rightarrow} \M \otimes e_{k}\A \overset{f}{\rightarrow} \N
$$
The inverse is given by sending $g \in \Homo(\M, \N_{k})$ to the composition
$$
\M \otimes e_{k}\A \overset{g \otimes e_{k}\A}{\rightarrow} {\N}_{k}\otimes e_{k}\A \overset{\mu}{\rightarrow}  \N.
$$
\end{lemma}

\begin{lemma} \label{lem.fg}
Let $\N$ be an object of ${\sf Mod }X$.  If $\N$ is noetherian, $\N \otimes e_{k}\A$ is a finitely generated graded $\A$-module.
\end{lemma}

\begin{proof}
Suppose there is an epimorphism
$$
{\bigoplus}_{i \in I} \M_{i} \rightarrow \N\otimes e_{k}\A
$$
of graded $\A$-modules with $I$ countable (the notation $\M_{i}$ here is the $i$th $\A$-module, not the $i$th graded piece of a graded $\A$-module $\M$).  If we identify $\N$ with $\N \otimes \A_{kk}$ and $\M_{i}$ with its image in $\N \otimes e_{k}\A$, we find that
$$
{\bigcup}_{i \in I}\M_{i} \cap \N=\N.
$$
Since $\N$ is noetherian, there exists a finite set $I' \subset I$ such that $\bigcup_{I'}\M_{i} \cap \N=\N$.  Thus, the restriction $\bigoplus_{i \in I'}\M_{i} \rightarrow \N\otimes e_{k}\A$ is an epimorphism of graded $\A$-modules, as desired.
\end{proof}

\begin{lemma} \label{lem.noethno}
Let $X$ be a scheme such that ${\sf Mod }X$ is locally noetherian with noetherian generators $\{\M_{i}\}_{i \in I}$.  If ${\sf Gr}\A$ is locally noetherian, $\{\M_{i}\otimes e_{k}\A\}_{i \in I, k \in {\mathbb{Z}}}$ is a set of noetherian generators of ${\sf Gr}\A$.
\end{lemma}

\begin{proof}
By Lemmas \ref{lem.fg} and Lemma \ref{lem.fgnoeth}, $\M_{i} \otimes e_{k}\A$ is noetherian.  By Lemma \ref{lem.change},
\begin{equation} \label{eqn.gens}
\Homo(\M_{i}, \N_{k}) \cong \Homa(\M_{i}\otimes e_{k}\A,\N).
\end{equation}
Thus, if $\N \neq 0$, then the left hand side of (\ref{eqn.gens}) is nonzero for some $i \in I,k \in \mathbb{Z}$ since $\{ \M_{i} \}_{i \in I}$ generates ${\sf Mod }X$.  Thus, the right hand side of (\ref{eqn.gens}) is nonzero for some $i \in I, k \in \mathbb{Z}$.
\end{proof}

\begin{definition}
A category $\sf{B}$ {\bf satisfies (Gen) with a set of objects} $\mathbf{\{\B_{i,m}\}_{i, m \in \mathbb{Z}}}$ (in $\sf{B}$) if, for any ${\C}^{*}$ in $D(\sf{B})$ with $h^{n}{\C}^{*} \neq 0$ then
$$
\mbox{Hom}_{{\sf D}(\sf{B})}(\B_{i,m}\{-n\},{\C}^{*}) \neq 0
$$
for $i,m>>0$, where $\{-\}$ is the suspension functor in $\sf{D}(\sf B)$.
\end{definition}
The following result is a variation of \cite[Lemma 2.2, p.712]{duality}.
\begin{proposition} \label{prop.gen}
If ${\sf Gr}\A$ is locally noetherian then ${\sf Proj }\A$ satisfies $\operatorname{(Gen)}$ with the collection $\{\pi(\mathcal{O}_{X}(i)\otimes e_{m}\A)\{n\}\}_{i, m,n \in {\mathbb{Z}}}$.
\end{proposition}

\begin{proof}
Let ${\C}^{*} \in {\sf D}({\sf Proj }\A)$ have $\operatorname{h}^{n}{\C}^{*} \neq 0$.  Given $M, I \in {\mathbb{Z}}$, we claim there exist $m \geq M$ and $i \geq I$ such that
$$
\operatorname{Hom}_{{\sf D}({\sf Proj }\A)}(\mathcal{O}_{X}(-i) \otimes e_{m}\A\{-n\},{\C}^{*}) \neq 0.
$$
To prove the claim, we first note that ${\C}^{*}$ is a complex over ${\sf Proj }\A$ so $\D^{*}=\omega{\C}^{*}$ is a complex over ${\sf Gr }\A$ and ${\C}^{*} \cong \pi \omega {\C}^{*}= \pi \D^{*}$.  Since $\pi$ is exact and $\operatorname{h}^{n}{\C}^{*} \neq 0$, $\operatorname{h}^{n}\D^{*}$ is not torsion.

By Corollary \ref{cor.not} there exists an $m \geq M$ and a noetherian submodule $\N$ of $h^{n}\D^{*}$ generated in degree $m$ which is not torsion.  Let $\psi$ be the composition
$$
(Z^{n}\D^{*})_{m} \overset{\delta}{\rightarrow} (h^{n}\D^{*})_{m} \rightarrow (h^{n}\D^{*}/\N)_{m}.
$$
whose left arrow, $\delta$, is the quotient map.  The image of the kernel of $\psi$ in $(h^{n}\D^{*})_{m}$ is $\N_{m}$.  Since $X$ is noetherian, $\operatorname{ker} \psi$ is the direct limit of its coherent submodules, and the image of some such submodule, $\M$, under $\delta$ must generate a non-torsion $\A$-submodule $\N'$ of $\N$.  By \cite[Corollary 5.18, p. 121]{hartshorne} there exists an integer $i_{0}$ such that for $i \geq i_{0}$ there is an epimorphism
$$
\underset{k}{\oplus}\mathcal{O}_{X}(-i) \rightarrow \M \rightarrow {\N'}_{m}.
$$
The image of some summand must generate a non-torsion $\A$-submodule $\N^{''}$ of $\N^{'}$.  Thus, there is a nonzero morphism
\begin{equation} \label{eqn.o1}
{\mathcal{O}_{X}}(-i) \rightarrow (\operatorname{Z}^{n}\D^{*})_{m} \rightarrow (\operatorname{h}^{n}\D^{*})_{m}
\end{equation}
whose image equals $\N^{''}_{m}$.  Tensoring (\ref{eqn.o1}) by $e_{m}\A$ gives a map of complexes
$$
\mathcal{O}_{X}(-i) \otimes e_{m}\A\{-n\} \rightarrow \D^{*}
$$
whose induced map in cohomology has non-torsion image.
\end{proof}

\vfill
\eject 

%% file: hom.tex
\section{Internal Tensor and Hom functors on ${\sf Gr}\A$} 
\begin{definition} \label{def.B}
Let $\sf{Bimod }\A-\A$ denote the category of $\A-\A$-bimodules.  Specifically:
\begin{itemize}

\item{}
an object of $\sf{Bimod }\A-\A$ is a triple 
$$
(\mathcal{C}=\{C_{ij}\}_{i,j \in \mathbb{Z}}, \{\mu_{ijk}\}_{i,j,k \in \mathbb{Z}}, \{\psi_{ijk}\}_{i,j,k \in \mathbb{Z}})
$$ 
where ${\mathcal{C}}_{ij}$ is an ${\mathcal{O}}_{X}$-bimodule and $\mu_{ijk}:\C_{ij} \otimes \A_{jk} \rightarrow \C_{ik}$ and $\psi_{ijk}: \A_{ij} \otimes \C_{jk} \rightarrow \C_{ik}$ are morphisms of $\mathcal{O}_{X^{2}}$-modules making $\C$ an $\A$-$\A$ bimodule.

\item{}
A morphism $\phi:  \mathcal{C} \rightarrow \mathcal{D}$ between objects in ${\sf Bimod }\A-\A$ is a  collection $\phi=\{\phi_{ij}\}_{i,j \in \mathbb{Z}}$ such that $\phi_{ij}:{\mathcal{C}}_{ij} \rightarrow {\mathcal{D}}_{ij}$ is a morphism of ${\mathcal{O}}_{X^{2}}$-modules, and such that $\phi$ respects the $\mathcal{A}-\mathcal{A}$-bimodule structure on $\mathcal{C}$ and $\mathcal{D}$.  
\end{itemize}

Let $\sf{Bimod }\mathcal{O}_{X}-\A$ denote the full subcategory of $\sf{Bimod }\A-\A$ consisting of objects $\mathcal{C}$ such that for some $n \in \mathbb{Z}$, ${\mathcal{C}}_{ij}=0$ for $i \neq n$ (we say $\C$ is {\bf left-concentrated in degree $\mathbf{n}$}).

Let $\mathbb{B}$ denote the full subcategory of $\sf{Bimod }\A-\A$ whose objects $\mathcal{C}=\{C_{ij}\}_{i,j \in \mathbb{Z}}$ have the property that $\C_{ij}$ is coherent and locally free for all $i,j \in \mathbb{Z}$. 

Let ${\mathbb{G}}{\sf r} \mathcal{A}$ denote the full subcategory of $\mathbb{B}$ consisting of objects $\mathcal{C}$ such that for some $n \in \mathbb{Z}$, ${\mathcal{C}}_{ij}=0$ for $i \neq n$ (we say $\C$ is left-concentrated in degree $n$).  
\end{definition}
In what follows, we usually omit indices on multiplication morphisms.

\subsection{The internal Tensor functor on ${\sf Gr }A$}
\begin{definition}
Let $\mathcal{C}$ be an object in ${\sf Bimod} \A-\A$ and let $\mathcal{M}$ be a graded right $\mathcal{A}$-module.  We define $\boldsymbol{\mathcal{M} \Ten \mathcal{C}}$ to be the $\mathbb{Z}$-graded $\mathcal{O}_{X}$-module whose $k$-th component is the coequalizer of the diagram

\begin{equation} 
\begin{CD}
\underset{l}{\oplus} \underset{m}{\oplus} ({\mathcal{M}}_{l} \otimes {\mathcal{A}}_{lm}) \otimes {\mathcal{C}}_{mk}   @>{{\mu}_{\mathcal{M}}\otimes \C}>> \underset{m}{\oplus} {\mathcal{M}}_{m}\otimes {\mathcal{C}}_{mk} \\
@V{\M \otimes {\mu}_{\mathcal{C}}}VV		@VV{=}V\\
\underset{l}{\oplus} {\mathcal{M}}_{l}\otimes {\mathcal{C}}_{lk}  @>>{=}> \underset{m}{\oplus} {\mathcal{M}}_{m}\otimes {\mathcal{C}}_{mk}.  
\end{CD}
\end{equation}

Let $\mathcal{C}$ be an object in ${\sf Bimod} \A-\A$ and let $\mathcal{D}$ be an object in ${\sf Bimod} \mathcal{O}_{X}-\A$ left-concentrated in degree $n$.  We define $\boldsymbol{\mathcal{D} \Ten \mathcal{C}}$ to be the $\mathbb{Z}$-graded $\mathcal{O}_{X^{2}}$-module whose $k$-th component is the coequalizer of the diagram

\begin{equation} 
\begin{CD}
\underset{l}{\oplus} \underset{m}{\oplus} ({\mathcal{D}}_{nl} \otimes {\mathcal{A}}_{lm}) \otimes {\mathcal{C}}_{mk}   @>{{\mu}_{\mathcal{D}}\otimes \C}>> \underset{m}{\oplus} {\mathcal{D}}_{nm}\otimes {\mathcal{C}}_{mk} \\
@V{\D \otimes {\mu}_{\mathcal{C}}}VV		@VV{=}V\\
\underset{l}{\oplus} {\mathcal{D}}_{nl}\otimes {\mathcal{C}}_{lk}  @>>{=}> \underset{m}{\oplus} {\mathcal{D}}_{nm}\otimes {\mathcal{C}}_{mk}.  
\end{CD}
\end{equation}
Since each component of $\mathcal{D} \Ten \mathcal{C}$ is a quotient of a direct limit of $\mathcal{O}_{X}$-bimodules, each component is an $\mathcal{O}_{X}$-bimodule.
\end{definition}

\begin{proposition}
If $\M$ is an object in $\sf{Gr }\A$, $\C$ is an object in ${\sf Bimod }\A-\A$ and $\D$ is an object in $\sf{Bimod }\mathcal{O}_{X}-\A$, then ${\mathcal{M}}\Ten \mathcal{C}$ and ${\mathcal{D}}\Ten \mathcal{C}$ inherit a graded right $\mathcal{A}$-module structure from the right $\mathcal{A}$-module structure of $\mathcal{C}$, making $-\Ten \mathcal{C}:{\sf{Gr }}\mathcal{A} \rightarrow {\sf{Gr }}\mathcal{A}$ and $-\Ten \mathcal{C}:{\sf{Bimod }}\mathcal{O}_{X}-\mathcal{A} \rightarrow {\sf{Gr }}\mathcal{A}$ functors.
\end{proposition}

\begin{proof}
We only prove the first claim.  We construct a map 
\begin{equation} \label{eqn.c}
\mu_{ij}: (\M \Ten \C)_{k} \otimes \A_{kj} \rightarrow (\M \Ten \C)_{j},
\end{equation}
i.e., we construct a map from the coequalizer of
\begin{equation} \label{eqn.a}
\begin{CD}
\underset{l}{\oplus}\underset{m}{\oplus}((\M_{l} \otimes \A_{lm}) \otimes \C_{mk}) \otimes \A_{kj} @>{\mu}>> \underset{m}{\oplus} (\M_{m} \otimes \C_{mk}) \otimes \A_{kj} \\
@V{\mu}VV @VV{=}V \\
\underset{l}{\oplus}(\M_{l} \otimes \C_{lk}) \otimes \A_{kj} @>>{=}> \underset{m}{\oplus} (\M_{m} \otimes \C_{mk}) \otimes \A_{kj}
\end{CD}
\end{equation}
to the coequalizer of
\begin{equation} \label{eqn.b}
\begin{CD}
\underset{l}{\oplus}\underset{m}{\oplus}((\M_{l} \otimes \A_{lm}) \otimes \C_{mj} @>{\mu}>> \underset{m}{\oplus}\M_{m} \otimes \C_{mj} \\
@V{\mu}VV @VV{=}V \\
\underset{l}{\oplus}\M_{l}\otimes \C_{lj} @>>{=}> \underset{m}{\oplus}\M_{m} \otimes \C_{mj}.
\end{CD}
\end{equation}
We think of (\ref{eqn.a}) and (\ref{eqn.b}) as the top and bottom, respectively, of a cube whose vertical arrows are induced by multiplication maps.  In order to construct a map (\ref{eqn.c}), it suffices, by the universal property of coequalizers, to prove that the vertical faces of this cube commute.  The diagram
$$
\begin{CD}
\underset{l}{\oplus}\underset{m}{\oplus}((\M_{l} \otimes \A_{lm}) \otimes \C_{mk}) \otimes \A_{kj} @>>> \underset{m}{\oplus} (\M_{m} \otimes \C_{mk}) \otimes \A_{kj} \\
@VVV @VVV \\
\underset{l}{\oplus}\underset{m}{\oplus}(\M_{l} \otimes \A_{lm}) \otimes \C_{mj} @>>> \underset{m}{\oplus} \M_{m} \otimes \C_{mj}
\end{CD}
$$
whose arrows are induced by multiplication, commutes by functoriality of $\otimes$.  In addition, the diagram
$$
\begin{CD}
\underset{l}{\oplus}\underset{m}{\oplus}((\M_{l} \otimes \A_{lm}) \otimes \C_{mk}) \otimes \A_{kj} @>>> \underset{l}{\oplus}(\M_{l} \otimes \C_{lk}) \otimes \A_{kj} \\
@VVV @VVV \\
\underset{l}{\oplus}\underset{m}{\oplus}(\M_{l}\otimes \A_{lm}) \otimes \C_{mj} @>>> \underset{l}{\oplus}\M_{l} \otimes \C_{lj}
\end{CD}
$$
whose arrows are induced by multiplication, commutes by the associativity of $\A-\A$-bimodule multiplication.  Thus, by the universal property of coequalizers, there exists a map
$$
\mu_{ij}:(\M \Ten \C)_{k} \otimes \A_{kj} \rightarrow (\M \Ten \C)_{j}.
$$
The fact that $\mu$ is associative and compatible with scalar multiplication follows from the fact that the right $\A$ multiplication on $\C$ is associative and compatible with scalar multiplication.

The proof of the functoriality of $-\Ten \C$ and $\M \Ten -$ is straightforward, and we omit it.
\end{proof}
We now establish some important properties of $\Ten$.
\begin{lemma} \label{lemma.asso}
If $\mathcal{L}$ is an object of $\sf{Mod }X$ and $\D$ is an object in ${\mathbb{G}{\sf{r}} \A}$, there is an isomorphism of right $\A$-modules 
$$
(\mathcal{L} \otimes \D) \Ten \C \rightarrow \mathcal{L} \otimes (\D \Ten \C).
$$
natural in $\C$.
\end{lemma}

\begin{proof}
The proof follows from the associativity of the ordinary tensor product of bimodules, and we omit the details.
\end{proof}

\begin{proposition} \label{prop.multiso}
If $\M$ is an object in ${\sf Gr} \A$, multiplication induces a natural isomorphism $\M \Ten \A \rightarrow \M$ in ${\sf Gr }\A$.  If $\C$ is an object in ${\sf Bimod }\A-\A$, multiplication induces a natural isomorphism $e_{i}\A \Ten \C \rightarrow e_{i}\C$.
\end{proposition}

\begin{proof}
The proof of the second result is similar to the proof of the first result, so we omit it.  We first prove $\mu_{\M}: \underset{k}{\oplus} \underset{m}{\oplus} {\M}_{m} \otimes {\A}_{mk} \rightarrow \underset{k}{\oplus}{\M}_{k}$ is the coequalizer of
\begin{equation} \label{eqn.tensor2}
\begin{CD}
\underset{k}{\oplus}\underset{l}{\oplus}\underset{m}{\oplus} ({\M}_{l} \otimes {\A}_{lm}) \otimes {\A}_{mk} @>{\alpha=\mu_{\M}\otimes \A}>> \underset{k}{\oplus} \underset{m}{\oplus} {\M}_{m} \otimes {\A}_{mk} \\
@V{\beta=\M \otimes \mu_{\A}}VV @VV{=}V \\
\underset{k}{\oplus}\underset{l}{\oplus} {\M}_{l} \otimes {\A}_{lk} @>>{=}> \underset{k}{\oplus} \underset{m}{\oplus} {\M}_{m} \otimes {\A}_{mk} 
\end{CD}
\end{equation}
which will prove that multiplication induces an isomorphism $\M \Ten \A \rightarrow \M$.  We will then prove that this isomorphism is natural.  
\newline
\newline
\noindent {\it Part 1:  We show $\mu_{\M}: \underset{k}{\oplus} \underset{m}{\oplus} {\M}_{m} \otimes {\A}_{mk} \rightarrow \underset{k}{\oplus}{\M}_{k}$ is the coequalizer of (\ref{eqn.tensor2}).}

In order to show $\mu_{\M}: \underset{k}{\oplus} \underset{m}{\oplus} {\M}_{m} \otimes {\A}_{mk} \rightarrow \underset{k}{\oplus}{\M}_{k}$ is the coequalizer of (\ref{eqn.tensor2}), it suffices to prove that given a diagram
\begin{equation} \label{eqn.tensor4}
\begin{CD}
\underset{k}{\oplus}\underset{l}{\oplus}\underset{m}{\oplus} ({\M}_{l} \otimes {\A}_{lm}) \otimes {\A}_{mk} @>{\alpha}>> \underset{k}{\oplus} \underset{m}{\oplus} {\M}_{m} \otimes {\A}_{mk} @>{f}>> \underset{k}{\oplus} {\D}_{k} \\
@V{\beta}VV @VV{\mu_{\M}}V @VV{=}V \\
\underset{k}{\oplus}\underset{l}{\oplus} {\M}_{l} \otimes {\A}_{lk} @>>{{\mu}_{\M}}> \underset{k}{\oplus} {\M}_{k} @>>{g}> \underset{k}{\oplus} {\D}_{k}
\end{CD}
\end{equation}
of right $\A$-modules such that $f \alpha = f \beta$, there exists a unique $g$ making the right square commute.  Notice that the right square in (\ref{eqn.tensor4}) commutes by associativity of right $\A$-module multiplication.

Let $\delta_{k}$ denote the composition
$$
{\M}_{k} \overset{{\mu_{\mathcal{O}}}^{-1}}{\rightarrow} \M_{k}\otimes \A_{kk} \rightarrow \underset{m}{\oplus}\M_{m}\otimes \A_{mk}
$$
whose right arrow is inclusion, and let $\delta = \underset{k}{\oplus}\delta_{k}$.  It is easy to see that $\mu_{\M} \circ \delta ={\operatorname{id}}_{\M}$.  

Let 
$$
\gamma:\underset{k}{\oplus}\underset{l}{\oplus}\underset{m}{\oplus} (\M_{l} \otimes \A_{lm}) \otimes \A_{mk} \rightarrow \underset{k}{\oplus} \underset{m}{\oplus} \M_{m} \otimes \A_{mk}
$$
denote the map $-\beta+\alpha-\delta\mu_{\M}\beta$.
\newline
\newline
\noindent {\it Part 1, Step 1:  We prove the map $\gamma$ is an epimorphism.}
To show that $\gamma$ is an epi, it suffices to construct a map
$$
\zeta: \underset{k}{\oplus}\underset{m}{\oplus} \M_{m} \otimes \A_{mk} \rightarrow \underset{k}{\oplus}\underset{l}{\oplus}\underset{m}{\oplus} (\M_{m} \otimes \A_{ml}) \otimes \A_{lk}
$$  
such that $\gamma \zeta = -\operatorname{id}$.  Let $\zeta_{k}$ be induced by the composition
$$
\A_{mk} \overset{{\mu_{\mathcal{O}}}^{-1}}{\rightarrow} \A_{mk}\otimes \A_{kk} \rightarrow \underset{l}{\oplus}\A_{ml}\otimes \A_{lk}
$$
whose right composite is inclusion of a summand.

Now, $\gamma \zeta = (-\beta+\alpha-\delta \mu_{\M} \beta)\zeta= -\beta \zeta + \alpha \zeta - \delta \mu_{\M}\beta \zeta$.  We notice that $(\beta \zeta)_{k}$ equals the composition
$$
\underset{m}{\oplus} \M_{m} \otimes \A_{mk} \overset{{\mu_{\mathcal{O}}}^{-1}}{\rightarrow} \underset{m}{\oplus}\M_{m} \otimes (\A_{mk} \otimes \A_{kk}) \rightarrow \underset{l}{\oplus} \underset{m}{\oplus}\M_{m} \otimes (\A_{ml} \otimes \A_{lk}) \overset{\mu}{\rightarrow}
$$
$$
\underset{m}{\oplus}\M_{m} \otimes \A_{mk}
$$
whose second arrow is inclusion in a direct summand.  This composition is the identity map.  Next, we compute $\alpha \zeta$.  Since the diagram
$$
\begin{CD}
\underset{m}{\oplus} \M_{m} \otimes \A_{mk} @>{{\mu_{\mathcal{O}}}^{-1}}>> \underset{m}{\oplus}\M_{m} \otimes (\A_{mk} \otimes \A_{kk}) @>>> \underset{l}{\oplus} \underset{m}{\oplus} \M_{m} \otimes (\A_{ml} \otimes \A_{lk}) \\
@V{\mu_{\M}}VV @VV{\mu_{\M} \otimes \A_{kk}}V @VV{\alpha_{k}}V \\
\M_{k} @>>{{\mu_{\mathcal{O}}}^{-1}}> \M_{k}\otimes \A_{kk} @>>> \underset{l}{\oplus}\M_{l} \otimes \A_{lk}
\end{CD}
$$
whose top horizontals compose to give $\zeta_{k}$ and whose bottom horizontals compose to give $\delta_{k}$, commutes, $\alpha \zeta = \delta \mu_{\M}$.  Thus, $\gamma \zeta=-\beta\zeta = -\operatorname{id}$ as desired.  Thus $\gamma$ is an epi.
\newline
\newline
\noindent{\it Part 1, Step 2:  We prove}
\begin{equation} \label{eqn.tensoriso}
\begin{CD} 
\underset{k}{\oplus}\underset{l}{\oplus}\underset{m}{\oplus} ({\M}_{l} \otimes {\A}_{lm}) \otimes {\A}_{mk} @>{\beta-\alpha}>> \underset{k}{\oplus} \underset{m}{\oplus} {\M}_{m} \otimes {\A}_{mk} \\
@V{\gamma}VV @VV{=}V \\
\underset{k}{\oplus}\underset{m}{\oplus} {\M}_{m} \otimes {\A}_{mk} @>>{\delta \mu_{\M}-\operatorname{id}}> \underset{k}{\oplus} \underset{m}{\oplus} {\M}_{m} \otimes {\A}_{mk} 
\end{CD}
\end{equation}
{\it commutes.}  

Using computations from Step 1, we have
\begin{align*}
(\delta \mu_{\M}-\operatorname{id})\gamma & = \delta \mu_{\M}(-\beta+\alpha-\delta \mu_{\M} \beta)-(-\beta+\alpha-\delta\mu_{\M}\beta) \\
& = -\delta \mu_{\M} \beta + \delta \mu_{\M} \alpha - \delta \mu_{\M}\delta\mu_{\M}\beta+\beta-\alpha+\delta\mu_{\M}\beta \\
& = -\delta \mu_{\M} \beta + \delta \mu_{\M} \alpha - \delta \mu_{\M}\beta+\beta-\alpha+\delta\mu_{\M}\beta \\
& = \delta \mu_{\M} \alpha - \delta \mu_{\M}\beta+\beta - \alpha
\end{align*}
since $\mu_{\M} \delta = \operatorname{id}$.

By the commutivity of the right square in (\ref{eqn.tensor4}) $\delta \mu_{\M}\alpha = \delta \mu_{\M} \beta$, i.e. $\delta \mu_{\M}(\alpha-\beta)=0$.  This establishes the commutivity of (\ref{eqn.tensoriso}).
\newline
\newline
\noindent{\it Part 1, Step 3:  We prove that the map $g=f \delta$ makes the right square in (\ref{eqn.tensor4}) commute.} 

We show that $g \mu_{\M}=f$, i.e. $f \delta \mu_{\M} = f$, or $f (\delta \mu_{\M} - {\operatorname{id}})=0$.  The fact that $f (\delta \mu_{\M} - {\operatorname{id}})=0$ follows from the commutivity of (\ref{eqn.tensoriso}) and the fact that $\gamma$ is an epi.  These results were established in Step 2 and Step 1, respectively.
\newline
\newline
\noindent{\it Part 1, Step 4:  We prove $g=f \delta:\M \rightarrow \D$ is unique such that $g \mu_{\M}=f$.}

Suppose $g' \mu_{\M}=f$.  Then $g' = g' \mu_{\M} \delta= f \delta = g$.  

We may conclude from Part 1 that multiplication $\mu_{\M}: \underset{k}{\oplus} \underset{m}{\oplus} {\M}_{m} \otimes {\A}_{mk} \rightarrow \underset{k}{\oplus}{\M}_{k}$ induces an isomorphism $\M \Ten \A \rightarrow \M$
\newline
\newline
\noindent{\it Part 2:  We show that the isomorphism $\M \Ten \A \rightarrow \M$ constructed in Part 1 is natural in $\M$.}  

Let $\phi:\M \rightarrow \mathcal{N}$ be a morphism in $\sf{Gr} \A$ and consider the induced diagram
$$
\begin{CD}
\underset{k}{\oplus}\underset{m}{\oplus}\M_{m} \otimes \A_{mk} @>{\underset{k}{\oplus}\underset{m}{\oplus}\phi_{m} \otimes \A_{mk}}>> \underset{k}{\oplus}\underset{m}{\oplus}\mathcal{N}_{m} \otimes \A_{mk} \\
@V{\psi}VV @VVV  \\
\M \Ten \A @>>> \N \Ten \A \\
@VVV @VVV \\
\M @>>{\phi}> \N  
\end{CD}
$$
whose bottom verticals are induced by multiplication, whose middle horizontal is induced by the top horizontal, and whose top verticals are the cokernels of the diagram defining $\Ten$.  Notice that the composition of the left verticals is $\mu_{\M}$ and the composition of the right verticals is $\mu_{\N}$.  We must show that the bottom square of this diagram commutes.  We know the outer square commutes since multiplication is natural, and we know that the top square commutes since the middle horizontal is induced by the top horizontal.  This implies that
$$
\begin{CD}
\underset{k}{\oplus}\underset{m}{\oplus}\M_{m} \otimes \A_{mk} @>{\psi}>> \M \Ten \A @>>> \M \\
& & @VVV @VVV \\
& & \N \Ten \A @>>> \N
\end{CD}
$$
commutes.  Thus, the bottom square of the first diagram commutes since $\psi$ is an epi.
\end{proof}

\begin{proposition} \label{lem.inclusion}
Let $\mathcal{L}$ be a coherent, locally free $\mathcal{O}_{X}$-module and let $\N$ be an object of $\sf{Gr }\A$.  The inclusion $\A_{\geq 1} \rightarrow \A$ induces an epi
$$
\Homa(\N \Ten \A, \mathcal{L}\otimes e_{l}\A) \rightarrow \Homa(\N \Ten \A_{\geq 1},\mathcal{L} \otimes e_{l}\A).
$$
\end{proposition}

\begin{proof}
Let $\phi: \N \Ten \A_{\geq 1} \rightarrow \mathcal{L} \otimes e_{l}\A$ be given.  It suffices, in light of Proposition \ref{prop.multiso}, to construct an $\A$-module morphism $\psi:\N \rightarrow \mathcal{L} \otimes e_{l}\A$ such that 
\begin{equation} \label{eqn.star1}
\begin{CD}
\N \Ten \A_{\geq 1} @>{\phi}>> \mathcal{L} \otimes e_{l}\A \\
@A{=}AA @AA{\psi}A \\
\N \Ten \A_{\geq 1} @>>{\mu}> \N
\end{CD}
\end{equation}
commutes, where we have abused notation by letting $\mu:\N \Ten \A_{\geq 1} \rightarrow \N$ denote the morphism {\it induced} by multiplication
$$
\underset{k}{\oplus} \underset{m<k}{\oplus} \N_{m} \otimes \A_{mk} \rightarrow \underset{k}{\oplus}\N_{k}.
$$
By the definition of the graded tensor product, the map $\phi_{k}$ corresponds to a map
$$
\gamma_{k}:\underset{m<k}{\oplus} \N_{m} \otimes \A_{mk} \rightarrow \mathcal{L} \otimes \A_{lk}
$$
such that
\begin{equation} \label{eqn.star2}
\begin{CD}
\underset{m<k}{\oplus} (\N_{l} \otimes \A_{lm}) \otimes \A_{mk} @>{\mu_{\N} \otimes \A}>> \underset{m<k}{\oplus}\N_{m} \otimes \A_{mk} \\
@V{\N \otimes \mu_{\A}}VV @VV{=}V \\
\underset{l<k}{\oplus}\N_{l}\otimes \A_{lk} @>>{=}> \underset{m<k}{\oplus} \N_{m} \otimes \A_{mk} @>>{\gamma_{k}}> \mathcal{L} \otimes \A_{l,k}
\end{CD}
\end{equation}
commutes.
\newline
\newline
\noindent{\it Step 1:  We construct the components of $\psi$.}  To define the $k$th graded component of $\psi$, notice that the diagram
$$
\begin{CD}
\N_{k} \otimes \Q_{k} @>>> \N_{k} \otimes (\A_{k,k+1} \otimes \A_{k+1,k+2}) @>>> \N_{k} \otimes \A_{k,k+2} \\
& & @V{\gamma_{k+1} \otimes \A_{k+1,k+2}}VV @VV{\gamma_{k+2}}V \\
& & \mathcal{L} \otimes \A_{l,k+1}\otimes \A_{k+1,k+2} @>>{\mu}> \mathcal{L} \otimes \A_{l,k+2}
\end{CD}
$$
whose right horizontals are multiplication maps and whose left horizontal is induced by the inclusion $\Q_{k} \rightarrow \A_{k,k+1}\otimes \A_{k+1,k+2}$ commutes, and hence factors as
\begin{equation} \label{eqn.factor1}
\begin{CD}
\N_{k} \otimes \Q_{k} @>>> \N_{k} \otimes \A_{k,k+1} \otimes \A_{k+1,k+2} @>>> \N_{k} \otimes \A_{k,k+2} \\
@VVV @V{\gamma_{k+1} \otimes \A_{k+1,k+2}}VV @VV{\gamma_{k+2}}V \\
\operatorname{ker} \mu @>>> \mathcal{L} \otimes \A_{l,k+1}\otimes \A_{k+1,k+2} @>>{\mu}> \mathcal{L} \otimes \A_{l,k+2}.
\end{CD}
\end{equation}
Let 
$$
\mathcal{R} = \underset{i=l}{\overset{k-1}{\Sigma}} \A_{l,l+1} \otimes \cdots \otimes \A_{i-1,i} \otimes \Q_{i} \otimes \A_{i+2,i+3} \otimes \cdots \otimes \A_{k,k+1} \otimes \A_{k+1,k+2}.
$$
By \cite[Corollary 3.18, p.38]{me}, there is a unique isomorphism 
$$
\zeta: \mbox{ker }\mu \rightarrow \frac{\mathcal{L} \otimes \mathcal{R} + \mathcal{L} \otimes \A_{l,l+1} \otimes \cdots \otimes \A_{k-1,k} \otimes \Q_{k}}{\mathcal{L} \otimes \mathcal{R}}
$$
making the diagram
$$
\begin{CD}
\operatorname{ker} \mu @>>> \mathcal{L} \otimes \A_{l,k+1} \otimes \A_{k+1,k+2} \\
@V{\zeta}VV @VVV \\
\frac{\mathcal{L} \otimes \mathcal{R} + \mathcal{L} \otimes \A_{l,l+1} \otimes \cdots \otimes \A_{k-1,k} \otimes \Q_{k}}{\mathcal{L} \otimes \mathcal{R}} @>>> \frac{\mathcal{L} \otimes \A_{l,l+1} \otimes \cdots \otimes \A_{k,k+1} \otimes \A_{k+1,k+2}}{\mathcal{L} \otimes \mathcal{R}}
\end{CD}
$$
whose right vertical is the canonical isomorphism from \cite[Corollary 3.18, p.38]{me} and whose horizontals are inclusions, commute.

Letting $\mathcal{R} \cap \A_{l,l+1} \otimes \cdots \otimes \A_{k-1,k} \otimes \Q_{k}$ denote the appropriate pullback, the canonical morphism 
$$
\underset{i=l}{\overset{k-2}{\Sigma}} \A_{l,l+1} \otimes \cdots \otimes \A_{i-1,i} \otimes \Q_{i} \otimes \A_{i+2,i+3} \otimes \cdots \otimes \A_{k-1,k} \otimes \Q_{k} \rightarrow \mathcal{R} \cap \A_{l,l+1} \otimes \cdots \otimes \A_{k-1,k} \otimes \Q_{k}
$$
is an isomorphism by Lemma \ref{lemma.intersect}.  Hence, there are isomorphisms
\begin{align} \label{eqn.secondisom} 
\frac{\mathcal{R} + \A_{l,l+1} \otimes \cdots \otimes \A_{k-1,k} \otimes \Q_{k}}{\mathcal{R}} & \overset{\cong}{\rightarrow} \frac{\A_{l,l+1} \otimes \cdots \otimes \A_{k-1,k} \otimes \Q_{k}}{\mathcal{R} \cap \A_{l,l+1} \otimes \cdots \otimes \A_{k-1,k} \otimes \Q_{k}} \\
& \overset{\cong}{\rightarrow} \A_{lk} \otimes \Q_{k} \notag.
\end{align}
Thus, the morphism
\begin{equation} \label{eqn.cross}
\N_{k} \otimes \Q_{k} \rightarrow \operatorname{ker }\mu \overset{\zeta}{\rightarrow} \frac{\mathcal{L} \otimes \mathcal{R} + \mathcal{L} \otimes \A_{l,l+1} \otimes \cdots \otimes \A_{k-1,k} \otimes \Q_{k}}{\mathcal{L} \otimes \mathcal{R}} \rightarrow  \mathcal{L} \otimes \A_{lk} \otimes \Q_{k}
\end{equation}
whose rightmost arrow is $\mathcal{L}$ tensored with (\ref{eqn.secondisom}), gives a morphism 
$$
\N_{k} \otimes \Q_{k} \rightarrow \mathcal{L} \otimes \A_{l,k} \otimes \Q_{k}.
$$  
Tensoring this morphism with $\Q_{k}^{*}$ gives a morphism
$$
\N_{k} \otimes \Q_{k} \otimes \Q_{k}^{*} \rightarrow \mathcal{L} \otimes \A_{lk} \otimes \Q_{k} \otimes \Q_{k}^{*}.
$$
Since $\Q_{k}$ is invertible, we thus get a morphism
$$
\psi_{k}:\N_{k} \rightarrow \mathcal{L} \otimes \A_{l,k}.
$$
\newline
\newline
\noindent{\it Step 2:  Let $\Q = \Q_{k}$, $\E = \A_{k,k+1}$ and let $\alpha: \Q \rightarrow \E \otimes \E^{*}$ denote inclusion.  We show the diagram}
$$
\begin{CD}
\N_{k} \otimes \Q @>{\alpha}>> \N_{k} \otimes \E \otimes \E^{*} @>{=}>> \N_{k} \otimes \E \otimes \E^{*} \\
@V{\psi_{k} \otimes \Q}VV & & @VV{\gamma_{k+1} \otimes \E^{*}}V \\
\mathcal{L} \otimes \A_{lk} \otimes \Q @>>{\alpha}> \mathcal{L} \otimes \A_{lk} \otimes \E \otimes \E^{*} @>>{\mu}> \mathcal{L} \otimes \A_{l,k+1} \otimes \E^{*}
\end{CD}
$$
{\it commutes.}  Consider the diagram
$$
\begin{CD}
\N_{k} \Q @>{\alpha}>> \N_{k} \E \E^{*} \\
@VVV @VV{\gamma_{k+1}}V \\
\operatorname{ker }\mu @>>> \mathcal{L}\A_{l,k+1}\E^{*} \\
@V{\zeta}VV @VV{\cong}V \\
\frac{\mathcal{L} \mathcal{R}+\mathcal{L} \A_{l,l+1}\cdots \A_{k-1,k}\Q}{\mathcal{L}\mathcal{R}} @>>> \frac{\mathcal{L}\A_{l,l+1}\cdots \E^{*}}{\mathcal{L}\mathcal{R}} @<{\cong}<< \mathcal{L}\A_{l,k+1}\E^{*} \\
@VVV & & @AA{\mu}A \\
\frac{\mathcal{L} \A_{l,l+1} \cdots \A_{k-1,k}\Q}{\mathcal{R} \cap \A_{l,l+1} \cdots \A_{k-1,k} \Q} @>>> \mathcal{L}\A_{lk} \Q @>>{\alpha}> \mathcal{L} \A_{lk}\E \E^{*}
\end{CD}
$$
whose two middle two horizontals are induced by inclusion, whose upper left vertical is the left vertical in (\ref{eqn.factor1}), whose lower left vertical and horizontal are induced by (\ref{eqn.secondisom}) and whose other unlabeled isomorphisms are both the canonical isomorphisms.  Since the left column of this diagram composed with the bottom left horizontal in the diagram equals $\psi_{k}$ tensored with $\Q_{k}$, it suffices to show the outer circuit of the diagram commutes.  Since the top two squares of the diagram commute, it suffices to show the bottom circuit commutes.  This follows from the commutivity of the diagram
\begin{equation} \label{eqn.associa}
\begin{CD}
\A_{l,l+1} \cdots \A_{k-1,k}\Q @>>> \mathcal{R}+\A_{l,l+1}\cdots \A_{k-1,k}\Q \\
@V{\mu}VV @VVV \\
\A_{l,k} & & \A_{l,l+1} \cdots \E^{*} \\
@V{\alpha}VV @VV{\mu}V \\
\A_{l,k}\E\E^{*} @>>{\mu}> \A_{l,k+1}\E^{*}
\end{CD}
\end{equation}
whose unlabeled arrows are canonical inclusions.  As the reader can check, the commutivity of (\ref{eqn.associa}) follows from the associativity of multiplication.
\newline
\newline
\noindent{\it Step 3:  Keep the notation as in the previous step, and let $\delta:\Q^{*} \otimes \E \rightarrow \E$ denote the isomorphism defined in Lemma \ref{lemma.added}.  We show the diagram}
\begin{equation} \label{eqn.digeq0}
\begin{CD}
\N_{k} \otimes \E @>{=}>> \N_{k} \otimes \E \\
@V{\psi_{k} \otimes \E}VV @VV{\gamma_{k+1}}V \\
\mathcal{L} \otimes \A_{lk} \otimes \E @>>{\mu}> \mathcal{L} \otimes \A_{l,k+1}
\end{CD}
\end{equation}
{\it commutes.}  We first prove the rectangle consisting of the diagram
\begin{equation} \label{eqn.digeq1}
\begin{CD}
\N \E @>{\eta}>> \N \Q\Q^{*}\E @>{\alpha}>> \N\E\E^{*}\Q^{*}\E @>{=}>> \N\E\E^{*}\Q^{*}\E \\
@V{\psi_{k} \E}VV @VV{\psi_{k} \Q\Q^{*} \E}V & & @VV{\gamma_{k+1} \E^{*}\Q^{*}\E}V \\
\mathcal{L} \A_{lk} \E @>>{\eta}> \mathcal{L} \A_{lk} \Q \Q^{*} @>>{\alpha}> \mathcal{L} \A_{lk} \E\E^{*}\Q^{*}\E @>>{\mu}> \mathcal{L} \A_{l,k+1} \E^{*}\Q^{*}\E
\end{CD}
\end{equation}
to the left of the diagram
\begin{equation} \label{eqn.digeq2}
\begin{CD}
\N \E\E^{*}\Q^{*}\E @>{\delta}>> \N \E \E^{*} \E @>{\epsilon}>> \N\E \\
@V{\gamma_{k+1} \E^{*} \Q^{*} \E}VV @VV{\gamma_{k+1}\E^{*}\E}V @VV{\gamma_{k+1}}V \\
\mathcal{L} \A_{l,k+1} \E^{*}\Q^{*} \E @>>{\delta}> \mathcal{L} \A_{l,k+1}\E^{*}\E @>>{\epsilon}> \mathcal{L} \A_{l,k+1}
\end{CD}
\end{equation}
commutes.  For, the left square in (\ref{eqn.digeq0}) commutes by the functoriality of the tensor product, the right rectangle in (\ref{eqn.digeq0}) commutes by Step 2, and the squares in (\ref{eqn.digeq1}) commute by the functoriality of the tensor product.  The top row of the rectangle consisting of (\ref{eqn.digeq1}) to the left of (\ref{eqn.digeq2}) is the identity morphism by Lemma \ref{lemma.added}.  Thus, in order to show (\ref{eqn.digeq0}) commutes, it suffices to show the bottom route in the rectangle consisting of (\ref{eqn.digeq1}) to the left of (\ref{eqn.digeq2}) is equal to the multiplication morphism $\mu:\mathcal{L} \otimes \A_{lk} \otimes \E \rightarrow \mathcal{L} \otimes \A_{l,k+1}$.  Thus, it suffices to prove that the diagram
\begin{equation} \label{eqn.digeq3}
\begin{CD}
\A_{lk}\E @>{\eta}>> \A_{lk}\Q\Q^{*} \E @>{\alpha}>> \A_{lk}\E\E^{*}\Q^{*}\E \\
@V{\mu}VV & & @VV{\mu}V \\
\A_{l,k+1} @<<{\epsilon}< \A_{l,k+1} \E^{*}\E @<<{\delta}< \A_{l,k+1} \E^{*}\Q^{*}\E
\end{CD}
\end{equation}
commutes.  By functoriality of the tensor product, the diagram
$$  
\begin{CD}
\A_{lk} \E \E^{*} \Q^{*} \E @>{\delta}>> \A_{lk}\E\E^{*}\E \\
@V{\mu}VV @VV{\mu}V \\
\A_{l,k+1}\E^{*}\Q^{*}\E @>>{\delta}> \A_{l,k+1}\E^{*}\E
\end{CD}
$$
commutes.  Thus, to prove (\ref{eqn.digeq3}) commutes, it suffices to prove
\begin{equation} \label{eqn.digeq4}
\begin{CD}
\A_{lk}\E @>{\eta}>> \A_{lk} \Q\Q^{*} \E @>{\alpha}>> \A_{lk}\E\E^{*}\Q^{*}\E \\
@V{\mu}VV & & @VV{\delta}V \\
\A_{lk} @<<{\epsilon}< \A_{l,k+1}\E^{*}\E @<<{\mu}< \A_{lk}\E\E^{*}\E
\end{CD}
\end{equation}
commutes.  By Lemma \ref{lemma.added}, the composition of the first three maps of the top route is induced by the counit $\eta:\E \rightarrow \E \E^{*} \E$.  Thus, it suffices to show the left square in
$$
\begin{CD}
\A_{lk}\E @>{\eta}>> \A_{lk}\E\E^{*}\E @>{\epsilon}>> \A_{lk}\E \\
@V{\mu}VV @VV{\mu}V @VV{\mu}V \\
\A_{l,k+1} @<<{\epsilon}< \A_{l,k+1}\E^{*}\E @>>{\epsilon}> \A_{l,k+1}
\end{CD}
$$
commutes.  The right square commutes by functoriality of the tensor product.  Since $\eta$ composed with the top route of the right square is $\mu:\A_{lk}\E \rightarrow \A_{l,k+1}$, while $\eta$ composed with the bottom route of the right square is the top route of the left square, the left square commutes.  Thus (\ref{eqn.digeq4}), and hence (\ref{eqn.digeq0}), commutes.
\newline
\newline
\noindent{\it Step 4:  We show $\psi$ is an $\A$-module morphism}.  Since (\ref{eqn.digeq0}) commutes by Step 3, the two right hand squares of
$$
\begin{CD}
\N_{k-1}\A_{k-1,k} @>{\mu}>> \N_{k} @>{\eta}>> \N_{k}\A_{k,k+1}\A_{k,k+1}^{*} @>{=}>> \N_{k}\A_{k,k+1}\A_{k,k+1}^{*} \\
@V{\gamma_{k}}VV @VV{\psi_{k}}V @VV{\psi_{k}\A_{k,k+1}\A_{k,k+1}^{*}}V @VV{\gamma_{k+1}\A_{k,k+1}^{*}}V \\
\mathcal{L}\A_{l,k} @>>{=}> \mathcal{L}\A_{l,k} @>>{\eta}> \mathcal{L} \A_{l,k}\A_{k,k+1}\A_{k,k+1}^{*} @>>{\mu}> \mathcal{L} \A_{l,k+1}\A_{k,k+1}^{*}
\end{CD}
$$
commute.  Since the bottom horizontal is monic by Lemma \ref{lemma.monomult}, in order to show the left square commutes, it suffices to show that 
\begin{equation} \label{eqn.fixit1}
\begin{CD}
\N_{k-1} \otimes \A_{k-1,k} @>{\mu}>> \N_{k} @>{\eta}>> \N_{k} \otimes \A_{k,k+1} \otimes \A_{k,k+1}^{*} \\
@V{\gamma_{k}}VV & & @VV{\gamma_{k+1}\otimes \A_{k,k+1}^{*}}V \\
\mathcal{L} \otimes \A_{l,k} @>>{\eta}> \mathcal{L} \otimes \A_{l,k} \otimes \A_{k,k+1} \otimes \A_{k,k+1}^{*} @>>{\mu}> \mathcal{L} \otimes \A_{l,k+1} \otimes \A_{k,k+1}^{*}
\end{CD}
\end{equation}
commutes.  Since (\ref{eqn.star2}) commutes, the diagram
$$
\begin{CD}
\N_{m} \otimes \A_{m,k-1} \otimes \A_{k-1,k} @>{\mu}>> \N_{m} \otimes \A_{mk} \\
@V{\mu}VV @VV{\gamma_{k}}V \\
\N_{k-1} \otimes \A_{k-1,k} @>>{\gamma_{k}}> \mathcal{L} \otimes \A_{lk}
\end{CD}
$$
commutes.  Thus
$$
\begin{CD}
\N_{k-1} \otimes \A_{k-1,k} \otimes \A_{k,k+1} @>{\mu_{\N}}>> \N_{k} \otimes \A_{k,k+1} \\
@V{\gamma_{k}\otimes \A_{k,k+1}}VV @VV{\gamma_{k+1}}V \\
\mathcal{L} \otimes \A_{l,k} \otimes \A_{k,k+1} @>>> \mathcal{L} \otimes \A_{l,k+1}
\end{CD}
$$
commutes.  By tensoring this diagram on the right by $\A_{k,k+1}^{*}$ and precomposing on the left with the diagram
$$
\begin{CD}
\N_{k-1}\otimes \A_{k-1,k} @>>> \N_{k-1} \otimes \A_{k-1,k}\otimes \A_{k,k+1} \otimes \A_{k,k+1}^{*} \\
@V{\gamma_{k}}VV @VV{\gamma_{k} \otimes \A_{k,k+1} \otimes \A_{k,k+1}^{*}}V \\
\mathcal{L} \otimes \A_{l,k} @>>> \mathcal{L} \otimes \A_{l,k} \otimes \A_{k,k+1} \otimes \A_{k,k+1}^{*}
\end{CD}
$$
whose horizontals are counits, it is readily seen that the left hand square in (\ref{eqn.fixit1}) commutes.  We note that $\psi$ is an $\A$-module map since the left hand square in (\ref{eqn.fixit1}) commutes and (\ref{eqn.digeq0}) commutes.  Since $\A$ is generated in degree 1, the result follows.
\newline
\newline
\noindent{\it Step 5:  We show (\ref{eqn.star1}) commutes.}  We must show 
\begin{equation} \label{eqn.bign}
\begin{CD}
\N_{m} \otimes \A_{mk} @>{\mu}>> \N_{k} \\
@V{=}VV @VV{\psi_{k}}V \\
\N_{m} \otimes \A_{mk} @>>{\gamma_{k}}> \mathcal{L} \otimes \A_{l,k} 
\end{CD}
\end{equation}
commutes.  In order to prove (\ref{eqn.bign}) commutes, it suffices to show
$$
\begin{CD}
\N_{m} \otimes \A_{m,k-1} \otimes \A_{k-1,k} @>{\mu}>> \N_{m} \otimes \A_{mk} @>{\mu}>> \N_{k} \\
& & @V{=}VV @VV{\psi_{k}}V \\
& & \N_{m} \otimes \A_{mk} @>>{\gamma_{k}}> \mathcal{L} \otimes \A_{l,k}
\end{CD}
$$
commutes, since the left horizontal is an epi.  Since $\A$-module multiplication is associative, the diagram
$$
\begin{CD}
\N_{m} \otimes \A_{m,k-1} \otimes \A_{k-1,k} @>>> \N_{m} \otimes \A_{mk} \\
@VVV @VVV \\
\N_{k-1}\otimes \A_{k-1,k} @>>> \N_{k}
\end{CD}
$$
whose arrows are multiplication maps, commutes.  Since (\ref{eqn.star2}) holds, the diagram
\begin{equation} \label{eqn.fixit2}
\begin{CD}
\N_{m} \otimes \A_{m,k-1} \otimes \A_{k-1,k} @>>> \N_{m} \otimes \A_{mk} \\
@VVV @VV{\gamma_{k}}V \\
\N_{k-1} \otimes \A_{k-1,k} @>>{\gamma_{k}}> \mathcal{L} \otimes \A_{lk}
\end{CD}
\end{equation}
whose unlabeled maps are multiplication maps, commutes.  Thus, to prove (\ref{eqn.bign}) commutes, it is enough to show 
\begin{equation} \label{eqn.bign2}
\begin{CD}
\N_{k-1} \otimes \A_{k-1,k} @>{\mu}>> \N_{k} \\
@V{=}VV @VV{\psi_{k}}V \\
\N_{k-1} \otimes \A_{k-1,k} @>>{\gamma_{k}}> \mathcal{L} \otimes \A_{lk}
\end{CD}
\end{equation}
commutes.  Since, by Step 4, $\psi$ is an $\A$-module morphism, the commutivity of (\ref{eqn.bign2}) follows from Step 3.
\end{proof}

\subsection{The internal Hom functor on $\sf{Gr }\A$}
\begin{definition} \label{def.hom}
Let $\mathcal{C}$ be an object in $\mathbb{B}$ and let $\mathcal{M}$ be a graded right $\mathcal{A}$-module.  We define $\boldsymbol {\HU(\mathcal{C},\mathcal{M})}$ to be the $\mathbb{Z}$-graded $\mathcal{O}_{X}$-module whose $k$th component is the equalizer of the diagram

\begin{equation} \label{eqn.homdef}
\begin{CD}
\underset{i}{\Pi}{\mathcal{M}}_{i} \otimes {\mathcal{C}}_{ki}^{*} @>\alpha>> \underset{j}{\Pi}{\mathcal{M}}_{j} \otimes {\mathcal{C}}_{kj}^{*}\\
@V{\beta}VV		@VV{\gamma}V\\
\underset{j}{\Pi}(\underset{i}{\Pi}({\mathcal{M}}_{j} \otimes {\mathcal{A}}_{ij}^{*}) \otimes {\mathcal{C}}_{ki}^{*}) @>>{\delta}> \underset{j}{\Pi}(\underset{i}{\Pi}{\mathcal{M}}_{j} \otimes ({\mathcal{C}}_{ki} \otimes {\mathcal{A}}_{ij})^{*}) 
\end{CD}
\end{equation}
where $\alpha$ is the identity map, $\beta$ is induced by the composition 
\begin{equation} \label{eqn.dualh}
{\mathcal{M}}_{i} \overset{\eta}{\rightarrow} {\mathcal{M}}_{i}\otimes {\mathcal{A}}_{ij} \otimes {\mathcal{A}}_{ij}^{*} \overset{\mu}{\rightarrow} {\mathcal{M}}_{j}\otimes {\mathcal{A}}_{ij}^{*},
\end{equation}
$\gamma$ is induced by the dual (Definition \ref{definition.dual}) of 
$$
{\mathcal{C}}_{ki} \otimes {\mathcal{A}}_{ij} \overset{\mu}{\rightarrow} {\mathcal{C}}_{kj},
$$ 
and $\delta$ is induced by the composition
$$
(\M_{j} \otimes \A_{ij}^{*}) \otimes \C_{ij}^{*} \rightarrow \M_{j}\otimes ({\mathcal{A}}_{ij}^{*}\otimes {\mathcal{C}}_{ki}^{*}) \rightarrow \M_{j} \otimes ({\mathcal{C}}_{ki} \otimes {\mathcal{A}}_{ij})^{*}
$$
whose left arrow is the associativity isomorphism and whose right arrow is induced by the canonical map (\ref{eqn.canonicalisom0}).  If $\mathcal{C}$ is an object of $\mathbb{G}\sf{r} \A$ left-concentrated in degree $k$, we define $\boldsymbol {\HomA (\mathcal{C},\mathcal{M})}$ to be the equalizer of (\ref{eqn.homdef}).
\end{definition}

\begin{remark}
If $\C$ is an object in $\mathbb{B}$ and $\D$ is an object in $\sf{Bimod }\A-\A$ then we can define $\boldsymbol {\HU(\mathcal{C},\mathcal{D})}$ to be the bigraded $\mathcal{O}_{X^{2}}$-module whose $l,k$th component is the equalizer of the diagram

$$
\begin{CD}
\underset{l}{\oplus} \underset{k}{\oplus}\underset{i}{\Pi}{\mathcal{D}}_{li} \otimes {\mathcal{C}}_{ki}^{*} @>\alpha>> \underset{l}{\oplus}\underset{k}{\oplus}\underset{j}{\Pi}{\mathcal{D}}_{lj} \otimes {\mathcal{C}}_{kj}^{*}\\
@V{\beta}VV		@VV{\gamma}V\\
\underset{l}{\oplus}\underset{k}{\oplus} \underset{j}{\Pi}(\underset{i}{\Pi}({\mathcal{D}}_{lj} \otimes {\mathcal{A}}_{ij}^{*}) \otimes {\mathcal{C}}_{ki}^{*}) @>>{\delta}> \underset{l}{\oplus}\underset{k}{\oplus} \underset{j}{\Pi}(\underset{i}{\Pi}{\mathcal{D}}_{lj} \otimes ({\mathcal{C}}_{ki} \otimes {\mathcal{A}}_{ij})^{*}) 
\end{CD}
$$
where $\alpha$ is the identity map, $\beta$ is induced by the composition 
$$
{\mathcal{D}}_{li} \overset{\eta}{\rightarrow} {\mathcal{D}}_{li}\otimes {\mathcal{A}}_{ij} \otimes {\mathcal{A}}_{ij}^{*} \overset{\mu}{\rightarrow} {\mathcal{D}}_{lj}\otimes {\mathcal{A}}_{ij}^{*},
$$
$\gamma$ is induced by the dual of 
$$
{\mathcal{C}}_{ki} \otimes {\mathcal{A}}_{ij} \overset{\mu}{\rightarrow} {\mathcal{C}}_{kj},
$$ 
and $\delta$ is induced by the composition 
$$
(\D_{lj} \otimes \A_{ij}^{*}) \otimes \C_{ij}^{*} \rightarrow \D_{lj}\otimes ({\mathcal{A}}_{ij}^{*}\otimes {\mathcal{C}}_{ki}^{*}) \rightarrow \D_{lj} \otimes ({\mathcal{C}}_{ki} \otimes {\mathcal{A}}_{ij})^{*}
$$
whose left arrow is the associativity isomorphism and whose right arrow is induced by the canonical map (\ref{eqn.canonicalisom0}).  Since we will not use this object in what follows, we will not study its properties.
\end{remark}

\begin{lemma} \label{lemma.homer}
Let $\F$ be an object of $\mathbb{G}\sf{r} \A$ left-concentrated in degree $k$ and let $\mathcal{L}$ be an $\mathcal{O}_{X}$-module.  Then $\Homa(\mathcal{O}_{X} \otimes \mathcal{F},\mathcal{M})$ is the equalizer of the diagram
\begin{equation} \label{eqn.homer}
\begin{CD}
\underset{i}{\Pi}\Homo(\mathcal{L} \otimes {\mathcal{F}}_{ki},{\mathcal{M}}_{i}) @>=>> \underset{j}{\Pi}\Homo(\mathcal{L} \otimes {\mathcal{F}}_{kj},{\mathcal{M}}_{j}) \\
@VVV @VVV \\
\underset{j}{\Pi}(\underset{i}{\Pi}\Homo(\mathcal{L} \otimes {\mathcal{F}}_{ki},{\mathcal{M}}_{j}\otimes {\mathcal{A}}_{ij}^{*})) @>>{\cong}> \underset{j}{\Pi}(\underset{i}{\Pi}\Homo(\mathcal{L} \otimes ({\mathcal{F}}_{ki}\otimes {\mathcal{A}}_{ij}),{\mathcal{M}}_{j})) 
\end{CD}
\end{equation}
whose left vertical is induced by the composition 
$$
{\mathcal{M}}_{i} \overset{\eta}{\rightarrow} {\mathcal{M}}_{i}\otimes {\mathcal{A}}_{ij} \otimes {\mathcal{A}}_{ij}^{*} \overset{\mu}{\rightarrow} {\mathcal{M}}_{j}\otimes {\mathcal{A}}_{ij}^{*},
$$
whose right vertical is induced by  
$$
{\mathcal{F}}_{ki} \otimes {\mathcal{A}}_{ij} \overset{\mu}{\rightarrow} {\mathcal{F}}_{kj},
$$
and whose bottom horizontal is induced by the composition of isomorphisms
\begin{align*}
\Homo(\mathcal{L} \otimes {\mathcal{F}}_{ki},{\mathcal{M}}_{j}\otimes {\mathcal{A}}_{ij}^{*}) & \cong \Homo((\mathcal{L} \otimes \F_{ki}) \otimes \A_{ij},\M_{j}) \\
& \cong \Homo(\mathcal{L} \otimes ({\mathcal{F}}_{ki}\otimes {\mathcal{A}}_{ij}),{\mathcal{M}}_{j})
\end{align*}
where the last isomorphism is induced by the associativity of the tensor product.
\end{lemma}

\begin{proof}
Suppose $\underset{i}{\Pi}f_{i} \in \underset{i}{\Pi}\Homo(\mathcal{L} \otimes {\mathcal{F}}_{ki},{\mathcal{M}}_{i})$.  Then $\underset{i}{\Pi}f_{i}$ maps, via the top route, to the product $\underset{j}{\Pi}(\underset{i}{\Pi}g_{ij})$, where $g_{ij}$ is the composition
$$
\mathcal{L} \otimes ({\mathcal{F}}_{ki} \otimes {\mathcal{A}}_{ij}) \overset{\mu}{\rightarrow} \mathcal{L} \otimes {\mathcal{F}}_{kj} \overset{f_{j}}{\rightarrow}{\mathcal{M}}_{j}.
$$
$\underset{i}{\Pi}f_{i}$ maps, via the bottom route, to $\underset{j}{\Pi}(\underset{i}{\Pi}h_{ij})$, where $h_{ij}$ is the composition
$$
\mathcal{L} \otimes (\F_{ki} \otimes \A_{ij}) \rightarrow (\mathcal{L} \otimes {\mathcal{F}}_{ki}) \otimes {\mathcal{A}}_{ij} \overset{f_{i} \otimes {\mathcal{A}}_{ij}}{\longrightarrow} {\mathcal{M}}_{i} \otimes {\mathcal{A}}_{ij} \overset{\eta}{\longrightarrow}
$$
$$
({\mathcal{M}}_{i} \otimes ({\A}_{ij} \otimes {\A}_{ij}^{*}) \otimes {\mathcal{A}}_{ij}) \overset{\mu \otimes {\A}_{ij}^{*} \otimes {\A}_{ij}}{\longrightarrow} ({\mathcal{M}}_{j} \otimes {\A}_{ij}^{*}) \otimes {\mathcal{A}}_{ij} \overset{\epsilon}{\longrightarrow} {\M}_{j}.
$$
Thus, $\underset{i}{\Pi}f_{i}$ is in the equalizer of (\ref{eqn.homer}) if and only if the diagram
\begin{equation} \label{eqn.homer2}
\begin{CD}
(\mathcal{L} \otimes {\mathcal{F}}_{ki}) \otimes {\mathcal{A}}_{ij} @>>> \mathcal{L} \otimes {\mathcal{F}}_{kj} @>{f_{j}}>> {\mathcal{M}}_{j} \\
@V{f_{i}\otimes {\A}_{ij}}VV & & @AA{\epsilon}A \\
{\mathcal{M}}_{i} \otimes {\mathcal{A}}_{ij} @>>{\eta}> {\mathcal{M}}_{i} \otimes ({\A}_{ij} \otimes {\A}_{ij}^{*}) \otimes {\mathcal{A}}_{ij} @>>{\mu \otimes {\A}_{ij}^{*} \otimes {\A}_{ij}}>{\mathcal{M}}_{j} \otimes {\A}_{ij}^{*} \otimes {\mathcal{A}}_{ij} 
\end{CD}
\end{equation}
whose top left arrow is the composition
\begin{equation} \label{eqn.refer}
(\mathcal{L} \otimes {\mathcal{F}}_{ki}) \otimes {\mathcal{A}}_{ij} \rightarrow \mathcal{L} \otimes ({\mathcal{F}}_{ki} \otimes \A_{ij}) \overset{\mu}{\rightarrow} \mathcal{L} \otimes \F_{kj}
\end{equation}
with left arrow induced by associativity of tensor product, commutes for all $i$ and $j$.  However, since 
$$
\begin{CD}
{\M}_{i} \otimes {\A}_{ij} @>{\eta}>> ({\mathcal{M}}_{i} \otimes ({\A}_{ij} \otimes {\A}_{ij}^{*}) \otimes {\mathcal{A}}_{ij}) @>{\epsilon}>> {\M}_{i} \otimes {\A}_{ij} \\
& & @V{\mu \otimes {\A}_{ij}^{*} \otimes {\A}_{ij}}VV @VV{\mu}V \\
& & ({\M}_{j} \otimes {\A}_{ij}^{*}) \otimes {\A}_{ij} @>>{\epsilon}> {\M}_{j}
\end{CD}
$$
commutes, $\underset{i}{\Pi}f_{i}$ is in the equalizer of (\ref{eqn.homer}) if and only if the diagram
$$
\begin{CD}
(\mathcal{L} \otimes {\mathcal{F}}_{ki}) \otimes {\mathcal{A}}_{ij} @>{f_{i} \otimes {\A}_{ij}}>> {\M}_{i} \otimes {\A}_{ij}\\
@VVV @VV{\mu}V \\
\mathcal{L} \otimes {\mathcal{F}}_{kj} @>>{f_{j}}> {\M}_{j}
\end{CD}
$$
whose left vertical is the morphism (\ref{eqn.refer}), commutes for all $i$ and $j$. 
\end{proof}
The following result provides some justification for Definition \ref{def.hom}.
\begin{proposition} \label{prop.global}
If $\F$ is an object of $\mathbb{G}\sf{r} \A$ which is left-concentrated in degree $k$, and if $\mathcal{L}$ is an $\mathcal{O}_{X}$-module, $\Homo(\mathcal{L}, \HomA (\mathcal{F}, \mathcal{M})) \cong \operatorname{Hom}_{{\sf{Gr }}\mathcal{A}}(\mathcal{L} \otimes \mathcal{F}, \mathcal{M}).$  In particular, if $\Gamma(X,-):{\sf Mod }X \rightarrow {\sf Mod }\Gamma(X,\mathcal{O}_{X})$ is the global sections functor, 
$$
\Gamma(X, \HomA(\F, \M)) \cong \Homa(\mathcal{O}_{X} \otimes \F, \M).
$$
\end{proposition}

\begin{proof}
Since $\Hom_{{\mathcal{O}}_{X}}(\mathcal{L},-)$ is left exact, $\Homo(\mathcal{L},\HomA(\F,\M))$ is isomorphic to the equalizer of the diagram
\begin{equation} \label{eqn.gamma}
\begin{CD}
\underset{i}{\Pi} (\Homo(\mathcal{L}, {\M}_{i} \otimes {\F}_{ki}^{*})) @>=>> \underset{j}{\Pi} (\Homo(\mathcal{L}, {\M}_{j} \otimes {\F}_{kj}^{*})) \\
@VVV @VVV \\
\underset{j}{\Pi}(\underset{i}\Pi (\Homo(\mathcal{L}, ({\M}_{j} \otimes {\A}_{ij}^{*})\otimes {\F}_{ki}^{*}))) @>>{\cong}> \underset{j}{\Pi} (\underset{i}{\Pi} (\Homo(\mathcal{L}, {\M}_{j} \otimes ({\F}_{ki}\otimes {\A}_{ij})^{*})))
\end{CD}
\end{equation}
whose left vertical is induced by the composition
$$
\M_{i} \overset{\eta}{\rightarrow} \M_{i} \otimes \A_{ij} \otimes \A_{ij}^{*} \overset{\mu}{\rightarrow} \M_{j} \otimes \A_{ij}^{*},
$$
whose right vertical is induced by the dual of
$$
\F_{ki} \otimes \A_{ij} \overset{\mu}{\rightarrow} \F_{kj},
$$
and whose bottom horizontal is induced by 
\begin{equation} \label{eqn.newstar}
(\M_{j} \otimes \A_{ij}^{*}) \otimes \F_{ki}^{*} \rightarrow \M_{j}\otimes ({\mathcal{A}}_{ij}^{*}\otimes {\mathcal{F}}_{ki}^{*}) \rightarrow \M_{j} \otimes ({\mathcal{F}}_{ki} \otimes {\mathcal{A}}_{ij})^{*}.
\end{equation}
To complete the proof, we must show the equalizer of (\ref{eqn.gamma}) is isomorphic to the equalizer of (\ref{eqn.homer}).  We think of (\ref{eqn.gamma}) and (\ref{eqn.homer}) as the top and bottom, respectively, of a cube whose vertical arrows are adjointness isomorphisms.  To show the equalizers of (\ref{eqn.gamma}) and (\ref{eqn.homer}) are isomorphic, it suffices to show the vertical faces of this cube commute.  By naturality of adjointness isomorphisms, three of these vertical faces obviously commute.  To complete the proof, we must show that the diagram
\begin{equation} \label{eqn.suffices}
\begin{CD}
\Homo(\mathcal{L},(\M_{j} \otimes \A_{ij}^{*}) \otimes \F_{ki}^{*}) @>>> \Homo(\mathcal{L},\M_{j} \otimes (\F_{ki} \otimes \A_{ij})^{*}) \\
@VVV @VVV \\
\Homo(\mathcal{L} \otimes \F_{ki}, \M_{j} \otimes \A_{ij}^{*}) @>>> \Homo(\mathcal{L} \otimes (\F_{ki} \otimes \A_{ij}), \M_{j})
\end{CD}
\end{equation}
whose left and right arrows are adjointness isomorphisms, whose top arrow is induced by (\ref{eqn.newstar}), and whose bottom arrow is the composition
\begin{align*}
\Homo(\mathcal{L} \otimes \F_{ki}, \M_{j} \otimes \A_{ij}^{*}) & \cong \Homo((\mathcal{L} \otimes \F_{ki}) \otimes \A_{ij},\M_{j}) \\
& \cong \Homo(\mathcal{L} \otimes (\F_{ki} \otimes \A_{ij}), \M_{j})
\end{align*}
whose first isomorphism is the adjointness isomorphism and whose second isomorphism is induced by associativity, commutes.  The diagram (\ref{eqn.suffices}) can be rewritten as the diagram
$$
\begin{CD}
\Homo((\mathcal{L} \otimes \F_{ki})\otimes \A_{ij},\M_{j}) @>>> \Homo(\mathcal{L}, (\M_{j}\otimes \A_{ij}^{*}) \otimes \F_{ki}^{*}) \\
@VVV @VVV \\
\Homo(\mathcal{L} \otimes (\F_{ki} \otimes \A_{ij}),\M_{j}) @>>> \Homo(\mathcal{L},\M_{j} \otimes (\F_{ki} \otimes \A_{ij})^{*})
\end{CD}
$$
whose top horizontal is the composition of adjointness isomorphisms
\begin{align*}
\Homo((\mathcal{L} \otimes \F_{ki})\otimes \A_{ij},\M_{j}) & \cong \Homo(\mathcal{L} \otimes \F_{ki},\M_{j}\otimes \A_{ij}^{*}) \\
& \cong \Homo(\mathcal{L}, (\M_{j}\otimes \A_{ij}^{*}) \otimes \F_{ki}^{*}), 
\end{align*}
whose bottom horizontal is the adjointness isomorphism, whose left vertical is induced by associativity, and whose right vertical is induced by (\ref{eqn.newstar}).  This commutes by Lemma \ref{lemma.commute}, and the result follows.
\end{proof}

\begin{proposition} \label{prop.homo}
If $\mathcal{M}$ is an object in ${\sf{Gr }}\mathcal{A}$ and $\mathcal{C}$ is an object in $\mathbb{B}$, $\HU (\mathcal{C},\mathcal{M})$ inherits a graded right $\mathcal{A}$-module structure from the left $\mathcal{A}$-module structure of $\mathcal{C}$, making $\HU (-,-):{\mathbb{B}}^{op} \times {\sf{Gr}}\mathcal{A} \rightarrow {\sf{Gr}}\mathcal{A}$ a bifunctor.
\end{proposition}

\begin{proof}
We first show that $\HU(\C,-)$ is a functor from $\sf{Gr }\A$ to $\sf{Gr }\A$.  We begin by constructing a map of $\mathcal{O}_{X}$-modules
\begin{equation} \label{eqn.muij}
\mu_{ij}:\HU(\C,\M)_{i} \otimes \A_{ij} \rightarrow \HU(\C,\M)_{j}.
\end{equation}
To this end, we note that $\HU(\C,\M)_{i}$ is an $\mathcal{O}_{X}$-submodule of $\underset{l}{\Pi}\M_{l}\otimes \C_{il}^{*}$.  Thus, we have a map
\begin{equation} \label{eqn.mu1}
\HU(\C,\M)_{i}\otimes \A_{ij} \rightarrow (\underset{l}{\Pi}\M_{l}\otimes \C_{il}^{*})\otimes \A_{ij} \overset{\cong}{\rightarrow} \underset{l}{\Pi}\M_{l}\otimes (\C_{il}^{*}\otimes \A_{ij})
\end{equation}
whose right composite is induced by associativity of the tensor product.  The left $\A$-module structure of $\C$ yields a multiplication map $\A_{ij} \otimes \C_{jl} \rightarrow \C_{il}$.  Dualizing  gives us a map
\begin{equation} \label{eqn.mu2}
\C_{il}^{*} \rightarrow (\A_{ij}\otimes \C_{jl})^{*} \overset{\cong}{\rightarrow} \C_{jl}^{*} \otimes \A_{ij}^{*}
\end{equation}
whose right composite is the canonical isomorphism.  Tensoring (\ref{eqn.mu2}) by $\A_{ij}$ on the right and composing with the map induced by the counit $\C_{jl}^{*} \otimes \epsilon:\C_{jl}^{*} \otimes \A_{ij}^{*} \otimes \A_{ij} \rightarrow \C_{jl}^{*}$ gives us a map
\begin{equation} \label{eqn.mu25}
\C_{il}^{*}\otimes \A_{ij} \rightarrow (\C_{jl}^{*} \otimes \A_{ij}^{*}) \otimes \A_{ij} \rightarrow \C_{jl}^{*}.
\end{equation}
Tensoring (\ref{eqn.mu25}) on the left with $\M_{l}$, taking the product over all $l$, and composing with (\ref{eqn.mu1}) gives us a map
\begin{equation} \label{eqn.mu3}
\HU(\C,\M)_{i} \otimes \A_{ij} \rightarrow \underset{l}{\Pi}(\M_{l} \otimes \C_{jl}^{*}).
\end{equation}
We must show this map factors through $\HU(\C,\M)_{j}$.  Let
\begin{equation} \label{eqn.dig1}
\begin{CD}
\underset{l}{\Pi}(\M_{l}\C_{il}^{*})\A_{ij} @>>> \underset{l}{\Pi}(\M_{l}(\C_{jl}^{*}\A_{ij}^{*}))\A_{ij} @>{\epsilon}>> \underset{l}{\Pi}\M_{l}\C_{jl}^{*} \\
@V{=}VV @VV{=}V @VV{=}V \\
\underset{m}{\Pi}(\M_{m} \C_{im}^{*}) \A_{ij} @>>> \underset{m}{\Pi}(\M_{m}(\C_{jm}^{*}\A_{ij}^{*})) \A_{ij} @>>{\epsilon}> \underset{m}{\Pi} \M_{m}\C_{jm}^{*} \\
@V{\mu^{*}}VV @VV{\mu^{*}}V @VV{\mu^{*}}V \\
\underset{m}{\Pi}\underset{l}{\Pi}(\M_{m}(\C_{il}\A_{lm})^{*})\A_{ij} @>>> \underset{m}{\Pi}\underset{l}{\Pi}(\M_{m}((\C_{jl}\A_{lm})^{*}\A_{ij}^{*})\A_{ij}) @>>> \underset{m}{\Pi}\underset{l}{\Pi}\M_{m}(\C_{jl}\A_{lm})^{*}
\end{CD}
\end{equation}
be the diagram whose top and middle rows are induced by (\ref{eqn.mu25}) and whose bottom row is induced by the composition
\begin{equation} \label{eqn.lump}
(\C_{il} \otimes \A_{lm})^{*} \overset{\mu^{*}}{\rightarrow} ((\A_{ij} \otimes \C_{jl}) \otimes \A_{lm})^{*} \rightarrow (\C_{jl} \otimes \A_{lm})^{*} \otimes \A_{ij}^{*}
\end{equation}
(whose right composite is the canonical isomorphism).  Let
\begin{equation} \label{eqn.dig2}
\begin{CD}
\underset{l}{\Pi}(\M_{l}\C_{il}^{*})\A_{ij} @>>> \underset{l}{\Pi}(\M_{l}(\C_{jl}^{*}\A_{ij}^{*})\A_{ij}) @>{\epsilon}>> \underset{l}{\Pi}\M_{l}\C_{jl}^{*} \\
@VVV @VVV @VVV \\
\underset{m}{\Pi}\underset{l}{\Pi}((\M_{m}\A_{lm}^{*})\C_{il}^{*})\A_{ij} @>>> \underset{m}{\Pi}\underset{l}{\Pi}((\M_{m}\A_{lm}^{*})(\C_{jl}^{*}\A_{ij}^{*})\A_{ij}) @>>{\epsilon}> \underset{m}{\Pi}\underset{l}{\Pi}(\M_{m}\A_{lm}^{*})\C_{jl}^{*} \\
@V{\cong}VV @VV{\cong}V @VV{\cong}V \\
\underset{m}{\Pi}\underset{l}{\Pi}(\M_{m}(\C_{il}\A_{lm})^{*})\A_{ij} @>>> \underset{m}{\Pi}\underset{l}{\Pi}((\M_{m}(\C_{jl}\A_{lm})^{*})\A_{ij}^{*}) \A_{ij} @>>{\epsilon}> \underset{m}{\Pi}\underset{l}{\Pi}\M_{m}(\C_{jl}\A_{lm})^{*}
\end{CD}
\end{equation}
be the diagram whose top and middle rows are induced by (\ref{eqn.mu25}), whose bottom left horizontal is induced by (\ref{eqn.lump}), whose top verticals are induced by the composition
$$
\M_{l} \overset{\eta}{\rightarrow} \M_{l} \otimes \A_{lm} \otimes \A_{lm}^{*} \rightarrow \M_{m} \otimes \A_{lm}^{*},
$$
and whose bottom verticals are induced by the canonical isomorphisms (\ref{eqn.canonicalisom0}).  We leave it as an exercise for the reader to prove, using the universal property of equalizer, that, in order to prove (\ref{eqn.mu3}) factors through $\HU(\C,\M)_{j}$, it suffices to show that the outer circuits of (\ref{eqn.dig1}) and (\ref{eqn.dig2}) commute.

The upper squares in (\ref{eqn.dig1}) commute since the vertical maps are identity maps.  The right squares in (\ref{eqn.dig1}) and (\ref{eqn.dig2}), and the upper left square in (\ref{eqn.dig2}), commute by the functoriality of the tensor product.  Thus, in order to show that (\ref{eqn.mu3}) factors through $\HU(\C,\M)_{j}$, we must show that the lower-left squares in (\ref{eqn.dig1}) and (\ref{eqn.dig2}) commute, i.e. we must show
\begin{equation} \label{eqn.dig3}
\begin{CD}
\C_{im}^{*} @>>> \C_{jm}^{*} \otimes \A_{ij}^{*} \\
@VVV @VVV \\
(\C_{il} \otimes \A_{lm})^{*} @>>> (\C_{jl} \otimes \A_{lm})^{*} \otimes \A_{ij}^{*}
\end{CD}
\end{equation} 
whose maps are induced by $\mu_{\C}$, commutes, and 
\begin{equation} \label{eqn.dig4}
\begin{CD}
\A_{lm}^{*} \otimes \C_{il}^{*} @>>> \A_{lm}^{*} \otimes (\C_{jl}^{*} \otimes \A_{ij}^{*}) \\
@VVV @VVV \\
(\C_{il}\otimes \A_{lm})^{*} @>>> (\C_{jl}\otimes \A_{lm})^{*} \otimes \A_{ij}^{*}
\end{CD}
\end{equation}
whose horizontal maps are induced by $\mu_{\C}$ and whose vertical maps are canonical isomorphisms, commutes.  By Corollary \ref{cor.tweed}, in order to show (\ref{eqn.dig3}) commutes, it suffices to show the diagram
$$
\begin{CD}
\C_{im}^{*} @>>> (\A_{ij}\otimes \C_{jm})^{*} \\
@VVV @VVV \\
(\C_{il}\otimes \A_{lm})^{*} @>>> ((\A_{ij} \otimes \C_{jl}) \otimes \A_{lm})^{*}
\end{CD}
$$
whose arrows are induced by multiplication, commutes.  This follows from the functoriality of $(-)^{*}$.  

Expanding (\ref{eqn.dig4}), we have a diagram
$$
\begin{CD}
\A_{lm}^{*} \otimes \C_{il}^{*} @>{\mu^{*}}>> \A_{lm}^{*} \otimes (\A_{ij} \otimes \C_{jl})^{*} @>>> \A_{lm}^{*} \otimes (\C_{jl}^{*}\otimes \A_{ij}^{*}) \\
@VVV @VVV @VVV\\
(\C_{il}\otimes \A_{lm})^{*} @>>{\mu^{*}}> ((\A_{ij} \otimes \C_{jl}) \otimes \A_{lm})^{*} @>>> (\C_{jl} \otimes \A_{lm})^{*}\otimes \A_{ij}^{*}
\end{CD}
$$
whose right vertical is a composition of an associativity isomorphism with the canonical isomorphism (\ref{eqn.canonicalisom0}).  The left square commutes by Corollary \ref{cor.tweed}, while the right square commutes by Corollary \ref{corollary.tweed3}.

We defer the proof that the map (\ref{eqn.muij}) is associative and compatible with scalar multiplication to Section 7. 

We omit the routine proofs that $\HU(\C,-)$ and $\HU(-,\M)$ are functorial, that $\HU(\C,\operatorname{id}_{\M})=\operatorname{id}_{\HU(\C,\M)}$ and that $\HU(\C,-)$ is compatible with composition.
\end{proof}

\subsection{Adjointness of internal Hom and Tensor}
\begin{definition}
We say that $F: \sf{A} \times \sf{B} \rightarrow \sf{A}$ and $G: {\sf{B}}^{op} \times \sf{A} \rightarrow \sf{A}$ are {\bf conjugate bifunctors} if, for every $\mathcal{C} \in \sf{B}$, $F(-,\mathcal{C})$ has a right adjoint $G(\mathcal{C},-)$ via an adjunction
$$
\phi:  {\operatorname{Hom}}_{\sf{A}}(F(\mathcal{M},\mathcal{C}),\mathcal{N}) \rightarrow {\operatorname{Hom}}_{\sf{A}}(\mathcal{M},G(\mathcal{C},\mathcal{N}))
$$ 
which is natural in $\mathcal{M}$, $\mathcal{N}$ {\it and} $\mathcal{C}$.
\end{definition}

\begin{proposition} \label{prop.tensor}
Let $\C$ be an object in $\mathbb{B}$.  
\begin{enumerate}
\item{}
There exist natural transformations
and
$$
\eta:  {\sf id}_{{\sf{Gr }}\mathcal{A}} \rightarrow \HU (\mathcal{C},-\Ten \mathcal{C})
$$
$$
\epsilon: \HU (\mathcal{C},-) \Ten \mathcal{C} \rightarrow {\sf{id }}_{{\sf{Gr }}\mathcal{A}}
$$
making $(-\Ten \mathcal{C}, \HU (\mathcal{C}, -) , \eta, \epsilon): {\sf{Gr }}\mathcal{A} \rightarrow {\sf{Gr }}\mathcal{A}$ an adjunction.

\item{}
There exists a unique way to make $-\Ten -:{\sf{Gr }}\mathcal{A} \times \mathbb{B} \rightarrow {\sf{Gr }}\mathcal{A}$ a bifunctor such that $-\Ten-$ and $\HU (-,-)$ are conjugate.
\end{enumerate}
\end{proposition}

\begin{proof}
We construct an $\A$-module map
$$
\eta_{\M}:\M \rightarrow \HU(\C,\M\Ten\C) 
$$
natural in $\M$.  To begin, for each $k$ we construct a map
\begin{equation} \label{eqn.startmap}
\M_{k} \rightarrow (\M \Ten \C)_{i} \otimes \C_{ki}^{*}
\end{equation}
for all $i$.  In order to construct (\ref{eqn.startmap}), we construct a map
\begin{equation} \label{eqn.etastart}
\M_{k} \rightarrow (\underset{m}{\oplus}\M_{m} \otimes \C_{mi})\otimes \C_{ki}^{*}.
\end{equation}
The map (\ref{eqn.etastart}) is just the composition
$$
\M_{k} \overset{\eta}{\rightarrow} \M_{k}\otimes \C_{ki}\otimes \C_{ki}^{*} \rightarrow (\underset{m}{\oplus} \M_{m} \otimes \C_{mi}) \otimes \C_{ki}^{*}
$$
whose right arrow is induced by inclusion of a direct summand.  Thus, we have a map
$$
\M_{k} \rightarrow (\underset{m}{\oplus}\M_{m} \otimes \C_{mi}) \otimes \C_{ki}^{*} \rightarrow (\M \Ten \C)_{i} \otimes \C_{ki}^{*}
$$
for every $i$, and hence a map
$$
\psi_{k}: \M_{k} \rightarrow \underset{i}{\Pi}(\M \Ten \C)_{i} \otimes \C_{ki}^{*}.
$$
In order to show this map induces a map
$$
\eta_{\M,k}:\M_{k} \rightarrow \HU(\C,\M \Ten \C)_{k},
$$
we must show the diagram
\begin{equation} \label{eqn.diga1}
\begin{CD}
\M_{k} \\
@V{\psi_{k}}VV \\
\underset{i}{\Pi}(\M \Ten \C)_{i} \otimes \C_{ki}^{*} @>>> \underset{j}{\Pi}(\M \Ten \C)_{j} \otimes \C_{kj}^{*} \\
@VVV @VVV \\
\underset{j}{\Pi}\underset{i}{\Pi}((\M \Ten \C)_{j} \otimes \A_{ij}^{*}) \otimes \C_{ki}^{*} @>>{\cong}> \underset{j}{\Pi}\underset{i}{\Pi}((\M \Ten \C)_{j} \otimes (\C_{ki} \otimes \A_{ij})^{*})
\end{CD}
\end{equation}
whose unlabeled maps are those in (\ref{eqn.homdef}), commutes.  In order to prove this, it suffices to show the projection of (\ref{eqn.diga1}) onto its $i$,$j$ component, the left square in
$$
\begin{CD}
\M_{k} @>{=}>> \M_{k} @>{=}>> \M_{k} \\
@V{\eta}VV @VV{\eta}V @VV{\eta}V \\
(\M_{k} \C_{ki}) \C_{ki}^{*} & & (\M_{k} \C_{kj})  \C_{kj}^{*} & & \M_{k}  (\C_{kj} \C_{kj}^{*})\\
@VVV @VV{\mu^{*}}V @VV{\mu^{*}}V \\
((\M_{k}  \C_{kj})  \A_{ij}^{*})  \C_{ki}^{*} @>>{\cong}> (\M_{k}  \C_{kj})  (\C_{ki}  \A_{ij})^{*} @>>> \M_{k}  (\C_{kj} (\C_{ki}  \A_{ij})^{*})
\end{CD}
$$
whose bottom right horizontal is induced by associativity, commutes.  Since the outer circuit of this diagram equals the outer circuit of the diagram 
$$
\begin{CD}
\M_{k} @>{=}>> \M_{k} @>{=}>> \M_{k} \\
@V{\eta}VV @VV{\eta}V @VV{\eta}V \\
\M_{k}  (\C_{ki}  (\A_{ij}  \A_{ij}^{*})  \C_{ki}^{*}) @>>> \M_{k}  ((\C_{ki}  \A_{ij})  (\C_{ki}  \A_{ij})^{*}) & & \M_{k}  (\C_{kj}  \C_{kj}^{*}) \\
@V{\mu}VV @VV{\mu}V @VV{\mu^{*}}V \\
\M_{k}  ((\C_{kj}  \A_{ij}^{*})  \C_{ki}^{*}) @>>> \M_{k}  (\C_{kj}  (\C_{ki}  \A_{ij})^{*}) @>>{=}> \M_{k}  (\C_{kj}  (\C_{ki}  \A_{ij})^{*})
\end{CD}
$$
whose unlabeled arrows are canonical isomorphisms, it suffices to show each circuit in this diagram commutes.  The right rectangle commutes by Corollary \ref{corollary.tweed}, the bottom-left square commutes by the functoriality of the tensor product and the upper-left square commutes by Lemma \ref{lem.com1}.

The fact that $\eta$ is natural follows from a straightforward computation, which we omit.  We defer the proof the $\eta$ is compatible with $\A$-module multiplication to Section 7.

We next construct a map $\epsilon_{\M}:\HU(\M,\C) \Ten \C \rightarrow \M$ natural in $\M$.  To begin, we construct a natural map 
$$
\epsilon_{k}:(\HU(\C,\M)\Ten \C)_{k} \rightarrow \M_{k}.
$$
To construct $\epsilon_{k}$, it suffices to construct a map 
$$
\delta_{k}: \underset{m}{\oplus}(\HU(\C,\M))_{m} \otimes \C_{mk} \rightarrow \M_{k}
$$
such that the diagram
\begin{equation} \label{eqn.epsilon1}
\begin{CD}
\underset{l}{\oplus}\underset{m}{\oplus}((\HU(\C,\M))_{l} \otimes \A_{lm}) \otimes \C_{mk} @>>> \underset{m}{\oplus}(\HU(\C,\M))_{m} \otimes \C_{mk} \\
@VVV @VV{=}V \\
\underset{l}{\oplus}(\HU(\C,\M))_{l}\otimes \C_{lk} @>>{=}> \underset{m}{\oplus}(\HU(\C,\M))_{m} \otimes \C_{mk} \\
& & @VV{\delta}V \\
& & \M_{k}
\end{CD}
\end{equation}
whose unlabeled arrows are induced by multiplication, commutes.  We define $\delta_{k}$ as the map induced by the composition
$$
(\underset{i}{\Pi}\M_{i} \otimes \C_{mi}^{*})\otimes \C_{mk} \overset{\epsilon}{\rightarrow} \M_{k}\otimes \C_{mk}^{*} \otimes \C_{mk} \rightarrow \M_{k}
$$
which we call $\gamma_{k}$.  To show (\ref{eqn.epsilon1}) commutes, it suffices to show the diagram
$$
\begin{CD}
\underset{l}{\oplus}\underset{m}{\oplus}((\underset{i}{\Pi}\M_{i} \otimes \C_{li}^{*}) \otimes \A_{lm}) \otimes \C_{mk} @>>> \underset{m}{\oplus}(\underset{i}{\Pi}\M_{i} \otimes \C_{mi}^{*})\otimes \C_{mk} \\
@VVV @VV{\gamma_{k}}V \\
\underset{l}{\oplus}(\underset{i}{\Pi}\M_{i} \otimes \C_{li}^{*}) \otimes \C_{lk} @>>{\gamma_{k}}> \M_{k}
\end{CD}
$$
whose left vertical is induced by multiplication and whose top horizontal is induced by the composition
\begin{equation} \label{eqn.composite}
\C_{li}^{*} \otimes \A_{lm} \overset{\mu^{*}}{\rightarrow} (\A_{lm} \otimes \C_{mi})^{*} \otimes \A_{lm} \overset{\cong}{\rightarrow} \C_{mi}^{*} \otimes \A_{lm}^{*} \otimes \A_{lm} \overset{\epsilon}{\rightarrow} \C_{mi}^{*}
\end{equation}
commutes for all $l$ and $m$.  Fixing $l,m$, we must show the diagram
$$
\begin{CD}
\underset{i}{\Pi}((\M_{i}\otimes \C_{li}^{*})\otimes \A_{lm}) \otimes \C_{mk} @>>> \underset{i}{\Pi}(\M_{i} \otimes \C_{mi}^{*}) \otimes \C_{mk} \\
@V{\mu}VV @VVV \\
\underset{i}{\Pi}(\M_{i}\otimes \C_{li}^{*})\otimes \C_{lk} & & (\M_{k}\otimes\C_{mk}^{*})\otimes \C_{mk} \\
@VVV @VV{\epsilon}V \\
(\M_{k} \otimes \C_{lk}^{*}) \otimes \C_{lk} @>>{\epsilon}> \M_{k}
\end{CD}
$$
whose top horizontal is induced by (\ref{eqn.composite}) and whose bottom left-vertical and top right vertical are projections, commutes.  This is equivalent to showing the diagram
$$
\begin{CD}
((\M_{k}\otimes \C_{lk}^{*})\otimes \A_{lm})\otimes \C_{mk} @>>> (\M_{k} \otimes \C_{mk}^{*}) \otimes \C_{mk} \\
@V{\mu}VV @VV{\epsilon}V \\
(\M_{k} \otimes \C_{lk}^{*}) \otimes \C_{lk} @>>{\epsilon}> \M_{k}
\end{CD}
$$
whose top horizontal is induced by (\ref{eqn.composite}) and whose left vertical is induced by multiplication, commutes.  This follows from Corollary \ref{corollary.tweed}.

The fact that $\epsilon$ is natural follows from a straightforward computation, which we omit.  We defer the proof the $\epsilon$ is compatible with $\A$-module multiplication to Appendix 2.

The proofs that 
$$
(\epsilon * (- \Ten \C))\circ ((-\Ten \C)*\eta)=(- \Ten \C)
$$
and 
$$
(\HU(\C,-)*\epsilon) \circ (\eta* \HU(\C,-))=\HU(\C,-)
$$
are straightforward but tedious, and we omit them.

We have established that $(-\Ten \C, \HU(\C,-), \eta, \epsilon)$ forms an adjunction.  Thus, the second part of the proposition follows from \cite[Theorem 3, p. 102]{mac}.
\end{proof}
We endow $- \Ten -$ with the bifunctor structure from (3) in the Proposition.

\begin{lemma} \label{lem.general}
Suppose $\sf{A}$ is a subcategory of an abelian category, $\sf{B}$ is a Grothendieck category with generator $\mathcal{O}$, and suppose 
$$
F:{\sf{B}}\times \sf{A} \rightarrow {\sf{B}}
$$ 
and 
$$
G:{\sf{A}}^{op}\times \sf{B} \rightarrow \sf{B}
$$ 
are conjugate bifunctors.  If $F(\mathcal{O},-)$ is exact, $G(-,\mathcal{E})$ is left exact for $\mathcal{E}$ in $\sf{B}$, and $G(-,\mathcal{E})$ is exact for $\mathcal{E}$ injective in $\sf{B}$.   
\end{lemma}

\begin{proof}
Let 
\begin{equation} \label{eqn.newexact}
0 \rightarrow \mathcal{H} \rightarrow \mathcal{I} \rightarrow \mathcal{J} \rightarrow 0
\end{equation}
be exact in $\sf{A}$.  Applying $G(-,\mathcal{E})$ to this sequence, we get a sequence
\begin{equation} \label{eqn.exactg}
\mathcal{K} \rightarrow G(\mathcal{J},\mathcal{E}) \rightarrow G(\mathcal{I},\mathcal{E}) \rightarrow G(\mathcal{H},\mathcal{E}) \rightarrow \mathcal{C}
\end{equation}
where $\mathcal{K}$ is the kernel of the leftmost arrow in (\ref{eqn.exactg}) and $\mathcal{C}$ is the cokernel of the rightmost arrow in (\ref{eqn.exactg}).  Since $\sf{B}$ is a Grothendieck category, there exists a quotient of $\sf{Mod } \Hom_{\sf{B}}(\mathcal{O},\mathcal{O})$ (with quotient functor $\pi$) such that $\pi \Hom_{\sf{B}}(\mathcal{O},-)$ is an equivalence \cite[Theorem 13.5, p. 88]{smith2}.  If we apply the functor $P= \pi \Hom_{\sf{B}}(\mathcal{O},-)$ to (\ref{eqn.exactg}), we get a sequence
\begin{equation} \label{eqn.exactg2}
P(\mathcal{K}) \rightarrow P(G(\mathcal{J},\mathcal{E})) \rightarrow P(G(\mathcal{I},\mathcal{E})) \rightarrow P(G(\mathcal{H},\mathcal{E})) \rightarrow P(\mathcal{C}).
\end{equation}
We will show that $P(\mathcal{K})=P(\mathcal{C})=0$ and that (\ref{eqn.exactg2}) is exact in the middle.  Since $P$ is an equivalence, the result will follow.  
Since $F$ and $G$ are conjugate bifunctors, the diagram
\begin{equation}
\begin{CD}
\Hom_{\sf{B}}(\mathcal{O},G(\mathcal{J},\mathcal{E}))  @>>> \Hom_{\sf{B}}(\mathcal{O},G(\mathcal{I},\mathcal{E})) @>>> \Hom_{\sf{B}}(\mathcal{O},G(\mathcal{H},\mathcal{E}))  \\
 @VVV @VVV	@VVV \\
\Hom_{\sf{B}}(F(\mathcal{O},\mathcal{J}), \mathcal{E})  @>>> \Hom_{\sf{B}}(F(\mathcal{O},\mathcal{I}),\mathcal{E}) @>>> \Hom_{\sf{B}}(F(\mathcal{O},\mathcal{K}), \mathcal{E}).   
\end{CD}
\end{equation}
with vertical arrows adjointness isomorphisms and horizontal arrows induced by (\ref{eqn.newexact}), commutes.  Since $F(\mathcal{O},-)$ is exact, the kernel of the leftmost map on the bottom row is $0$ and the bottom row is exact in the middle.  If $\mathcal{E}$ is injective, the rightmost bottom map is surjective as well.  Thus the kernel of the leftmost map on the top row is $0$ and the top row is exact in the middle. If  $\mathcal{E}$ is injective, the rightmost top map is surjective as well.
\end{proof}

\subsection{Elementary properties of internal Hom}

\begin{proposition} \label{prop.homisom}
The composition
\begin{equation} \label{eqn.homisom}
\M \overset{\eta_{\M}}{\rightarrow} \HU(\A,\M \Ten \A) \rightarrow \HU(\A,\M)
\end{equation}
whose right-most map is induced by the multiplication map $\M \Ten \A \rightarrow \A$ is a natural isomorphism in ${\sf Gr} \A$.
\end{proposition}

\begin{proof}
By Yoneda's Lemma, if $\alpha:\Homa(-,\mathcal{N}) \rightarrow \Homa(-, \mathcal{N}')$ is an isomorphism of functors, $\alpha_{\N}(\operatorname{id}):\N \rightarrow \N'$ is an isomorphism.  

We construct an isomorphism $\alpha:  \Homa(-, \M) \rightarrow \HU(-,\HomA(\A,\M))$.  By Proposition \ref{prop.multiso}, multiplication induces an isomorphism of functors $- \Ten \A \rightarrow \operatorname{id}_{\sf{Gr }\A}$, which induces an isomorphism $\Homa(-, \M) \rightarrow \Homa(-\Ten \A, \M)$.  Composing this map with the adjoint isomorphism gives the desired isomorphism $\alpha$:
$$
\Homa(-, \M) \rightarrow \Homa(- \Ten \A, \M) \rightarrow \Homa(-, \HU(\A,\M)).
$$
Thus, the image of the identity $\operatorname{id}_{\M}:\M \rightarrow \M$ under $\alpha$ gives an isomorphism $\M \rightarrow \HU(\A,\M)$.  Computing explicitly, we find $\alpha(\operatorname{id}_{\M})$ is the morphism (\ref{eqn.homisom}).
\end{proof}

\begin{theorem} \label{theorem.hom}
Let $\mathcal{M}$ be an object in ${\sf{Gr }}\mathcal{A}$ and let $\mathcal{C}$ be an object in $\mathbb{B}$.
\begin{enumerate}

\item{} 
$\HU (-,-)$ is left exact in either factor.

\item{}
If $\mathcal{M}$ is an injective object in ${\sf{Gr}} \mathcal{A}$, $\HU (-,\mathcal{M}):{\mathbb{B}}^{op} \rightarrow {\sf{Gr}\mathcal{A}}$ is exact.

\item{}
$\HomA (-,-):{{\mathbb{G}}{\sf r} \mathcal{A}}^{op} \times {\sf{Gr}} \mathcal{A} \rightarrow {\sf{Mod }}{\mathcal{O}}_{X}$ is a bifunctor which is left exact in either factor.  Furthermore, if $\mathcal{M}$ is an injective object in ${\sf{Gr}} \mathcal{A}$, $\HomA (-,\mathcal{M}):{{\mathbb{G}}{\sf r} \mathcal{A}}^{op} \rightarrow {\sf{Qcoh }}{\mathcal{O}}_{X}$ is exact.

\item{}
If $\mathcal{Q}$ is a coherent, locally free ${\mathcal{O}}_{X}$-bimodule, there is a natural isomorphism
$$
\HomA (\mathcal{Q}\otimes e_{m}\mathcal{A},\mathcal{M}) \cong {\mathcal{M}}_{m} \otimes {\mathcal{Q}}^{*}.
$$
In particular, $\HomA (\mathcal{Q}\otimes e_{m}\mathcal{A},-)$ is exact.
\end{enumerate}
\end{theorem}

\begin{proof} 
The functor $\HU(\C,-)$ is left exact since, by Proposition \ref{prop.tensor} (1), it has a left adjoint.  For $\M$ an object of $\sf{Gr }\A$, it remains to show $\HU(-,\M)$ is left exact.  To this end, we note that 
$$
\underset{i,m \in \mathbb{Z}}{\oplus} ({\mathcal{O}}_{X}(i) \otimes e_{m}\mathcal{A})
$$ 
is a generator of $\sf{Gr }\A$ by Lemma \ref{lem.noethno}.  In addition, 
\begin{align*}
\underset{i,m \in \mathbb{Z}}{\oplus} ({\mathcal{O}}_{X}(i) \otimes e_{m}\mathcal{A}) \Ten - & \cong \underset{i,m \in \mathbb{Z}}{\oplus} ({\mathcal{O}}_{X}(i) \otimes  (e_{m}\mathcal{A} \Ten -)) \\
& \cong \underset{i,m \in \mathbb{Z}}{\oplus} ({\mathcal{O}}_{X}(i) \otimes  e_{m}(-))
\end{align*}
where the first isomorphism is a consequence of Lemma \ref{lemma.asso} and the last isomorphism is a consequence of Proposition \ref{prop.multiso}.  Thus, by Proposition \ref{prop.tensor} (2), the hypothesis of Lemma \ref{lem.general} hold and $\HU(-,\M)$ is left exact.  Thus, the first part of the theorem holds.  A similar argument, with $\M$ injective, proves that the second part of the theorem holds.  The third part of the theorem follows from the first two parts of the theorem since taking the $k$th graded component of an object in $\sf{Gr }\A$ is an exact functor.

We now prove the fourth part of the theorem.  We note that the $\mathcal{O}_{X}$-module $\HomA (\mathcal{Q}\otimes e_{m}\mathcal{A},\mathcal{M})$ is the equalizer of the diagram
\begin{equation} \label{eqn.hom8}
\begin{CD}
\underset{i}{\Pi}{\mathcal{M}}_{i} \otimes (\Q \otimes {\mathcal{A}}_{mi})^{*} @>>> \underset{j}{\Pi}{\mathcal{M}}_{j} \otimes (\Q \otimes {\mathcal{A}}_{mj})^{*}\\
@VVV		@VVV\\
\underset{j}{\Pi}(\underset{i}{\Pi}({\mathcal{M}}_{j} \otimes {\mathcal{A}}_{ij}^{*}) \otimes (\Q \otimes {\mathcal{A}}_{mi})^{*}) @>>> \underset{j}{\Pi}(\underset{i}{\Pi}{\mathcal{M}}_{j} \otimes (\Q \otimes ({\mathcal{A}}_{mi} \otimes {\mathcal{A}}_{ij}))^{*}) 
\end{CD}
\end{equation}
whose arrows are defined as in (\ref{eqn.homdef}).  We claim this diagram is isomorphic to the diagram
\begin{equation} \label{eqn.hom9}
\begin{CD}
\underset{i}{\Pi}({\mathcal{M}}_{i} \otimes {\mathcal{A}}_{mi}^{*}) \otimes \Q^{*} @>>> \underset{j}{\Pi}({\mathcal{M}}_{j} \otimes {\mathcal{A}}_{mj}^{*}) \otimes \Q^{*} \\
@VVV		@VVV\\
\underset{j}{\Pi}(\underset{i}{\Pi}(({\mathcal{M}}_{j} \otimes {\mathcal{A}}_{ij}^{*}) \otimes {\mathcal{A}}_{mi}^{*}) \otimes \Q^{*}) @>>> \underset{j}{\Pi}(\underset{i}{\Pi}({\mathcal{M}}_{j} \otimes ({\mathcal{A}}_{mi} \otimes {\mathcal{A}}_{ij})^{*})\otimes \Q^{*})
\end{CD}
\end{equation}
whose arrows are tensor products of the arrows in (\ref{eqn.homer}) with $\Q^{*}$.  To prove the claim, we think of these diagrams as the top and bottom, respectively, of a cube whose vertical edges are canonical isomorphisms (\ref{eqn.canonicalisom0}).  We orient the cube so that the bottom horizontals of (\ref{eqn.hom8}) and (\ref{eqn.hom9}) are the top and bottom horizontals of the facing side of the cube.  The far face of this cube obviously commutes.  The left face commutes by functoriality of the tensor product.

To show the closest face commutes, it suffices to show that the diagram
$$
\begin{CD}
\A_{ij}^{*} \otimes (\Q \otimes \A_{mi})^{*} @>>> (\Q \otimes (\A_{mi} \otimes \A_{ij}))^{*} \\
@VVV @VVV \\
\A_{ij}^{*} \otimes (\A_{mi}^{*} \otimes \Q^{*}) @>>> (\A_{mi} \otimes \A_{ij})^{*} \otimes \Q^{*}
\end{CD}
$$
whose maps are induced by the canonical isomorphisms (\ref{eqn.canonicalisom0}), commutes.  This follows from Corollary \ref{corollary.tweed3}.

Finally, to show the right face commutes, it suffices to show that the diagram
$$
\begin{CD}
(\Q \otimes \A_{mj})^{*} @>>> \A_{mj}^{*} \otimes \Q^{*} \\
@VVV @VVV \\
(\Q \otimes (\A_{mi} \otimes \A_{ij}))^{*} @>>> (\A_{mi} \otimes \A_{ij})^{*} \otimes \Q^{*}
\end{CD}
$$
whose arrows are induced by the canonical isomorphisms (\ref{eqn.canonicalisom0}), commutes.  This follows from Corollary \ref{cor.tweed}, and the claim follows.  

The claim, together with the fact that $- \otimes \Q^{*}$ commutes with products (\cite[Theorem 1, p. 118]{mac}), implies that $\HomA (\mathcal{Q}\otimes e_{m}\mathcal{A},\mathcal{M}) \cong \HomA (e_{m}\mathcal{A},\mathcal{M}) \otimes \Q^{*}$.  The result now follows from Proposition \ref{prop.homisom}.
\end{proof}

\subsection{Torsion}
In this section, we write the torsion functor $\tau:{\sf Gr}\A \rightarrow {\sf Gr }\A$ in terms of the internal Hom functor.

\begin{lemma} \label{lemma.tor1}
If $\mathcal{N}$ is an object of ${\sf Gr}\mathcal{A}$ such that ${\mathcal{N}}_{i}=0$ for all $i \geq m$, then $\dlim \HU(\A_{\geq n}, \mathcal{N})=0$.
\end{lemma}

\begin{proof}
Let $k \in \mathbb{Z}$.  We claim $\HU(\A_{\geq n},\N)_{k}=0$ whenever $n \geq m-k$.  To prove the claim, we note that
$$
\HU(\A_{\geq n},\N)_{k} \subset \underset{i<m}{\Pi}\N_{i} \otimes (\A_{\geq n})_{ki}^{*}.
$$
On the other hand,
$$
(A_{\geq n})_{ki}=
\begin{cases} {\mathcal{A}}_{ki} & \mbox{if $i \geq k+n$,} \\ 
0 & \mbox{otherwise.}  
\end{cases}
$$
Thus, $\HU(\A_{\geq n},\N)_{k}=0$ whenever $k+n \geq m$ as desired.
\end{proof}

\begin{lemma} \label{lemma.tor2}
If $\M$ is an object of $\sf{Gr }\A$, $\HU(\A/\A_{\geq n}, \M)$ is a direct limit of torsion submodules. 
\end{lemma}

\begin{proof}
We claim the multiplication map
$$
\mu_{k, k+n}:\HU(\A/\A_{\geq n},\M)_{k} \otimes \A_{k,k+n} \rightarrow \HU(\A/\A_{\geq n},\M)_{k+n}
$$
is the zero map for all $k$.  This will imply that $\HU(\A/\A_{\geq n},\M)^{k}$ is torsion.  The lemma will then follow from Lemma \ref{lemma.torneeds}.

To prove the claim, recall that the diagram
$$
\begin{CD}
\HU(\A/\A_{\geq n},\M)_{k} \otimes \A_{k,k+n} @>{\mu}>> \HU(\A/\A_{\geq n},\M)_{k+n} \\
@VVV @VVV \\
(\underset{i}{\Pi}\M_{i} \otimes (\A/\A_{\geq n})_{ki}^{*}) \otimes \A_{k,k+n} @>>> \underset{i}{\Pi}\M_{i} \otimes (\A/\A_{\geq n})_{k+n,i}^{*}
\end{CD}
$$
whose verticals are inclusions, and whose bottom map is induced by the multiplication map
\begin{equation} \label{eqn.mullt}
\A_{k,k+n}\otimes (\A/\A_{\geq n})_{k+n,i} \rightarrow (\A/\A_{\geq n})_{k,i}.
\end{equation}
commutes by definition of the right $\A$-module structure on $\HU(\A/\A_{\geq n},\M)$.  If $i<k+n$, the domain of (\ref{eqn.mullt}) is $0$.  On the other hand, if $i \geq k+n$, the codomain of (\ref{eqn.mullt}) is $0$.  The claim follows.
\end{proof}

\begin{proposition} \label{prop.tau}
There is a natural equivalence
$$
\tau \cong \underset{n \to \infty}{\lim} \HU (\A/\A_{\geq n}, -).
$$
\end{proposition}

\begin{proof}
Let $\mathcal{M}$ be an object in $\sf{Gr } \mathcal{A}$.  We prove that there is a natural isomorphism
$$
\tau \mathcal{M} \cong \underset{n \to \infty}{\lim} \HU (\A/\A_{\geq n}, \mathcal{M}).
$$
To this end, we first apply the functor $\HU (-, \mathcal{M})$ to the exact sequence 
$$
0 \rightarrow \A_{\geq n} \rightarrow \A \rightarrow \A/\A_{\geq n} \rightarrow 0
$$
in the category $\mathbb{B}$.  Since $\HU (-, \mathcal{M})$ is left exact, we get an exact sequence
$$
0 \rightarrow \HU (\A/\A_{\geq n},\mathcal{M}) \rightarrow \HU (\A,\mathcal{M}) \rightarrow \HU (\A_{\geq n},\mathcal{M}).
$$
Since ${\sf Gr }\mathcal{A}$ is a Grothendieck category, the induced sequence
$$
0 \rightarrow \dlim \HU (\A/\A_{\geq n},\mathcal{M}) \rightarrow \dlim \HU (\A,\mathcal{M})\overset{\psi}{\rightarrow} \dlim \HU (\A_{\geq n},\mathcal{M}).
$$
is exact.  Let $i$ be the composition $\tau \M \rightarrow \M \overset{\gamma}{\rightarrow} \dlim \HU (\A, \M)$, where $\gamma$ is induced by (\ref{eqn.homisom}).  We claim ($\tau \M, i$) is a kernel of $\psi$.  This will induce an isomorphism of kernels $\tau \M \rightarrow \dlim \HU (\A/\A_{\geq n}, \M)$, and naturality of the isomorphism follows in a straightforward way from the universal property of kernels.

To prove the claim, we first prove that $\psi i = 0$.  Suppose $\mathcal{N} \subset \mathcal{M}$ has the property that there exists an $m$ such that ${\mathcal{N}}_{i}=0$ for all $i \geq m$ and consider the commutative diagram
$$
\begin{CD}
\mathcal{N} @>\cong>> \dlim \HU (\mathcal{A}, \mathcal{N}) @>>> \dlim \HU (\A_{\geq n}, \mathcal{N}) \\
@VVV @VVV @VVV \\
\mathcal{M} @>>{\gamma}> \dlim \HU (\mathcal{A}, \mathcal{M}) @>>{\psi}> \dlim \HU (\A_{\geq n}, \mathcal{M}) 
\end{CD}
$$
whose verticals are induced by inclusion, whose right horizontals are induced by inclusion $\A_{\geq n} \rightarrow \A$ and whose upper left horizontal is induced by (\ref{eqn.homisom}).  By Lemma \ref{lemma.tor1}, $\dlim \HU (\A_{\geq n}, \mathcal{N})=0$.  Thus, the bottom route of the outer circuit is $0$.  Since $\tau \mathcal{M}$ is the direct limit of such $\mathcal{N}$, the composition
\begin{equation} \label{eqn.tauya}
\tau\mathcal{M} \rightarrow \mathcal{M} \overset{\gamma}{\rightarrow} \dlim \HU (\mathcal{A}, \mathcal{M}) \overset{\psi}{\rightarrow} \dlim \HU (\A_{\geq n}, \mathcal{M}) 
\end{equation}
is zero as well.  

Next, suppose
$$
\mathcal{N} \overset{\phi}{\rightarrow} \dlim \HU (\mathcal{A}, \mathcal{M}) \overset{\psi}{\rightarrow} \dlim \HU (\A_{\geq n}, \mathcal{M}) 
$$
equals $0$.  To complete the proof that $(\tau \M, i)$ is the kernel of $\psi$, we must show that there exists a unique map $\phi':\mathcal{N} \rightarrow \tau \mathcal{M}$ such that the diagram
$$
\begin{CD}
\mathcal{N} @>{\phi'}>> \tau \mathcal{M} \\
@V{\phi}VV @VVV \\
\dlim \HU (\mathcal{A}, \mathcal{M}) @>>{=}> \dlim \HU (\mathcal{A}, \mathcal{M}) 
\end{CD}
$$
commutes.  Since the diagram
$$
\begin{CD}
\operatorname{im}(\gamma^{-1}\phi) @>>> \dlim \HU (\mathcal{A}, \mbox{im}(\gamma^{-1}\phi)) @>\delta>> \dlim \HU (\A_{\geq n}, \mbox{im}(\gamma^{-1}\phi))\\
@VVV @VVV @VVV \\
\M @>>{\gamma}> \dlim \HU (\mathcal{A}, \mathcal{M}) @>>{\psi}> \dlim \HU (\A_{\geq n}, \mathcal{M})
\end{CD}
$$
whose verticals are monomorphisms induced by inclusion and whose upper left horizontal is induced by (\ref{eqn.homisom}), commutes, $\delta=0$.  In particular, the quotient map $\mathcal{A} \rightarrow \mathcal{A}/\A_{\geq n} \rightarrow 0$ induces an isomorphism
$$
\dlim \HU(\mathcal{A}/\A_{\geq n}, \mbox{im }(\gamma^{-1}\phi)) \rightarrow \dlim \HU(\mathcal{A}, \mbox{im }(\gamma^{-1}\phi)) \rightarrow \mbox{im }(\gamma^{-1}\phi)
$$
whose right arrow is induced by (\ref{eqn.homisom}).  By Lemma \ref{lemma.tor2}, 
$$
\mbox{im }(\gamma^{-1} \phi) \cong \dlim \HU(\mathcal{A}/\A_{\geq n}, \mbox{im }(\gamma^{-1}\phi))
$$ 
is a direct limit of torsion submodules.  Thus, the inclusion $\mbox{im }(\gamma^{-1} \phi) \rightarrow \mathcal{M}$ factors through $\tau \mathcal{M}$, and hence, $\gamma^{-1} \phi$ factors through $\tau \mathcal{M}$ via the left-most vertical in the commutative diagram

$$
\begin{CD}
\mathcal{N} @>>> \mbox{im }(\gamma^{-1}\phi) @>\cong>> \dlim \HU(\mathcal{A}, \mbox{im }(\gamma^{-1}\phi)) \\
@V{\phi'}VV @VVV @VVV \\
\tau \mathcal{M} @>>> \mathcal{M} @>>{\gamma}> \dlim \HU(\mathcal{A}, \mathcal{M}). 
\end{CD}
$$
whose right two verticals are induced by inclusion and whose top right horizontal is induced by (\ref{eqn.homisom}).  To show $\phi'$ is unique, suppose $\zeta:\mathcal{N} \rightarrow \tau\mathcal{M}$ is another map making the left square commute.  Since $\tau \mathcal{M} \rightarrow \mathcal{M}$ is a monomorphism, $\phi'=\zeta$.
\end{proof}

\vfill
\eject

%% file: cohomology.tex
\section{Cohomology}

In light of Theorem \ref{theorem.hom}, the right derived functors of $\HomA (\mathcal{C},-)$ exist.  Let ${\Ext}^{i}(\mathcal{C},-)$ denote the functor ${\operatorname{R}}^{i}\HomA (\mathcal{C},-)$.

\subsection{Elementary properties of internal Ext}

\begin{proposition} \label{prop.cohomology}
Suppose there is an exact sequence 
$$
\ldots \rightarrow {\mathcal{L}}_{1} \rightarrow {\mathcal{L}}_{0} \rightarrow \mathcal{C} \rightarrow 0
$$
in ${\mathbb{G}}\sf{r} \mathcal{A}$ such that $\mathcal{L}_{i}$ is $\HomA(-,\M)$ acyclic for some $\mathcal{M}$ in $\sf{Gr} \mathcal{A}$.  Then there is a natural isomorphism
$$
\Ext^{i}(\mathcal{C},\mathcal{M}) \cong h^{i}(\HomA (\mathcal{L}{\bf .},\mathcal{M})).
$$  
\end{proposition}

\begin{proof}
By Theorem \ref{theorem.hom}, the hypothesis of \cite[Proposition 8.2, p. 809]{lang} are satisfied.  The result follows.
\end{proof}

\begin{lemma} \label{lemma.cohomology}
Let 
$$
0 \rightarrow \C' \rightarrow \C \rightarrow \C'' \rightarrow 0
$$
be a short exact sequence in $\mathbb{B}$.  Then, for fixed $\M$ in $\sf{Gr }\A$, we have a long exact sequence
\begin{align*}
0 & \rightarrow \HU(\C'',\M) \rightarrow \HU(\C,\M) \rightarrow \HU(\C',\M) \rightarrow \\
& \rightarrow \ExtU^{1}(\C'',\M) \rightarrow \ExtU^{1}(\C,\M) \rightarrow \ExtU^{1}(\C',\M) \rightarrow \cdots
\end{align*}
such that the association
$$
\mathcal{D} \longmapsto \ExtU^{n}(\D,\M)
$$
is a $\delta$-functor.
\end{lemma}

\begin{proof}
By Theorem \ref{theorem.hom}, the hypothesis of \cite[Lemma 8.3, p.809]{lang} are satisfied.
\end{proof}
Notice that we are not claiming $\delta$ is universal.  We shall use Lemma \ref{lemma.cohomology} without comment in the sequel.

\begin{lemma} \label{lemma.delta}
If $\Q$ is a coherent, locally free $\mathcal{O}_{X}$-bimodule, there is a natural isomorphism 
$$
\Ext^{i}(\mathcal{Q} \otimes \mathcal{C}, \mathcal{M}) \cong \Ext^{i}(\mathcal{C},\mathcal{M})\otimes {\mathcal{Q}}^{*}
$$
for all $i \geq 0$.
\end{lemma}

\begin{proof}
The $i=0$ case of the isomorphism follows from Theorem \ref{theorem.hom} (4).

Since any $\delta$-functor composed with an exact functor is a $\delta$-functor, all terms are $\delta$-functors.  For $\M$ injective and $i>0$, all terms vanish by Theorem \ref{theorem.hom} (2), and the result follows from \cite[Theorem 1.3A, p. 206]{hartshorne}.
\end{proof}
For the remainder of this paper, we assume $\A_{i,i+1}$ is locally free of rank $2$.  

\subsection{$\A$ is Gorenstein}
The next result tells us that $\A$ is Gorenstein.
\begin{theorem} \label{theorem.gor} 
Let $\mathcal{L}$ be a coherent, locally free $\mathcal{O}_{X}$-module.  Then 
$$
\ExtU^{i}(\A_{0},\mathcal{L} \otimes e_{l}\A)=0 \mbox{ for $i \neq 2$}
$$
and
$$
\ExtU^{2}(\A_{0},\mathcal{L} \otimes e_{l}\A)_{j}  \cong 
\begin{cases}
\mathcal{L} \otimes \Q_{l-2}^{*}& \text{if $j = l -2$}, \\
0& \text{otherwise}.
\end{cases}
$$
\end{theorem}

\begin{proof}
The $k$th graded piece of $\ExtU^{i}(\A_{0},\mathcal{L} \otimes e_{l}\A)$ is $\Ext^{i}(e_{k}\A_{0},\mathcal{L} \otimes e_{l}\A)$.  By \cite[Theorem 7.1.2]{p1bundles}, the sequence
\begin{equation} \label{eqn.sequence1}
0 \rightarrow \Q_{k} \otimes e_{k+2}\A \rightarrow \mathcal{E}_{k} \otimes e_{k+1}\A \overset{\phi}{\rightarrow} e_{k}\A \rightarrow e_{k}\A/e_{k}\A_{\geq 1} \rightarrow 0
\end{equation}
where $\phi$ is multiplication is exact.  In addition, by Theorem \ref{theorem.hom}(4), the sequence is a $\HomA(-,\mathcal{L} \otimes e_{l}\A)$-acyclic resolution of $e_{k}\A/e_{k}\A_{\geq 1}$.  Applying the functor $\HomA(-,\mathcal{L} \otimes e_{l}\A)$ to the sequence obtained by omitting $e_{k}\A/e_{k}\A_{\geq 1}$ from (\ref{eqn.sequence1}) gives a sequence
\begin{equation} \label{eqn.sequence2}
\HomA(e_{k}\A,\mathcal{L}  \otimes e_{l}\A)  \overset{d_{0}}{\rightarrow} \HomA(\mathcal{E}_{k} \otimes  e_{k+1}\A,\mathcal{L} \otimes e_{l}\A) 
\end{equation}
$$
\overset{d_{1}}{\rightarrow} \HomA(\mathcal{Q}_{k}\otimes  e_{k+2}\A, \mathcal{L} \otimes e_{l}\A).
$$
By Proposition \ref{prop.cohomology}, 
$$
\ExtU^{0}(\A_{0},\mathcal{L} \otimes e_{l}\A)_{k} \cong \operatorname{ker }d_{0},
$$ 
$$\ExtU^{1}(\A_{0}, \mathcal{L} \otimes e_{l}\A)_{k} \cong \frac{\operatorname{ker }d_{1}}{\operatorname{im }d_{0}},
$$ 
$$
\ExtU^{2}(\A_{0}, \mathcal{L} \otimes e_{l}\A)_{k} \cong \frac{\HomA(\Q_{k} \otimes e_{k+2}\A, \mathcal{L} \otimes e_{l}\A)}{\operatorname{im }d_{1}}
$$ 
and $\ExtU^{i}(\A_{0},\mathcal{L} \otimes e_{l}\A) = 0$ for all $i > 2$.  The terms of (\ref{eqn.sequence2}) are isomorphic to $\mathcal{L} \otimes \A_{lk}$, $\mathcal{L} \otimes \A_{l,k+1} \otimes \mathcal{E}_{k}^{*}$ and $\mathcal{L} \otimes \A_{l,k+2} \otimes \Q_{k}^{*}$ respectively.  Thus, by Theorem \ref{theorem.hom} (4), all terms vanish if $k < l-2$.  If $k=l-2$, the left and center terms are zero, and we conclude that
$$
\ExtU^{0}(\A_{0},\mathcal{L} \otimes e_{l}\A)_{l-2}=\ExtU^{1}(\A_{0},\mathcal{L} \otimes e_{l}\A)_{l-2}=0
$$ 
and
$$
\ExtU^{2}(\A_{0},\mathcal{L} \otimes e_{l}\A)_{l-2} \cong \mathcal{L} \otimes \A_{l,l} \otimes \Q_{l-2}^{*}.
$$  
We prove that (\ref{eqn.sequence2}) is exact for $k > l-2$, which will establish the result.  We first show (\ref{eqn.sequence2}) is exact in the middle.  Since the sequence 
$$
0 \rightarrow \Q_{k} \otimes e_{k+2}\A \rightarrow \mathcal{E}_{k} \otimes e_{k+1}\A \overset{\phi'}{\rightarrow} e_{k}\A_{\geq 1} \rightarrow 0
$$
is exact, where $\phi'$ is the restriction of $\phi$ to degrees $\geq 1$, and $\HomA(-,\mathcal{L} \otimes e_{l}\A)$ is left exact by Theorem \ref{theorem.hom}(1), this sequence induces an exact sequence
$$
0 \rightarrow \HomA(e_{k}\A_{\geq 1},\mathcal{L} \otimes e_{l}\A) \overset{d_{0}'}{\rightarrow} \HomA(\mathcal{E}_{k} \otimes e_{k+1}\A,\mathcal{L} \otimes  e_{l}\A) 
$$
$$
\overset{d_{1}}{\rightarrow} \HomA(\mathcal{Q}_{k} \otimes e_{k+2}\A, \mathcal{L} \otimes e_{l}\A).
$$
To prove (\ref{eqn.sequence2}) is exact in the middle, it suffices to show $\operatorname{im }d_{0}=\operatorname{im }d_{0'}$.  To prove this, it suffices to prove the morphism $\gamma:\HU(\A,\mathcal{L} \otimes e_{l}\A) \rightarrow \HU(\A_{\geq 1},\mathcal{L} \otimes e_{l}\A)$ induced by the inclusion $\A_{\geq 1} \rightarrow \A$ is epi.  

We proceed to prove that $\gamma$ is an epi.  By Proposition \ref{lem.inclusion} and the adjointness of $\HU$ and $\Ten$, we may conclude that, for all $\N$ in $\sf{Gr }\A$, the map induced by $\gamma$,
$$
f:\Homa(\N,\HU(\A,\mathcal{L} \otimes e_{l}\A)) \rightarrow \Homa(\N,\HU(\A_{\geq 1}, \mathcal{L} \otimes e_{l}\A))
$$
is epi.  If $\N$ is a generator of $\sf{Gr }\A$ and $\pi: {\sf Mod} \Homa(\N,\N) \rightarrow \sf{B}$ is a quotient functor which makes $\pi \Homa(\N,-):{\sf Gr}\A \rightarrow \sf{B}$ an equivalence, $\pi(f)=\pi \Homa(\N,\gamma)$ is an epi.  Thus $\gamma$ is epi and (\ref{eqn.sequence2}) is exact in the middle.

We next show that $d_{0}$ is a monomorphism.  This is the same as showing that $\HomA(e_{k}\A/A_{\geq 1},\mathcal{L} \otimes e_{l}\A)=0$, which is the same as showing that the equalizer of the diagram
\begin{equation} \label{eqn.route}
\begin{CD}
\mathcal{L} \otimes \A_{lk} \otimes (\A/\A_{\geq 1})_{kk}^{*} @>>> \underset{j}{\Pi}\mathcal{L} \otimes \A_{l,j+k} \otimes (\A/\A_{\geq 1})_{k,j+k}^{*} \\
@VVV @VVV \\
\underset{j}{\Pi}(\mathcal{L} \otimes \A_{l,j+k} \otimes \A_{k,j+k}^{*}) \otimes (\A/\A_{\geq 1})_{kk}^{*} @>>{\cong}> \underset{j}{\Pi}(\mathcal{L} \otimes \A_{lj} \otimes ((\A/\A_{\geq 1})_{kk} \otimes \A_{kj})^{*}
\end{CD}
\end{equation}
whose maps are those in (\ref{eqn.homdef}), equals $0$.  Now, the morphism
$$
\A_{lk} \rightarrow \A_{lk} \otimes \A_{k,k+1} \otimes \A_{k,k+1}^{*} \rightarrow \A_{l,k+1} \otimes \A_{k,k+1}^{*}
$$
is monic by Lemma \ref{lemma.monomult}, i.e. $j=1$ component of the left-hand route of (\ref{eqn.route}) is a monomorphism, while the $j=1$ component of the right-hand route of (\ref{eqn.route}) is the zero morphism, whence the claim.

We defer the proof that $d_{1}$ is epi to Section 7.5.
\end{proof}

\begin{definition}
Let $\C$ be an object of $\mathbb{B}$.  Let $l, r, m \in \mathbb{Z}$.  We say $\C \subset [l,r]$ if $\C_{ij}$ is zero whenever $i,j \in \mathbb{Z}$ do not satisfy $l \leq j-i \leq r$.  We say $\C$ is concentrated in degree $m$ if $\C \subset [m,m]$. 
\end{definition}

\begin{corollary} \label{cor.vanish}
If $\C$ is an object of ${\mathbb{G} \sf{r}} \A$ concentrated in degree $m$ and $\mathcal{L}$ is a coherent, locally free $\mathcal{O}_{X}$-module and $\M$ is an object in ${\sf Gr} \A$,
$$
\ExtU^{i}(\C,\M)_{j} \cong \ExtU^{i}(\A_{0},\M)_{j+m} \otimes \C_{j,j+m}^{*},
$$ 
$\ExtU^{i}(\C, \mathcal{L} \otimes e_{l}\A) = 0$ for all integers $i \neq 2$ and $\ExtU^{i}(\C, \M)=0$ for all $i>2$.
\end{corollary}

\begin{proof}
By hypothesis, $\C = \underset{k}{\oplus}\C_{k,k+m}$.  Thus
\begin{align*}
\ExtU^{i}(\C,\M)_{j} & = \ExtU^{i}(\underset{k}{\oplus}\C_{k,k+m},\M)_{j} \\
& \cong \Ext^{i}(e_{j}\C,\M) \\
& \cong \Ext^{i}(\C_{j,j+m} \otimes e_{j+m}\A_{0},\M) \\
& \cong \ExtU^{i}(\A_{0},\M)_{j+m} \otimes \C_{j,j+m}^{*}
\end{align*}
where we have used Lemma \ref{lemma.delta} to establish the final isomorphism.  The second result follows from Theorem \ref{theorem.gor} and the third result follows from Proposition \ref{prop.cohomology}.
\end{proof}

\begin{corollary} \label{cor.vanish2}
If $\mathcal{L}$ is a coherent, locally free $\mathcal{O}_{X}$-module, $\ExtU^{i}(\A/\A_{\geq n}, \M) = 0$ for $i>2$ and $\ExtU^{i}(\A/\A_{\geq n}, \mathcal{L} \otimes e_{l}\A) = 0$ for $i \neq 2$.
\end{corollary}

\begin{proof}
We prove the result by induction on $n$.  When $n=1$, the result follows from Corollary \ref{cor.vanish}.  Next, we note that the exact sequence
$$
0 \rightarrow \A_{\geq n}/\A_{\geq n+1} \rightarrow \A/\A_{\geq n+1} \rightarrow \A/\A_{\geq n} \rightarrow 0
$$
induces a long exact sequence, of which
$$
\ExtU^{i}(\A/\A_{\geq n},\M) \rightarrow \ExtU^{i}(\A/\A_{\geq n+1},\M) \rightarrow \ExtU^{i}(\A_{\geq n}/\A_{\geq n+1},\M)
$$
is a part.  If $i>2$ the left term is zero by induction hypothesis while the right term is zero by Corollary \ref{cor.vanish}.  If $i=0$ or $i=1$ and $\M=\mathcal{L} \otimes e_{l}\A$ the same reasoning ensures that the center is zero.
\end{proof}

The following Lemma and proof are slight modifications of \cite[Lemma 2.3, p. 21]{intro}
\begin{lemma} \label{lem.correctgrade}
Let $\C$ be an object of $\mathbb{B}$ such that $\C \subset [l,r]$ and let $\M$ be an object of ${\sf{Gr }}\A$ such that 
$$
\ExtU^{j}(\A_{0}, \M) \subset [\lambda, \rho]
$$
for all $j \leq i$.  Then, for $j \leq i$,
$$
\ExtU^{j}(\C, \M) \subset [\lambda-r, \rho-l].
$$
\end{lemma}

\begin{proof}
Assume $\C$ is concentrated in degree $m$.  By Corollary \ref{cor.vanish}, 
\begin{equation} \label{eqn.isomya}
\ExtU^{i}(\C,\M)_{k} \cong \ExtU^{i}(\A_{0},\M)_{k+m} \otimes \C_{k,k+m}^{*}.
\end{equation}
Since $\ExtU^{j}(\A_{0}, \M) \subset [\lambda, \rho]$, (\ref{eqn.isomya}) implies $\ExtU^{i}(\C,\M) \subset [\lambda-m,\rho-m]$ as desired.

Now, consider the general situation $\C \subset [l,r]$, and define $e_{k}\C'=\C_{k,k+r}$.  We have a short exact sequence
$$
0 \rightarrow \C' \rightarrow \C \rightarrow \C/\C' \rightarrow 0
$$
with $\C/\C' \subset [l,r-1]$.  We shall use this for induction on $r-l$, the case $r=l$ being dealt with above.

Suppose that $r-l=m$ and that the result is correct when $r-l=m-1$.  The long-exact sequence for ${\ExtU}(-,\M)$ on the above short exact sequence contains
$$
\ExtU^{j}(\C/\C',\M) \rightarrow \ExtU^{j}(\C,\M) \rightarrow \ExtU^{j}(\C',\M).
$$
But the left hand module is in $[\lambda -r+1, \rho-l]$, whereas the right module is in $[\lambda-r, \rho-r]$.
\end{proof}

\begin{lemma} \label{lemma.summand}
Let $\mathcal{L}$ be a coherent, locally free $\mathcal{O}_{X}$-module.  In the category of graded $\mathcal{O}_{X}$-modules, $\ExtU^{2}(\A_{0},\mathcal{L} \otimes e_{l}\A)$ is a direct summand of 
$$
\underset{n \to \infty}{\lim} \ExtU^{2}(\quota,\mathcal{L} \otimes e_{l}\A),
$$
and
\begin{equation} \label{eqn.conjecture}
(\underset{n \to \infty}{\lim}\ExtU^{2}(\A/\A_{\geq n},\mathcal{L} \otimes e_{l}\A))_{l-2-i} \cong 
\begin{cases}
\mathcal{L} \otimes Q_{l-2}^{*} \otimes \A_{l-2-i,l-2}^{*}& \text{if $i \geq 0$}, \\
0& \text{otherwise}.
\end{cases}
\end{equation}
\end{lemma}

\begin{proof}
For all $n \geq 1$, the exact sequence
\begin{equation} \label{eqn.needed}
0 \rightarrow \A_{\geq n}/\A_{\geq n+1} \rightarrow \A/\A_{\geq n+1} \rightarrow \A/\A_{\geq n} \rightarrow 0
\end{equation}
induces a sequence
\begin{equation} \label{eqn.split}
0 \rightarrow  \ExtU^{2}(\A/\A_{\geq n},\mathcal{L} \otimes e_{l}\A) \rightarrow \ExtU^{2}(\A/\A_{\geq n+1}, \mathcal{L} \otimes e_{l}\A) 
\end{equation}
$$
\rightarrow \ExtU^{2}(\A_{\geq n}/\A_{\geq n+1}, \mathcal{L} \otimes e_{l}\A) \rightarrow 0
$$
of graded $\mathcal{O}_{X}$-modules which is exact by Corollary \ref{cor.vanish2}.  In order to prove the lemma, we prove, by induction on $n$, that $\ExtU^{2}(\A/\A_{n+1},\mathcal{L} \otimes e_{l}\A)$ is contained in degrees $[l-2-n,l-2]$.

When $n=1$, we have an exact sequence
$$
0 \rightarrow  \ExtU^{2}(\A/\A_{\geq 1}, \mathcal{L} \otimes e_{l}\A) \rightarrow \ExtU^{2}(\A/\A_{\geq 2},\mathcal{L} \otimes e_{l}\A) 
$$
$$
\rightarrow \ExtU^{2}(\A_{\geq 1}/\A_{\geq 2},\mathcal{L} \otimes e_{l}\A) \rightarrow 0.
$$
Since the left hand term is concentrated in degree $l-2$, while the right hand term is concentrated in degree $l-3$ by Lemma \ref{lem.correctgrade}, the central term is contained in $[l-3,l-2]=[l-2-1,l-2]$.  

Next, suppose $\ExtU^{2}(\A/A_{\geq n}, \mathcal{L} \otimes e_{l}\A) \subset [l-2-n+1,l-2]$.  Since the right hand term of the exact sequence (\ref{eqn.split}) is concentrated in degree $l-2-n$ by Lemma \ref{lem.correctgrade}, the central term is contained in $[l-2-n,l-2]$ as desired.  The isomorphism (\ref{eqn.conjecture}) now follows from Corollary \ref{cor.vanish}.
\end{proof}

\begin{corollary} \label{cor.cd}
$$
\operatorname{cd }\underset{n \to \infty}{\lim} \HU (\A/\A_{\geq n}, -)=2.
$$  
\end{corollary}

\begin{proof}
Since direct limits are exact in $\sf{Gr }\A$, 
$$
\operatorname{R}^{i}(\underset{n \to \infty}{\lim} \HU (\A/\A_{\geq n}, -)) \cong \underset{n \to \infty}{\lim} \ExtU^{i}(\quota,-).
$$  
By Corollary \ref{cor.vanish2}, this is zero whenever $i > 2$.   

On the other hand, by Lemma \ref{lemma.summand}, $\ExtU^{2}(\A_{0},\mathcal{L} \otimes e_{l}\A)$ is a direct summand of $\underset{n \to \infty}{\lim} \ExtU^{2}(\quota,\mathcal{L} \otimes e_{l}\A)$ in the category of $\mathcal{O}_{X}$-modules.  By Theorem \ref{theorem.gor}, 
$$
\ExtU^{2}(\A_{0},\mathcal{L} \otimes e_{l}\A) \neq 0
$$ 
whence the result.
\end{proof}

\subsection{Application to Non-commutative ${\mathbb{P}}^{1}$-bundles}
In this section, we prove a vanishing result for the cohomology of ${\sf Proj }\A$ and we compute the cohomological dimension of ${\sf Proj }\A$ when $X$ is a smooth projective variety of dimension $d$ over $K$. 

\begin{theorem} \label{theorem.exacts}
For $i \geq 1$, the right-derived functors of $\tau:{\sf Gr }\A \rightarrow {\sf Gr }\A$ and $\omega: {\sf Proj }\A \rightarrow {\sf Gr }\A$ satisfy
$$
\operatorname{R}^{i+1}\tau(-) \cong \operatorname{R}^{i}\omega(\pi(-))
$$
and there is an exact sequence
$$
\begin{CD}
0 @>>> \tau\M @>>> \M @>{\eta}>> \omega \pi \M @>>> \operatorname{R}^{1}\tau \M @>>> 0.
\end{CD}
$$
\end{theorem}

\begin{proof}
By Lemma \ref{lemma.smithuse}, the hypothesis of \cite[Theorem 14.15 (3), p.99]{smith2} hold with $\sf{A}={\sf Gr }\A$ and $\sf{T}={\sf Tors }\A$.  The result follows.
\end{proof}

\begin{corollary} \label{cor.exacts}
The counit
$$
\mathcal{O}_{X}(-i) \otimes e_{m} \A \rightarrow \omega \pi (\mathcal{O}(-i)_{X} \otimes e_{m}\A) 
$$
is an isomorphism for all $i,m \in \mathbb{Z}$.
\end{corollary}

\begin{proof}
By Theorem \ref{theorem.exacts}, it suffices to prove that $\tau(\mathcal{O}_{X}(-i) \otimes e_{m} \A) =0$ and $\operatorname{R}^{1}\tau (\mathcal{O}_{X}(-i) \otimes e_{m} \A) =0$.  Both of these equalities follow from Corollary \ref{cor.vanish2} in light of the fact that 
$$
\operatorname{R}^{i}\tau \cong \underset{n \to \infty}{\lim} \ExtU^{i} (\A/\A_{\geq n}, -).
$$
\end{proof}

\begin{theorem} \label{theorem.cohom}
Let $d$ be the cohomological dimension of $X$.  For $0 \leq j \leq d$, 
\begin{equation} \label{eqn.basedual1}
\operatorname{R}^{j}(\Gamma \circ \omega(-)_{0}) (\pi(\mathcal{O}_{X}(-i) \otimes e_{m}\A)) =0
\end{equation}
whenever $i, m >> 0$.
\end{theorem}

\begin{proof}
Since the counit $\mathcal{O}_{X}(-i) \otimes e_{m}\A \rightarrow \omega \pi (\mathcal{O}_{X}(-i) \otimes e_{m}\A)$ is an isomorphism by Corollary \ref{cor.exacts}, (\ref{eqn.basedual1}) holds when $j=0$ since $(\mathcal{O}_{X}(-i) \otimes e_{m}\A)_{0} = 0$ for $m>0$.  In case $0<j \leq d$, consider the exact sequence
\begin{equation} \label{eqn.exactgs}
\operatorname{R}^{j}\Gamma(\omega(\pi(\mathcal{O}_{X}(-i) \otimes e_{m}\A))_{0}) \rightarrow \operatorname{R}^{j}(\Gamma \omega(-)_{0})(\pi(\mathcal{O}_{X}(-i) \otimes e_{m}\A)) \rightarrow 
\end{equation}
$$
\operatorname{R}^{j-1}\Gamma \operatorname{R}^{1}(\omega(-)_{0})(\pi(\mathcal{O}_{X}(-i) \otimes e_{m}\A))
$$
arising from the Grothendieck spectral sequence.  The left hand term of (\ref{eqn.exactgs}) is zero since the counit $\mathcal{O}_{X}(-i) \otimes e_{m}\A \rightarrow \omega \pi (\mathcal{O}_{X}(-i) \otimes e_{m}\A)$ is an isomorphism by Corollary \ref{cor.exacts}.  Thus, to prove the central term of (\ref{eqn.exactgs}) is zero, it suffices to show the final term of (\ref{eqn.exactgs}) is zero.  Since $(\operatorname{R}^{1}\omega) \pi \cong \operatorname{R^{2}\tau}$ by Theorem \ref{theorem.exacts}, it in fact suffices to prove that    
$$
\operatorname{R}^{j-1}\Gamma((\operatorname{R}^{2}\tau)(\mathcal{O}_{X}(-i) \otimes e_{m}\A)_{0})=0
$$
for $i,m >> 0$.  By Lemma \ref{lemma.summand}, it thus suffices to prove that 
$$
\operatorname{R}^{j-1}\Gamma(\mathcal{O}_{X}(-i) \otimes \Q_{m-2}^{*}\otimes \A_{0,m-2}^{*}) = 0
$$ 
for $i >> 0$.  Let $\mathcal{F} =  \Q_{m-2}^{*}\otimes \A_{0,m-2}^{*}$.  Then, for $1 \leq j \leq d$, Serre duality on $X$ gives us an isomorphism
$$
\operatorname{Ext}^{j}_{X}(\mathcal{O}_{X}(-i) \otimes  \mathcal{F},\omega_{X}) \cong \operatorname{R}^{d-j}\Gamma(X, \mathcal{O}_{X}(-i) \otimes \mathcal{F})'
$$
where $\omega_{X}$ is a dualizing sheaf on $X$.  On the other hand, for $i>>0$,
\begin{align*}
\operatorname{Ext}^{j}_{X}(\mathcal{O}_{X}(-i) \otimes  \mathcal{F},\omega_{X}) & \cong \operatorname{Ext}^{j}_{X}(\mathcal{O}_{X}(-i),\omega_{X} \otimes \mathcal{F}^{*}) \\
& \cong \Gamma(X,\mathcal{E}{\it xt}_{\mathcal{O}_{X}}^{j}(\mathcal{O}_{X}(-i), \omega_{X} \otimes \mathcal{F}^{*})) \\
& = 0
\end{align*}
where the first isomorphism is due to a variant of \cite[Proposition 6.7, p. 235]{hartshorne}, the second isomorphism is a consequence of \cite[Proposition 6.9, p. 236]{hartshorne}, and the final equality is due to the fact that $\mathcal{O}_{X}(-i)$ is locally free.
\end{proof}

\begin{lemma} \label{lem.spectral}
Suppose ${\sf A}$, ${\sf B}$ and ${\sf C}$ are abelian categories such that $\sf A$ and $\sf B$ have enough injectives and $G:\sf{A} \rightarrow \sf{B}$ and $F:\sf{B} \rightarrow \sf{C}$ are left exact functors.  If, for $\mathcal{I}$ injective in $\sf A$, $G(\mathcal{I})$ is $F$-acyclic and if $\operatorname{cd }G$=1 and $\operatorname{cd }F$=d, then
$$
\operatorname{cd }(FG) \leq d+1
$$
and
$$
\operatorname{R}^{d+1}(FG) \cong (\operatorname{R}^{d}F) (\operatorname{R}^{1}G).
$$
\end{lemma}

\begin{proof}
The Grothendieck Spectral Sequence is a convergent first quadrant spectral sequence
$$
E_{2}^{pq}=(\operatorname{R}^{p}F)(\operatorname{R}^{q}G)(\mathcal{A}) \Longrightarrow \operatorname{R}^{p+q}(FG)(\mathcal{A}).
$$
Since cd $G$=1, the only nonzero rows on the $E_{2}$ page are the $q=0$ and $q=1$ rows.  The lemma is a straightforward consequence of these facts, and we leave the rest of the proof as an exercise.
\end{proof}

\begin{lemma} \label{lemma.scheme}
If $X$ is an integral, proper $K$-scheme then $\Gamma(X, {\mathcal{O}}_{X}) \cong K$.
\end{lemma}

\begin{theorem} \label{theorem.cd} 
If $X$ is a smooth, projective variety of dimension $d$, the cohomological dimension of the functor $\Hom_{{\sf Proj }\mathcal{A}}(\pi(\mathcal{O}_{X}(-i) \otimes e_{m}\A),-)$ equals $\operatorname{dim }X + 1$ for all $i,m \in \mathbb{Z}$.
\end{theorem}

\begin{proof}
Since $X$ is projective and irreducible, Lichtenbaum's theorem implies $\operatorname{cd }X=d$ \cite[Theorem 3.1, p.416]{harts}.    

If $\Gamma:{\sf Mod} X \rightarrow {\sf Mod} \Gamma(X, \mathcal{O}_{X})$ is the global sections functor,
\begin{align*}
\Hom_{{\sf Proj }\mathcal{A}}(\pi(\mathcal{O}_{X}(-i) \otimes e_{m}\A),-) & \cong \Hom_{{\sf Gr }\mathcal{A}}(\mathcal{O}_{X}(-i) \otimes e_{m}\A,\omega(-)) \\
& \cong \Homo(\mathcal{O}_{X}(-i), \HomA(e_{m}\A,\omega(-))) \\
& \cong \Homo(\mathcal{O}_{X}(-i), (\omega(-))_{m}) \\
& \cong \Homo(\mathcal{O}_{X}, \mathcal{O}_{X}(i) \otimes (\omega(-))_{m}) \\
& \cong \Gamma \circ \mathcal{O}_{X}(i) \otimes (\omega(-))_{m}
\end{align*}
where we have used Proposition \ref{prop.global} for the second isomorphism and Proposition \ref{prop.homisom} for the third isomorphism.

To complete the proof of the theorem, it suffices, by Lemma \ref{lem.spectral}, to find an object $\mathcal{N}$ in ${\sf Proj}\A$ such that ${\operatorname{R}}^{d}\Gamma (\mathcal{O}_{X}(i) \otimes {\operatorname{R}}^{1}\omega (\mathcal{N})_{m}) \neq 0$.  Suppose $\pi(\M) \cong \mathcal{N}$.  Then we must find an object $\M$ in ${\sf Gr }\A$ such that
$$
{\operatorname{H}}^{d}(X,\mathcal{O}_{X}(i) \otimes R^{1}\omega(\pi \M)_{m}) \neq 0.
$$
Let $\M=\mathcal{O}_{X}(-i) \otimes {\omega}_{X} \otimes e_{m+2}\A$, where ${\omega}_{X}$ is the canonical sheaf on $X$.  By Serre duality on $X$, 
$$
{\operatorname{H}}^{d}(X,\mathcal{O}_{X}(i) \otimes {\operatorname{R}}^{1}\omega(\pi \M)_{m})' \cong \operatorname{Hom}_{{\mathcal{O}}_{X}}(\mathcal{O}_{X}(i) \otimes {\operatorname{R}}^{1}\omega(\pi \M)_{m}, {\omega}_{X}).
$$ 
By Theorem \ref{theorem.exacts}, ${\operatorname{R}}^{1}\omega(\pi \M) \cong {\operatorname{R}}^{2}\tau(\M)$.  By Proposition \ref{prop.tau},
$$
\operatorname{R}^{2}\tau \M \cong \underset{n \to \infty}{\lim} \ExtU^{2}(\quota,\M).
$$   
Thus
\begin{align*} 
\operatorname{Hom}_{{\mathcal{O}}_{X}}(\mathcal{O}_{X}(i) \otimes {\operatorname{R}}^{1}\omega(\pi \M)_{m}, {\omega}_{X}) & \cong \operatorname{Hom}_{{\mathcal{O}}_{X}}(\mathcal{O}_{X}(i) \otimes \underset{n \to \infty}{\lim} \ExtU^{2}(\quota,\M)_{m}, {\omega}_{X}) \\
& \cong \operatorname{Hom}_{\mathcal{O}_{X}}(\underset{n \to \infty}{\lim}\ExtU^{2}(\quota, \M)_{m}, \mathcal{O}_{X}(-i) \otimes {\omega}_{X}).
\end{align*}
By Lemma \ref{lemma.summand}, this is isomorphic to a product of
\begin{equation} \label{eqn.cddd}
\Hom_{{\mathcal{O}}_{X}}(\ExtU^{2}(\A_{0},\M)_{m}, \mathcal{O}_{X}(-i) \otimes \omega_{X})
\end{equation} 
with some other $K$-vector space.  To prove the result, it suffice to show (\ref{eqn.cddd}) is nonzero.  By Theorem \ref{theorem.gor},  
\begin{align*}
\ExtU^{2}(\A_{0},\mathcal{O}_{X}(-i) \otimes \omega_{X} \otimes e_{m+2}\A)_{m}  & \cong {\mathcal{O}}_{X}(-i) \otimes \omega_{X} \otimes \Q_{m}^{*} \\
& \cong \mathcal{O}_{X}(-i) \otimes \omega_{X}  
\end{align*}
as $\mathcal{O}_{X}$-modules.  Thus 
\begin{align*}
\Hom_{{\mathcal{O}}_{X}}(\ExtU^{2}(\A_{0},\M)_{m}, \mathcal{O}_{X}(-i) \otimes \omega_{X}) & \cong \Hom_{\mathcal{O}_{X}}(\mathcal{O}_{X}(-i) \otimes \omega_{X},\mathcal{O}_{X}(-i) \otimes \omega_{X}) \\
& \cong \Hom_{{\mathcal{O}}_{X}}(\mathcal{O}_{X},\mathcal{O}_{X}) \\
& \cong K 
\end{align*}
where the last isomorphism is from Lemma \ref{lemma.scheme}.  Thus (\ref{eqn.cddd}) is nonzero as desired.
\end{proof}

\vfill
\eject

%% file: serre.tex
\section{Serre duality}
Most of the results in this section are modest generalizations of the results in \cite{duality}.
\subsection{Quotient categories and derived functors}

\begin{lemma} \cite[Exercise 1, p. 370]{pop}
Let $\sf{A}$ be an abelian category.  If $\mathcal{M}$ is a noetherian $\sf{A}$-module, the functor $\operatorname{Hom}_{\sf{A}}(\mathcal{M},-)$ commutes with direct sums.
\end{lemma}

\begin{lemma} \label{lem.co}
If $\sf{A}$ is a Grothendieck category, taking (co-)homology of $\sf{A}$-complexes commutes with taking small direct sums of such complexes.
\end{lemma}

\begin{proof}
Let $f_{i}:M_{i} \rightarrow N_{i}$ be a family of morphisms in $\sf{A}$.  Since direct sums are exact in $\sf{A}$, ker $\oplus f_{i} \cong \oplus$ ker $f_{i}$ and coker $\oplus f_{i} \cong \oplus$ coker $f_{i}$.  Thus, $\oplus$ im $f_{i} \cong $ im $\oplus f_{i}$, and the assertion follows. 
\end{proof}

\begin{lemma} \label{lem.inj}
If $\sf{A}$ is a locally noetherian category, and $\sf{C}$ is a localizing subcategory of $\sf{A}$, the direct sum of a small family of injective objects in $\sf{A}/\sf{C}$ is injective.
\end{lemma}

\begin{proof}
Since $\sf{A}$ is locally noetherian, the direct sum of a small family of injective objects is injective \cite[Corollaire 1, p. 358]{Gab}.  Let $\{{\mathcal{A}}_{i}\}$ be a small family of injectives from $\sf{A}/\sf{C}$.  For each $i$, ${\mathcal{A}}_{i} \cong \pi \omega {\mathcal{A}}_{i}$ \cite[Proposition 11.20(5), p. 78]{smith2}.  In addition, $\omega {\mathcal{A}}_{i}$ is torsion-free \cite[Proposition 11.20(1), p. 78]{smith2} and, since $\sf{A}$ has enough injectives \cite[Th$\acute{\mbox{e}}$or$\acute{\mbox{e}}$m 2, p. 362]{Gab}, $\omega {\mathcal{A}}_{i}$ is also injective \cite[Theorem 11.25(1), p. 81]{smith2}.  Thus, $\oplus \omega {\mathcal{A}}_{i}$ is injective.  Since $\omega$ commutes with direct limits \cite[Corollaire 1, p. 379]{Gab}, $\oplus \omega {\mathcal{A}}_{i} \cong \omega \oplus {\mathcal{A}}_{i}$ is also torsion-free.  Thus, $\pi \omega \oplus {\mathcal{A}}_{i} \cong \oplus {\mathcal{A}}_{i}$ is injective \cite[Theorem 11.25(2)]{smith2} as desired.
\end{proof}

From now on, we let $\sf{A}$ denote a locally noetherian category, and we assume a localizing subcategory $\sf{C}$ of $\sf{A}$ has been chosen.  By \cite[Corollaire 1, p.81]{Gab}, the quotient category $\sf{A}/\sf{C}$ is locally noetherian.  We let $\sf{Q}$ denote the quotient $\sf{A}/\sf{C}$.

\begin{corollary} \label{cor.commute}
If $\sf{U}$ is an abelian category and $T:\sf{Q} \rightarrow \sf{U}$ is an additive left exact functor which commutes with small direct sums, then each right-derived functor $\operatorname{R}^{i}T$ also commutes with small direct sums.
\end{corollary}

\begin{corollary} \label{cor.derived}
If $\sf{U}$ is an abelian category and $\M$ is an object of $\sf{A}$ such that $\operatorname{Hom}_{\sf{A}}(\M,-)$ commutes with small direct sums, then each functor
$$
\operatorname{R}^{i}\operatorname{Hom}_{\sf{Q}}(\pi\mathcal{M},-):\sf{Q} \rightarrow \sf{U}
$$  
commutes with small direct sums.
\end{corollary}

\begin{proof}
Since $\mbox{Hom}_{\sf{Q}}(\mathcal{M},-)$ is additive, it suffices, by Corollary \ref{cor.commute} to show that $\mbox{Hom}_{\sf{Q}}(\mathcal{M},-)$ commutes with small direct sums.  Since $\omega$ commutes with direct sums, this follows from \cite[Proposition 11.20(3), p. 78]{smith2}.
\end{proof}

\begin{definition}
Let $R$ be a commutative ring.  Then $\sf{A}$ is {\bf $\mathbf{R}$-linear} if there is a map of rings $\lambda:R \rightarrow \operatorname{Hom}_{\sf A}(\M,\M)$ for all $\M$ in $\sf{A}$ such that if $f:\M \rightarrow \mathcal{N}$ is a morphism in $\sf{A}$, the diagram
$$
\begin{CD}
R @>>> \operatorname{Hom}_{\sf A}(\M,\M) \\
@VVV @VVV \\
\operatorname{Hom}_{\sf A}(\N,\N) @>>> \operatorname{Hom}_{\sf A}(\M,\N)
\end{CD}
$$
whose right vertical and bottom horizontal are induced by $f$ commutes.  
\end{definition}

We omit the straightforward proof of the following
\begin{lemma} \label{lem.linearity}
Suppose $\sf{A}$ is $R$-linear.  The categories ${\sf{Ch}}$ and $\sf{K}$ inherit an $R$-linear structure from $\sf{A}$.  If $R$ is a field, $\sf{D}$ inherits an $R$-linear structure from $\sf{K}$ such that the localization functor $Q:\sf{K} \rightarrow \sf{D}$ is $R$-linear.  
\end{lemma}
From now on, we assume $\sf{A}$ is $K$-linear for some field $K$.

The proof of the following lemma is straightforward, so we omit it.
\begin{lemma}
Let $\sf{B}$ be a category with arbitrary direct sums and let $S$ be a multiplicative system in $\sf{B}$ which is closed under direct sums.  Then the localization of $\sf{B}$ at $S$, ${\sf{B}}_{S}$ inherits direct sums from $\sf{B}$.
\end{lemma}

\begin{corollary} \label{cor.sums}
If $\sf B$ is an abelian category closed under direct sums, $\sf{K}(\sf{B})$ has direct sums and ${\sf{D}}(\sf{B})$ inherits direct sums from ${\sf{K}(\sf{B})}$.
\end{corollary}

\begin{proof}
Direct sums of complexes in $\sf{B}$ descend to direct sums in ${\sf{K}}$ \cite[Lemma 1.1, p.211]{homotopy}.  Since cohomology of complexes commutes with direct sums, quasi-isomorphisms in ${\sf{K}}$ are closed under direct sums and the previous lemma applies.
\end{proof}

\begin{proposition} \label{prop.commutewith}
Let $\sf{U}$ be an abelian category closed under direct sums.  If $T:\sf{Q} \rightarrow \sf{U}$ is a left exact functor with finite cohomological dimension $D$, then ${\sf{D}}(\sf{Q})$ and ${\sf{D}}(\sf{U})$ have direct sums and the derived functor
$$
{\bf R}T:\sf{D}(\sf{Q}) \rightarrow \sf{D}(\sf{U})
$$
exists and preserves small direct sums.
\end{proposition}

\begin{proof}
The first assertion follows from Corollary \ref{cor.sums}, and the proof of the second assertion is virtually identical to the proof of \cite[Proposition 2.1, p. 712]{duality}, relying on Lemma \ref{lem.co} and Corollary \ref{cor.commute}.
\end{proof}

\begin{lemma} \label{lem.sat}
Let $\mathcal{B}$ be an object of $\sf{Q}$ such that 
$$
T(-)= \operatorname{Hom}_{\sf{Q}}(\mathcal{B},-)
$$ 
has finite cohomological dimension, and let ${\sf U}={\sf Mod }K$.  Then the derived functor
$$
{\bf R}T:D(\sf{Q}) \rightarrow D(\sf{U})
$$
exists, preserves small direct sums, and there is an isomorphism of functors
\begin{equation} \label{eqn.isom}
h^{0}{\bf R}T(-) \overset{\cong}{\rightarrow} \operatorname{Hom}_{D(\sf{Q})}(\mathcal{B},-).
\end{equation}
\end{lemma}

\begin{proof}
The first two assertions follow from Proposition \ref{prop.commutewith}.  The proof of the third claim proceeds exactly like the proof of \cite[Lemma 2.3, p. 713]{duality}, where the isomorphism (\ref{eqn.isom}) on $T$-acyclics ${\mathcal{Q}}^{*}$ is the localization map
$$
\operatorname{Hom}_{{\sf K}(\sf{Q})}({\mathcal{B}},{\mathcal{Q}}^{*}) \rightarrow \operatorname{Hom}_{{\sf D}(\sf{Q})}(\mathcal{B},{\mathcal{Q}}^{*}).
$$
This map is a $K$-module homomorphism by Lemma \ref{lem.linearity}.  
\end{proof}

\begin{proposition} \label{prop.main}
If $\sf{Q}$ satisfies (Gen) with the set $\{{\mathcal{B}}_{i,m}\}_{i,m \in \mathbb{Z}}$ and $T(-)=\operatorname{Hom}_{\sf{Q}}(\mathcal{B}_{i,m},-)$ has finite cohomological dimension, then $\sf{D(Q)}$ is compactly generated.
\end{proposition}

\begin{proof}
By Corollary \ref{cor.sums}, $\sf{D(Q)}$ has direct sums.  Since $\sf Q$ satisfies (Gen), and since an acyclic complex is isomorphic to the zero complex in $\sf{D(Q)}$, $\{\B_{i,m}\{n\}\}_{i,m,n \in \mathbb{Z}}$ is a generating set for $\sf{D}(\sf{Q})$.

To complete the proof that $\sf{D}(\sf{Q})$ is compactly generated, we must prove 
$$
\{\B_{i,m}\{n\}\}_{i,m,n \in \mathbb{Z}}
$$ 
is a set of compact objects, i.e., we must prove that
$$
\operatorname{Hom}_{D(\sf{Q})}({\mathcal{B}}_{i,m}\{n\},-)
$$ 
commutes with direct sums.  By Lemma \ref{lem.sat}, ${\bf R}T$ exists and commutes with direct sums.  Thus, $h^{0}{\bf R}T(-)$ preserves direct sums, and so $\operatorname{Hom}_{\sf{D}(\sf{Q})}(\B_{i,m},-)$ preserves direct sums by Lemma \ref{lem.sat}.  By the observation following \cite[Definition 1.6, p. 210]{homotopy}, $\operatorname{Hom}_{\sf{D}(\sf{Q})}(\B_{i,m}\{n\},-)$ is compact as well.  
\end{proof}

\subsection{Serre duality for quotient categories}
In this section we assume $\sf{Q}=\sf{A}/\sf{C}$ satisfies (Gen) with the set $\{{\mathcal{B}}_{i,m}\}_{i,m \in \mathbb{Z}}$ and $T(-)=\operatorname{Hom}_{\sf{Q}}(\mathcal{B}_{i,m},-)$ has finite cohomological dimension $D$ for all $i,m \in \mathbb{Z}$.  

We let $\Gamma(-)=\operatorname{Hom}_{\sf{Q}}(\B_{0,0},-)$, and we let $H^{n}(-)$ denote the $n$'th right derived functor of $\Gamma(-)$.  Finally, we let $\operatorname{Hom}^{*}(-,-)$ denote the complex-hom.  We will sometimes abuse notation in this section by writing $K$ for ${\sf Mod }K$.

\begin{theorem} \label{theorem.brown}
The derived functor $\bf{R} \Gamma$ has a right-adjoint
$$
G:{\sf{D}}(K) \rightarrow {\sf{D}}(\sf{B}).
$$
\end{theorem}

\begin{proof}
Since the hypothesis of Proposition \ref{prop.main} hold, the result follows from the Brown Adjoint Functor Theorem \cite[Theorem 4.1, p. 223]{N}.
\end{proof}

\begin{definition}
We set ${\omega}^{*}=G(K)$, ${\omega}^{\circ}=h^{-D}({\omega}^{*})$ and call ${\omega}^{*}$ the {\bf dualizing complex}.
\end{definition}

\begin{corollary} \label{cor.serre}
There is an isomorphism of abelian groups
$$
\operatorname{Hom}_{{\sf D}(\sf{Q})}({\mathcal{C}}^{*},{\omega}^{*}) \cong \operatorname{Hom}_{K}(h^{0}{\bf{R}}\Gamma{\mathcal{C}}^{*},K)
$$
natural in ${\mathcal{C}}^{*}$.
\end{corollary}

\begin{proof}
By Theorem \ref{theorem.brown} there is a natural isomorphism of abelian groups
$$
\operatorname{Hom}_{{\sf D}(\sf{Q})}({\mathcal{C}}^{*},G{M}^{*}) \cong \operatorname{Hom}_{{\sf{D}}(K)}({\bf R}\Gamma{\mathcal{C}}^{*},{M}^{*}). 
$$
Setting $M^{*}=K$ gives an isomorphism 
$$
\operatorname{Hom}_{{\sf D}(\sf{Q})}({\mathcal{C}}^{*},{\omega}^{*}) \cong \operatorname{Hom}_{{\sf D}(K)}({\bf R}\Gamma{\mathcal{C}}^{*},K). 
$$
It is straightforward to show that $\operatorname{h}^{0}(\operatorname{Hom}_{K}^{*}({\bf R}\Gamma{\mathcal{C}}^{*},K)) \cong \operatorname{Hom}_{{\sf K}(K)}({\bf R}\Gamma{\mathcal{C}}^{*},K)$ as $K$-modules.  On the other hand, the localization map
$$
\operatorname{Hom}_{{\sf K}(K)}({\bf R}\Gamma{\mathcal{C}}^{*},K) \rightarrow \operatorname{Hom}_{{\sf D}(K)}({\bf R}\Gamma{\mathcal{C}}^{*},K)
$$
is an isomorphism of $K$-modules by \cite[Corollary 10.4.7, p.388]{W}.  Thus, to prove the result, we must establish an isomorphism 
$$
h^{0}(\operatorname{Hom}_{K}^{*}({\bf R}\Gamma{\mathcal{C}}^{*},K)) \rightarrow \operatorname{Hom}_{K}(\operatorname{h}^{0}{\bf R}\Gamma{\mathcal{C}}^{*},K).
$$  
Such an isomorphism exists by exactness of the functor $\operatorname{Hom}_{K}(-,K)$.
\end{proof}

\begin{lemma}
The only non-vanishing cohomology of ${\omega}^{*}$ is in degrees $-D,-D+1, \ldots, 0$.  Thus we may assume that ${\omega}^{*}$ has ${\omega}^{-D-r}=0$ for every $r>0$ and that ${\omega}^{*}$ consists of injective $\sf{Q}$-objects.
\end{lemma}

\begin{proof}
By Corollary \ref{cor.serre},
$$
\mbox{Hom}_{{\sf D}(\sf{Q})}({\mathcal{B}}_{i,m}\{-n\},{\omega}^{*}) \cong \mbox{Hom}_{K}(\operatorname{h}^{0}{\bf R}\Gamma({\mathcal{B}}_{i,m}\{-n\}),K).
$$
Since ${\bf R}\Gamma$ is triangulated, it commutes with shifts so that
\begin{align*}
\operatorname{h}^{0}{\bf R}\Gamma({\mathcal{B}}_{i,m}\{-n\}) & \cong \operatorname{h}^{-n}({\bf R}\Gamma({\mathcal{B}}_{i,m})) \\
& \cong \operatorname{R}^{-n}\Gamma({\mathcal{B}}_{i,m}).
\end{align*}
Thus
\begin{align*}
\mbox{Hom}_{K}(\operatorname{h}^{0}{\bf R}\Gamma({\mathcal{B}}_{i,m}\{-n\}),K) & \cong \mbox{Hom}_{K}(\operatorname{R}^{-n}\Gamma({\mathcal{B}}_{i,m}),K). 
\end{align*}
The first claim follows from the fact that $\sf{Q}$ satisfies (Gen) with $\{\B_{i,m}\}_{i,m \in \mathbb{Z}}$.

To prove the second assertion, note that by the first part of the Lemma, ${\omega}^{*}$ is quasi-isomorphic to
$$
\cdots \rightarrow 0 \rightarrow {\omega}^{-D}/B^{-D}({\omega}^{*}) \rightarrow {\omega}^{-D+1} \rightarrow \cdots \rightarrow {\omega}^{0} \rightarrow \cdots
$$
The fact that we can assume ${\omega}^{*}$ consists of injective $\sf{Q}$-objects now follows from the proof of \cite[Theorem 10.4.8, p.388]{W}.
\end{proof}

\begin{lemma}
Let $\M$ be an object of $\sf{Q}$.  The adjointness isomorphism of abelian groups
\begin{equation} \label{eqn.ad1}
\operatorname{Hom}_{{\sf D}(\sf{Q})}(\mathcal{M}\{n\},{\omega}^{*}) \cong \operatorname{Hom}_{{\sf D}(K)}({\bf R}\Gamma \M\{n\},K)
\end{equation}
is an isomorphism which induces an isomorphism of $K$-modules 
\begin{equation} \label{eqn.rmod}
h^{0}\operatorname{Hom}^{*}(\M \{n\},\omega^{*}) \cong \operatorname{h}^{0}\operatorname{Hom}^{*}({\bf R}\Gamma \M \{n\}, K).
\end{equation}
functorial in $\M \{n\}$.
\end{lemma}

\begin{proof}
Since $\M$ is quasi-isomorphic to an injective resolution of $\M$, and since (\ref{eqn.ad1}) is functorial, it suffices to prove the result after replacing $\M\{n\}$ with a left-bounded complex of injectives, $\mathcal{I}^{*}$.  By \cite[Corollary 10.5.11, p.394]{W} and \cite[Existence Theorem 10.5.6, p. 392]{W}, the diagram
$$
\begin{CD}
\operatorname{Hom}_{\sf K}(\mathcal{I}^{*},\omega^{*}) @>{\Gamma}>> \operatorname{Hom}_{\sf K}(\Gamma\mathcal{I}^{*},\Gamma\omega^{*}) @>>> \operatorname{Hom}_{\sf K}(\Gamma \mathcal{I}^{*},K) \\
@VVV @VVV @VVV \\
\operatorname{Hom}_{\sf D}(\mathcal{I}^{*},\omega^{*}) @>>{{\bf R}\Gamma}> \operatorname{Hom}_{\sf D}({\bf R}\Gamma\mathcal{I}^{*},{\bf R}\Gamma\omega^{*}) @>>> \operatorname{Hom}_{\sf D}({\bf R}\Gamma \mathcal{I}^{*},K)
\end{CD}
$$
whose verticals are localizations, whose bottom-right horizontal is induced by the counit of the adjoint pair $({\bf {R}}\Gamma, G)$, and whose top-right horizontal is the map guaranteed to exist by \cite[Corollary 10.4.7, p. 388]{W}, commutes.  Since $\sf K$ and $\Gamma$ are $K$-linear, and since the verticals are isomorphisms by \cite[Corollary 10.4.7, p. 388]{W} the bottom horizontals must be $K$-linear maps since the top verticals are.  Thus, (\ref{eqn.ad1}) is a morphism of $K$-modules.

To complete the proof, we note that 
$$
\operatorname{Hom}_{\sf{D}(\sf{Q})}(\mathcal{M}\{n\},{\omega}^{*}) \overset{Q^{-1}}{\rightarrow} \operatorname{Hom}_{\sf{K}(\sf{Q})}(\mathcal{M}\{n\},{\omega}^{*}) \cong \operatorname{h}^{0}\operatorname{Hom}^{*}(\M\{n\},\omega^{*})
$$
and
$$
\operatorname{Hom}_{{\sf D}(K)}({\bf R}\Gamma\M \{n\},K) \overset{Q^{-1}}{\rightarrow} \operatorname{Hom}_{{\sf K}(K)}({\bf R}\Gamma\M \{n\},K) \cong \operatorname{h}^{0}\operatorname{Hom}^{*}({\bf R}\Gamma \M\{n\}, K) 
$$
are isomorphisms of $K$-modules.
\end{proof}
The proof of the following two Corollaries are essentially the same as the proofs of \cite[Corollary 3.5, p. 715]{duality} and \cite[Corollary 3.6, p. 716]{duality}.
\begin{corollary} \label{cor.serre1}
The dualizing complex ${\omega}^{*}$ has the property that
$$
\operatorname{Hom}_{K}(\operatorname{H}^{n}(\mathcal{M}),K) \cong \operatorname{h}^{-n}\operatorname{Hom}_{\sf{Q}}(\mathcal{M},{\omega^{*}}).
$$
Furthermore, this isomorphism is $K$-linear and functorial in $\M$.
\end{corollary}

\begin{proof}
First, notice that 
$$
\operatorname{Hom}_{K}(\operatorname{H}^{n}(\mathcal{M}),K) \cong (\operatorname{h}^{0}{\bf R}\Gamma \M\{n\})'
$$
as $K$-modules.  On the other hand, (\ref{eqn.rmod}) tells us the right-hand side is isomorphic, as $K$-modules, to $\operatorname{h}^{0}\operatorname{Hom}^{*}(\M\{n\},\omega^{*})$.  
\end{proof}

\begin{corollary} \label{cor.serreya}
Let $r \in \mathbb{N}$, suppose
\begin{equation} \label{eqn.serrya}
\operatorname{H}^{D-1}({\mathcal{B}}_{i,m})= \cdots = \operatorname{H}^{D-r}({\mathcal{B}}_{i,m})=0
\end{equation}
for all $i>>0$ and all $m>>0$.  Then, for $n=D,D-1, \ldots, D-r$ there is a $K$-linear isomorphism
$$
\operatorname{H}^{n}(\M)' \rightarrow \operatorname{Ext}_{\sf{Q}}^{D-n}(\mathcal{M},{\omega^{0}})
$$
functorial in $\mathcal{M}$.
\end{corollary}

\begin{proof}
By Corollary \ref{cor.serre}, 
\begin{align*}
\operatorname{Hom}_{{\sf D}(\sf{Q})}({\mathcal{B}}_{i,m}\{D-n\},\omega^{*}) & \cong \operatorname{h}^{0}{\bf R}\Gamma(\B_{i,m}\{D-n\})' \\
& \cong \operatorname{H}^{D-n}({\mathcal{B}}_{i,m})'.
\end{align*}
Thus, since $\sf{Q}$ satisfies (Gen) with $\{\B_{i,m}\}_{i,m \in \mathbb{Z}}$, (\ref{eqn.serrya}) implies
$$
\operatorname{h}^{-D+1}\omega^{*}= \cdots = \operatorname{h}^{-D+r}\omega^{*}=0.
$$
This implies $\omega^{*}$ is the beginning of $\omega^{0}$'s injective resolution, twisted by $\{D\}$.  Thus, by Corollary \ref{cor.serre1},
$$
\operatorname{H}^{n}(\M)'=\operatorname{h}^{-n}\operatorname{Hom}_{\sf{Q}}(\M,\omega^{*})=\operatorname{Ext}_{\sf{Q}}^{D-n}(\M,\omega^{0}).
$$
\end{proof}

\begin{theorem} \label{theorem.serre} 
Let $X$ be a smooth projective variety of dimension $d$ over a field $K$, and let $\mathcal{A}$ be a non-commutative symmetric algebra generated by a locally free ${\mathcal{O}}_{X}$-bimodule of rank two.  Then there exists an object $\omega_{\A}$ in $\sf{Proj }\A$ such that for $0 \leq i \leq d+1$ there is an isomorphism
$$
\operatorname{Ext}_{{\sf Proj} \mathcal{A}}^{i}(\pi(\mathcal{O}_{X} \otimes e_{0}\mathcal{A}),\mathcal{M})' \cong \operatorname{Ext}_{{\sf Proj} \mathcal{A}}^{d+1-i}(\mathcal{M},\omega_{\A})
$$
natural in $\mathcal{M}$.  The prime denotes dualization with respect to $K$.
\end{theorem}

\begin{proof}
If we set ${\sf A}={\sf Gr}\A$, ${\sf C}={\sf Tors }\A$, $\B_{i,m} = \pi(\mathcal{O}_{X}(-i) \otimes e_{m}\A)$ and $r=D$, Theorem \ref{theorem.cd}, Proposition \ref{prop.gen} and Theorem \ref{theorem.cohom} ensure that the hypothesis of Corollary \ref{cor.serreya} are satisfied.  The result follows. 
\end{proof}

\vfill
\eject

%% file: compats.tex
\section{Compatibilities}
The purpose of this section is to prove various compatibilities between the duality functor $(-)^{*}$ and the tensor product functor $\otimes$.  

Given a diagram of categories, functors, and natural transformations:
$$
\xymatrix{
{\sf{X}}\ruppertwocell^F{\Delta} \rlowertwocell_{F''}{\Delta'} \ar[r]_{\hskip -.2in F'} & {\sf{Y}} \ruppertwocell^G{\Theta} \rlowertwocell_{G''}{\Theta'} \ar[r]_{\hskip -.2in G'} & {\sf{Z}}
}
$$
we have 
\begin{equation} \label{eqn.horiz}
(\Theta' \circ \Theta)*(\Delta'\circ \Delta)=(\Theta' * \Delta') \circ (\Theta * \Delta).
\end{equation}

The proof of the following lemma is straightforward, so we omit it.
\begin{lemma} \label{lemma.commute}
For $i=1,2$, suppose $F_{i}:\sf{C} \rightarrow \sf{D}$ and $G_{i}:\sf{D} \rightarrow \sf{C}$ are functors such that $(F_{i},G_{i}, \eta_{i}, \epsilon_{i})$ is an adjoint pair.  Suppose $\psi:F_{1} \rightarrow F_{2}$ is an isomorphism.  Then $(F_{1},G_{2}, G_{2}*\psi^{-1} \circ \eta_{2}, \epsilon_{2} \circ \psi*G_{2})$ is an adjoint pair, and the isomorphism $\phi:G_{2} \rightarrow G_{1}$ resulting from the uniqueness of right adjoints makes the diagram
$$
\begin{CD}
\Hom_{\sf{D}}(F_{2}\mathcal{A},\mathcal{B}) @>{\cong}>> \Hom_{\sf{C}}(\mathcal{A},G_{2}\mathcal{B}) \\
@V{- \circ \psi_{\A}}VV @VV{\phi_{\B} \circ -}V \\
\Hom_{\sf{D}}(F_{1}\mathcal{A},\mathcal{B}) @>>{\cong}> \Hom_{\sf{C}}(\mathcal{A},G_{1}\mathcal{B}) 
\end{CD}
$$
whose horizontals are adjointness isomorphisms, commute.
\end{lemma}

\begin{lemma} \label{lem.com1}
Keep the notation as in Lemma \ref{lemma.commute}.  If $\sf{C}=\sf{D}$, the diagram
$$
\begin{CD}
F_{1}G_{1} @>{\epsilon_{1}}>> \operatorname{id} @>{\eta_{1}}>> G_{1}F_{1} \\
@V{\psi * \phi^{-1}}VV @VV{=}V @VV{\phi^{-1}*\psi}V \\
F_{2}G_{2} @>>{\epsilon_{2}}> \operatorname{id} @>>{\eta_{2}}> G_{2}F_{2}
\end{CD}
$$
commutes. 
\end{lemma}

\begin{proof}
The commutivity of the diagram follows from the universal property of a right adjoint to $F_{1}$.
\end{proof}

\begin{corollary} \label{cor.tweed}
If $\A$, $\B$ and $\B'$ are coherent, locally free $\mathcal{O}_{X}$-bimodules and $f:\B \rightarrow \B'$ is a morphism, the diagram
$$
\begin{CD}
(\A \otimes \B')^{*} @>{(\A \otimes f)^{*}}>> (\A \otimes \B)^{*} \\
@V{\cong}VV @VV{\cong}V \\
{\B'}^{*} \otimes \A^{*} @>>{f^{*} \otimes \A^{*}}> \B^{*} \otimes \A^{*}
\end{CD}
$$
whose verticals are the maps (\ref{eqn.canonicalisom0}), commutes.
\end{corollary}

\begin{proof}
Let $(A,A^{*})$, $(B,B^{*})$, $(B',{B'}^{*})$, $(C,C^{*})$ and $(C',{C'}^{*})$ denote the adjoint pairs $(- \otimes \A,- \otimes \A^{*})$, $(- \otimes \B, - \otimes \B^{*})$, $(- \otimes \B', - \otimes {\B'}^{*})$, $(-\otimes (\A \otimes \B),- \otimes (A \otimes \B)^{*})$ and $(- \otimes (\A \otimes \B'),- \otimes (\A \otimes \B')^{*})$ respectively.  

Associativity of the bimodule tensor product induces isomorphisms $\alpha:C \rightarrow BA$ and $\beta:C' \rightarrow B'A$.  Thus, $(C,A^{*}B^{*})$ and $(C',A^{*}{B'}^{*})$ are canonically adjoint pairs as in Lemma \ref{lemma.commute}.  By uniqueness of right adjoints, there are isomorphisms $\gamma:C^{*} \rightarrow A^{*}B^{*}$ and $\delta:{C'}^{*} \rightarrow A^{*}{B'}^{*}$.  In the proof that follows, we will abuse notation by using the symbols $\eta$ and $\epsilon$ to denote the canonical units and counits of adjoint pairs.  We will also omit horizontal products with identity functors where no confusion arises.  

To prove the Corollary, we need to show the outer circuit in the diagram
\begin{equation} \label{eqn.adjointss}
\begin{CD}
{C'}^{*} @>{\eta}>> C^{*}C{C'}^{*} @>{f}>> C^{*}C'{C'}^{*} @>{\epsilon}>> C^{*} \\
@V{\delta}VV @VV{\gamma*C*\delta}V @VV{\gamma*D*\delta}V @VV{\gamma}V \\
A^{*}{B'}^{*} @>>{\eta}> A^{*}B^{*}CA^{*}{B'}^{*} @>>{f}> A^{*}B^{*}C'A^{*}{B'}^{*} @>>{\epsilon}> A^{*}B^{*} \\
@V{=}VV @VV{\alpha}V @VV{\beta}V @VV{=}V \\
A^{*}{B'}^{*} & & A^{*}B^{*}BAA^{*}{B'}^{*} @>>{f}> A^{*}B^{*}{B'}AA^{*}{B'}^{*} & & A^{*}B^{*} \\
@V{=}VV @VV{\epsilon}V @VV{\epsilon}V @VV{=}V \\
A^{*}{B'}^{*} @>>{\eta}> A^{*}B^{*}B{B'}^{*} @>>{f}> A^{*}B^{*}{B'}{B'}^{*} @>>{\epsilon}> A^{*}B^{*}
\end{CD}
\end{equation}
commutes.  Thus, it suffices to show each square in (\ref{eqn.adjointss}) commutes.  The upper left square in (\ref{eqn.adjointss}) is the outer circuit of the diagram
$$
\begin{CD}
{C'}^{*} @>{\eta}>> {C}^{*}C{C'}^{*} \\
@V{=}VV @VV{\gamma*C{C'}^{*}}V \\
{C'}^{*} @>>{\eta}> A^{*}B^{*}C{C'}^{*} \\
@V{\delta}VV @VV{\delta}V \\
A^{*}{B'}^{*} @>>{\eta}> A^{*}B^{*}CA^{*}{B'}^{*}.
\end{CD}
$$
The top square commutes by Lemma \ref{lem.com1} while the bottom square commutes by (\ref{eqn.horiz}).  

The upper middle square of (\ref{eqn.adjointss}) commutes by (\ref{eqn.horiz}).

The upper right square in (\ref{eqn.adjointss}) is the outer circuit in the diagram
$$
\begin{CD}
C^{*}C'{C'}^{*} @>{\epsilon}>> C^{*} \\
@V{\gamma}VV @VV{\gamma}V \\
A^{*}B^{*}C'{C'}^{*} @>>{\epsilon}> A^{*}B^{*} \\
@V{\beta}VV @AA{\epsilon}A \\
A^{*}B^{*}{B'}A{C'}^{*} @>>{\delta}> A^{*}B^{*}BAA^{*}{B'}^{*}.
\end{CD}
$$
The upper square in this diagram commutes by (\ref{eqn.horiz}) while the lower square commutes by Lemma \ref{lem.com1}.

The central square in (\ref{eqn.adjointss}) commutes by the naturality of the associativity isomorphisms, while the lower middle square commutes by (\ref{eqn.horiz}).  

Finally, we show the lower left rectangle in (\ref{eqn.adjointss}) commutes.  The proof that the lower right rectangle commutes is similar, so we omit it.  We must show the square
$$
\begin{CD}
A^{*} @>{\eta*A^{*}}>> A^{*}B^{*}CA^{*} \\
@V{\eta}VV @VV{\alpha}V \\
A^{*}B^{*}B @<<{\epsilon}< A^{*}B^{*}BAA^{*}
\end{CD}
$$
commutes.  Since
$$
\begin{CD}
A^{*} @>{\eta*A^{*}}>> A^{*}B^{*}CA^{*} \\
@V{\eta*A^{*}}VV @VV{\alpha}V \\
A^{*}AA^{*} @>>{A^{*}*\eta*AA^{*}}> A^{*}B^{*}BAA^{*}
\end{CD}
$$
commutes by Lemma \ref{lem.com1}, we need only show
$$
\begin{CD}
A^{*} @>{\eta*A^{*}}>> A^{*}AA^{*} \\
@V{\eta}VV @VV{\eta}V \\
A^{*}B^{*}B @<<{\epsilon}< A^{*}B^{*}BAA^{*}
\end{CD}
$$
commutes.  This, in turn, follows from the commutivity of 
$$
\begin{CD}
A^{*} @<{\epsilon}<< A^{*}AA^{*} \\
@V{\eta}VV @VV{\eta}V \\
A^{*}B^{*}B @<<{\epsilon}< A^{*}B^{*}BAA^{*}
\end{CD}
$$
which follows from (\ref{eqn.horiz}).
\end{proof}

\begin{lemma} \label{lemma.tweed3}
For $1 \leq i \leq 4$, let $(F_{i},G_{i},\eta_{i},\epsilon_{i})$ be an adjunction and suppose
$$
\begin{CD}
F_{1} @>{\psi_{12}}>> F_{2} \\
@V{\psi_{13}}VV @VV{\psi_{24}}V \\
F_{3} @>>{\psi_{34}}> F_{4}
\end{CD}
$$
is a commutative diagram of isomorphisms and let $\phi_{ji}:G_{j} \rightarrow G_{i}$ denote the isomorphism of functors resulting from the uniqueness of right adjoints.  Then, for each arrow $\psi_{ij}$ in the diagram, 
$$
(F_{i},G_{j}, G_{j}*{\psi_{ij}}^{-1} \circ \eta_{j}, \epsilon_{j} \circ \psi_{ij}*G_{j})
$$ 
is an adjoint pair and the diagram
$$
\begin{CD}
G_{1} @<{\phi_{21}}<< G_{2} \\
@A{\phi_{31}}AA @AA{\phi_{42}}A \\
G_{3} @<<{\phi_{43}}< G_{4}
\end{CD}
$$
commutes.
\end{lemma}

\begin{proof}
The proof is a consequence of the universal property of the arrows $\phi_{ji}$ and we leave the details to the reader.
\end{proof}

\begin{corollary} \label{corollary.tweed3}
If $\E$, $\F$ and $\G$ are locally free, coherent $\mathcal{O}_{X}$-bimodules, then the diagrams
\begin{equation} \label{eqn.tweed3}
\begin{CD}
\G^{*} \otimes (\E \otimes \F)^{*} @>>> \G^{*} \otimes (\F^{*} \otimes \E^{*}) \\
@VVV @VVV \\
((\E \otimes \F) \otimes \G)^{*} @>>> (\F \otimes \G)^{*} \otimes \E^{*}
\end{CD}
\end{equation}
and
$$
\begin{CD}
\G^{*} \otimes (\E \otimes \F)^{*} @>>> (\G^{*} \otimes \F^{*}) \otimes \E^{*} \\
@VVV @VVV \\
(\E \otimes (\F \otimes \G))^{*} @>>> (\F \otimes \G)^{*} \otimes \E^{*}
\end{CD}
$$
whose arrows are induced by the canonical isomorphism (\ref{eqn.canonicalisom0}), commutes.
\end{corollary}

\begin{proof}
We only prove the first diagram commutes since the proof that the other diagram commutes follows in a similar way.  The diagram of functors from ${\sf Mod }X$ to ${\sf Mod }X$
$$
\begin{CD}
(- \otimes (\E \otimes \F)) \otimes \G @<<< ((- \otimes \E)\otimes \F)\otimes \G \\
@AAA @AAA \\
- \otimes ((\E \otimes \F)\otimes \G) @<<< (- \otimes \E) \otimes (\F \otimes \G)
\end{CD}
$$
whose arrows are associativity isomorphisms, commutes by \cite[Propostion 2.5, p. 442]{translations}.  The arrows in this diagram induce the arrows of the diagram of adjoint functors
$$
\begin{CD}
(- \otimes \G^{*}) \otimes (\E \otimes \F)^{*} @>>> ((- \otimes \G^{*}) \otimes \F^{*}) \otimes \E^{*} \\
@VVV @VVV \\
- \otimes ((\E \otimes \F) \otimes \G)^{*} @>>> (- \otimes (\F \otimes \G)^{*}) \otimes \E^{*}
\end{CD}
$$
as in Lemma \ref{lemma.tweed3}, and hence the diagram commutes by Lemma \ref{lemma.tweed3}.  The fact that (\ref{eqn.tweed3}) commutes follows from the pentagon axiom (see \cite[Propostion 2.5, p. 442]{translations}) and we leave the details to the reader.
\end{proof}

\begin{lemma} \label{lem.tweed}
Let $(A,A^{*}, \eta_{A}, \epsilon_{A})$ and $(B,B^{*}, \eta_{B}, \epsilon_{B})$ be pairs of adjoint functors, where $A,B:{\sf C} \rightarrow {\sf D}$.  If $f:A \rightarrow B$ is a natural transformation and $f^{*}$ is the dual morphism,
$$
B^{*} \overset{\eta_{A}*B^{*}}{\longrightarrow} A^{*}AB^{*} \overset{A^{*}*f*B^{*}}{\longrightarrow} A^{*}BB^{*} \overset{A^{*}*\epsilon_{B}}{\longrightarrow} A^{*}
$$
the diagrams 
\begin{equation} \label{eqn.tweedee}
\begin{CD}
AB^{*} @>{f^{*}}>> AA^{*} \\
@V{f}VV @VV{\epsilon_{A}}V \\
BB^{*} @>>{\epsilon_{B}}> \operatorname{id}
\end{CD}
\end{equation}
and
\begin{equation} \label{eqn.tweedum}
\begin{CD}
\operatorname{id} @>{\eta_{A}}>> A^{*}A \\
@V{\eta_{B}}VV @VV{f}V \\
B^{*}B @>>{f^{*}*B}> A^{*}B
\end{CD}
\end{equation}
commute.
\end{lemma}

\begin{proof}
By definition of $f^{*}$, (\ref{eqn.tweedee}) equals the diagram
$$
\begin{CD}
AB^{*} @>{\eta_{A}}>> AA^{*}AB^{*} @>{f}>> AA^{*}BB^{*} @>{\epsilon_{B}}>> AA^{*} \\
@V{f*B^{*}}VV & & @VV{\epsilon_{A}*BB^{*}}V  @VV{\epsilon_{\A}}V \\
BB^{*} @>>{=}> BB^{*} @>>{=}> BB^{*} @>>{\epsilon_{B}}> \operatorname{id}.
\end{CD}
$$ 
The right commutes by (\ref{eqn.horiz}).  Thus, we must show
$$
\begin{CD}
AB^{*} @>{A*\eta_{A}*B^{*}}>> AA^{*}AB^{*} \\
@V{f*B^{*}}VV @VV{AA^{*}*f*B^{*}}V \\
BB^{*} @<<{\epsilon_{A}*BB^{*}}< AA^{*}BB^{*}
\end{CD}
$$
commutes, i.e.
$$
\begin{CD}
A @>{A*\eta_{A}}>> AA^{*}A \\
@V{f}VV @VV{AA^{*}f}V \\
B @<<{\epsilon_{A}*B}< AA^{*}B
\end{CD}
$$
commutes.  But
$$
\begin{CD}
A @<{\epsilon_{A}*A}<< AA^{*}A \\
@V{f}VV @VV{AA^{*}f}V \\
B @<<{\epsilon_{A}*B}< AA^{*}B
\end{CD}
$$
commutes by (\ref{eqn.horiz}).  Thus, the right square in the diagram
$$
\begin{CD}
A @>{A*\eta_{A}}>> AA^{*}A @>{\epsilon_{A}*A}>> A \\
@V{f}VV @VV{AA^{*}f}V @VV{f}V \\
B @<<{\epsilon_{A}*B}< AA^{*}B @>>{\epsilon_{A}*B}> B
\end{CD}
$$
commutes.  The fact that (\ref{eqn.tweedee}) commutes follows.

The outer circuit of the diagram
$$
\begin{CD}
\operatorname{id} @>{\eta_{A}}>> A^{*}A @>{A^{*}f}>> A^{*}B \\
@V{\eta_{B}}VV @VV{A^{*}A*\eta_{B}}V @AA{A^{*}*\epsilon_{B}*B}A \\
B^{*}B @>>{\eta_{A}*B^{*}B}> A^{*}AB^{*}B @>>{A^{*}*f*B^{*}B}> A^{*}BB^{*}B
\end{CD}.
$$ 
equals (\ref{eqn.tweedum}).  The left commutes by (\ref{eqn.horiz}), while the right commutes since
$$
\begin{CD}
A^{*}A @>{A^{*}*f}>> A^{*}B \\
@V{AA^{*}*\eta_{B}}VV @VV{A^{*}B*\eta_{B}}V \\
A^{*}AB^{*}B @>>{A^{*}*f*B^{*}B}> A^{*}BB^{*}B
\end{CD}
$$
commutes.
\end{proof}
We omit the straightforward proof of the following 
\begin{corollary} \label{corollary.tweed}
If $\A$ and $\B$ are locally free coherent $\mathcal{O}_{X}$-bimodules and $f:\A \rightarrow \B$ is a morphism, the diagrams
$$
\begin{CD}
\B^{*}\otimes \A @>{f^{*}}>> \A^{*}\otimes \A \\
@V{f}VV @VV{\epsilon}V \\
\B^{*}\otimes \B @>>{\epsilon}> \mathcal{O}_{X}
\end{CD}
$$
and 
$$
\begin{CD}
\mathcal{O}_{X} @>{\eta}>> \A \otimes \A^{*} \\
@V{\eta}VV @VV{f}V \\
\B \otimes \B^{*} @>>{f^{*}}> \B \otimes \A^{*}
\end{CD}
$$
commute.
\end{corollary}

\vfill
\eject

%% file: proofs.tex
\section{Proofs}
In this section, all arrows labeled with ``$\cong$" are canonical isomorphisms (\ref{eqn.canonicalisom0}).  In addition, we will sometimes implicitly invoke both the naturality of the associativity of bimodules and the fact that any diagram whose arrows are associativity isomorphisms commutes.
\subsection{Proposition \ref{prop.homo}, associativity}
We prove that the right $\A$-module multiplication on $\HU(\C,\M)$ is associative and compatible with scalar multiplication.  

To show $\mu: \HU(\C,\M)\otimes \A \rightarrow \HU(\C,\M)$ is associative, we must show the diagram
$$
\begin{CD}
(\HU(\C,\M)_{i}\otimes \A_{ij}) \otimes \A_{jk} @>>> \HU(\C,\M)_{i}\otimes \A_{ik} \\
@VVV @VVV \\
\HU(\C,\M)_{j} \otimes \A_{jk} @>>> \HU(\C,\M)_{k}
\end{CD}
$$
whose arrows are multiplications, commutes.  Since $\HU(\C,\M)_{k}$ is a kernel, the commutivity of this diagram follows from the commutivity of the diagram
$$
\begin{CD}
((\underset{l}{\Pi}\M_{l} \otimes \C_{il}^{*}) \otimes \A_{ij}) \otimes \A_{jk} @>>> (\underset{l}{\Pi}\M_{l} \otimes \C_{il}^{*}) \otimes \A_{ik} \\
@VVV @VVV \\
(\underset{l}{\Pi}\M_{l} \otimes \C_{jl}^{*}) \otimes \A_{jk} @>>> \underset{l}{\Pi}\M_{l} \otimes \C_{kl}^{*}
\end{CD}
$$
whose verticals are induced by (\ref{eqn.mu25}).  Expanding the diagram and dropping products, we must establish the commutivity of
\begin{equation} \label{eqn.ass1}
\begin{CD}
(\C_{il}^{*} \otimes \A_{ij}) \otimes \A_{jk} @>{\mu}>> \C_{il}^{*} \otimes \A_{ik} \\
@V{\mu^{*}}VV @VV{\mu^{*}}V \\
((\A_{ij} \otimes \C_{jl})^{*} \otimes \A_{ij}) \otimes \A_{jk} & & (\A_{ik} \otimes \C_{kl})^{*} \otimes \A_{ik} \\
@V{\cong}VV @VV{\cong}V \\
((\C_{jl}^{*} \otimes \A_{ij}^{*}) \otimes \A_{ij}) \otimes \A_{jk} & & (\C_{kl}^{*} \otimes \A_{ik}^{*}) \otimes \A_{ik} \\
@V{\epsilon}VV @VV{\epsilon}V \\
\C_{jl}^{*} \otimes \A_{jk} & & \C_{kl}^{*} \\
@V{\mu^{*}}VV @AA{\epsilon}A \\
(\A_{jk} \otimes \C_{kl})^{*} \otimes \A_{jk} @>>{\cong}> (\C_{kl}^{*} \otimes \A_{jk}^{*}) \otimes \A_{jk}.
\end{CD}
\end{equation}
Consider the diagram
\begin{equation} \label{eqn.ass2}
\begin{CD}
\C_{il}^{*}\A_{ik} @>{\mu^{*}}>> (\A_{ik}\C_{kl})^{*}\A_{ik} @>{\cong}>> (\C_{kl}^{*}\A_{ik}^{*})\A_{ik} @>{\epsilon}>> \C_{kl}^{*} \\
@A{\mu}AA @AA{\mu}A & & &  \\
(\C_{il}^{*}\A_{ij})\A_{jk} @>>{\mu^{*}}> ((\A_{ik}\C_{kl})^{*}\A_{ij})\A_{jk} & & & \\
@V{\mu^{*}}VV @VV{\mu^{*}}V & & & \\
((\A_{ij}\C_{jl})^{*}\A_{ij})\A_{jk} @>>{\mu^{*}}> ((\A_{ij}(\A_{jk}\C_{kl}))^{*}\A_{ij})\A_{jk} & & & & \parallel \\
@V{\cong}VV @VV{\cong}V & & & \\
((\C_{jl}^{*}\A_{ij}^{*})\A_{ij})\A_{jk} @>>{\mu^{*}}> (((\A_{jk}\C_{kl})^{*}\A_{ij}^{*})\A_{ij})\A_{jk} & & & \\
@V{\epsilon}VV @VV{\epsilon}V & & & \\
\C_{jl}^{*}\A_{jk} @>>{\mu^{*}}> (\A_{jk}\C_{kl})^{*}\A_{jk} @>>{\cong}> (\C_{kl}^{*}\A_{jk}^{*})\A_{jk} @>>{\epsilon}> \C_{kl}^{*}.
\end{CD}
\end{equation}
The outer top path of (\ref{eqn.ass2}) (starting at the second row and first column) is the right-hand route of (\ref{eqn.ass1}), while the lower route of (\ref{eqn.ass2}) is the left-hand route of (\ref{eqn.ass1}).  Thus, in order to prove (\ref{eqn.ass1}) commutes, it suffices to show every square in (\ref{eqn.ass2}) commutes.  By the functoriality of $\otimes$ and $(-)^{*}$, the top left square, the second left square and the bottom left square commute.  The third left square commutes by Corollary \ref{cor.tweed}.  Thus, to establish the result, we need to show the right-hand square in (\ref{eqn.ass2}) commutes.  Define $\beta$ as the composition
$$
\A_{ik}^{*} \overset{\mu^{*}}{\rightarrow} (\A_{ij}\A_{jk})^{*} \overset{\cong}{\rightarrow} \A_{jk}^{*}\A_{ij}^{*}.
$$
Since the outer circuit of the diagram
\begin{equation} \label{eqn.ass3}
\begin{CD}
(\A_{ik}\C_{kl})^{*}\A_{ik} @>{\cong}>> (\C_{kl}^{*}\A_{ik}^{*})\A_{ik} @>{\epsilon}>> \C_{kl}^{*} \\
@A{\mu}AA @A{\mu}AA @AA{=}A \\
((\A_{ik}\C_{kl})^{*}\A_{ij})\A_{jk} @>>{\cong}> ((\C_{kl}^{*}\A_{ik}^{*})\A_{ij})\A_{jk} & & \C_{kl}^{*} \\
@V{\mu^{*}}VV @V{\beta}VV @AA{\epsilon}A \\
(((\A_{ij}\A_{jk})\C_{kl})^{*}\A_{ij})\A_{jk} @>>{\cong}> ((\C_{kl}^{*}(\A_{jk}^{*}\A_{ij}^{*}))\A_{ij})\A_{jk} @>>{\epsilon}> (\C_{kl}^{*}\A_{jk}^{*})\A_{jk} \\
@V{\cong}VV @A{\cong}AA @A{\cong}AA \\
(((\A_{jk}\C_{kl})^{*}\A_{ij}^{*})\A_{ij})\A_{jk} @>>{=}> (((\A_{jk}\C_{kl})^{*}\A_{ij}^{*})\A_{ij})\A_{jk} @>>{\epsilon}> (\A_{jk}\C_{kl})^{*}\A_{jk} 
\end{CD}
\end{equation}
is the right-hand square in (\ref{eqn.ass2}), the result will follow from its commutivity.

The upper-left and bottom right squares commute by the functoriality of $\otimes$, the middle-left square commutes by Corollary \ref{cor.tweed}, while the bottom-left square commutes by Corollary \ref{corollary.tweed3}.  We leave it to the reader to check, using the naturality of associativity, that the commutivity of the right square in (\ref{eqn.ass3}) is implied by the commutivity of 
$$
\begin{CD}
\A_{ik}^{*}(\A_{ij}\A_{jk}) @>{\mu}>> \A_{ik}^{*}\A_{ik} \\
@V{\mu^{*}}VV @VV{\epsilon}V \\
(\A_{ij}\A_{jk})^{*}(\A_{ij}\A_{jk}) @>>{\epsilon}> \mathcal{O}_{X}.
\end{CD}
$$
The commutivity of this diagram follows from Corollary \ref{corollary.tweed}.

\subsection{Proposition \ref{prop.homo}, scalar multiplication}
To show that 
$$
\mu:\HU(\C,\M) \otimes \A \rightarrow \HU(\C,\M)
$$ 
is compatible with scalar multiplication, it suffices to show
\begin{equation} \label{eqn.mult1}
\begin{CD}
\C_{il}^{*} \otimes \A_{ii} @>^{\mu^{*}}>> (\A_{ii} \otimes \C_{il})^{*} \otimes \A_{ii} \\
@V{\mu_{\mathcal{O}}}VV @VVV \\
\C_{il}^{*} @<<{\mu_{\mathcal{O}} \circ \epsilon}< \C_{il}^{*} \otimes \A_{ii}^{*} \otimes \A_{ii}.
\end{CD}
\end{equation}
commutes, which will follow from the fact that
\begin{equation} \label{eqn.mult2}
\begin{CD}
\C_{il}^{*} \otimes \mathcal{O}_{X} @>{=}>> \C_{il}^{*} \otimes \A_{ii} \\
@A{\epsilon}AA @VV{\mu^{*}}V \\
(\C_{il}^{*} \otimes \A_{ii}^{*}) \otimes \A_{ii} @<<{\cong}< (\A_{ii} \otimes \C_{il})^{*} \otimes \A_{ii}
\end{CD}
\end{equation}
commutes.  To simplify notation, we write $\C$ for $\C_{il}$ and $\mathcal{O}$ for $\mathcal{O}_{\Delta}=\A_{ii}$.  Expanding the right vertical of (\ref{eqn.mult2}), we must show that the diagram
\begin{equation} \label{eqn.itya}
\begin{CD}
\C^{*}\mathcal{O} @>{\mu_{\mathcal{O}}^{-1}\mathcal{O}}>> (\C^{*}\mathcal{O})\mathcal{O} @>{\eta}>> (\C^{*}((\mathcal{O}\C)(\mathcal{O}\C)^{*}))\mathcal{O} \\
@A{=}AA & & @VV{_{\mathcal{O}}\mu}V \\
\C^{*}\mathcal{O} & & & & (\C^{*}(\C(\mathcal{O}\C)^{*}))\mathcal{O} \\
@A{\epsilon}AA & & @VV{\epsilon}V \\
(\C^{*}\mathcal{O}^{*})\mathcal{O} @<<{\cong}< (\mathcal{O}\C)^{*} \mathcal{O} @<<{_{\mathcal{O}}\mu}< (\mathcal{O}(\mathcal{O}\C)^{*})\mathcal{O}
\end{CD}
\end{equation}
commutes.  Let $\beta$ denote the composition
$$
\begin{CD}
\mathcal{O} @>{\eta}>> \mathcal{O}\mathcal{O}^{*} @>{\mu_{\mathcal{O}}^{-1}\mathcal{O}^{*}}>> \mathcal{O}\mathcal{O}\mathcal{O}^{*} @>{\eta}>> \mathcal{O}\C\C^{*}\mathcal{O}^{*}.
\end{CD}
$$
We leave it as an exercise to check, using the functoriality of $\otimes$ and Corollary \ref{cor.tweed}, to prove that the commutivity of (\ref{eqn.itya}) follows from the commutivity of 
\begin{equation} \label{eqn.mult3}
\begin{CD}
\C^{*}\mathcal{O} @>{\mu_{\mathcal{O}}^{-1}\mathcal{O}}>> (\C^{*}\mathcal{O})\mathcal{O} @>{\beta}>> (\C^{*}((\mathcal{O}(\C\C^{*}))\mathcal{O}^{*}))\mathcal{O} \\
@A{\epsilon}AA & & @VV{_{\mathcal{O}}\mu}V \\
(\C^{*}\mathcal{O}^{*})\mathcal{O} @<<{_{\mathcal{O}}\mu}< ((\mathcal{O}\C^{*})\mathcal{O}^{*})\mathcal{O} @<<{\epsilon}< (\C^{*}((\C\C^{*})\mathcal{O}^{*}))\mathcal{O}.
\end{CD}
\end{equation}
By the functoriality of $\otimes$, the diagram
\begin{equation} \label{eqn.mult35}
\begin{CD}
(\C^{*}\mathcal{O})\mathcal{O} @<{\mu_{\mathcal{O}}^{-1}\mathcal{O}}<< \C^{*}\mathcal{O} & & & & & & & \\
@V{\beta}VV \\
(\C^{*}((\mathcal{O}(\C\C^{*}))\mathcal{O}^{*}))\mathcal{O} @>{_{\mathcal{O}}\mu}>> (\C^{*}((\C\C^{*})\mathcal{O}^{*}))\mathcal{O} @>{\epsilon}>> ((\mathcal{O}\C^{*})\mathcal{O}^{*})\mathcal{O} @>{_{\mathcal{O}}\mu}>> (\C^{*}\mathcal{O}^{*})\mathcal{O} \\
& & @V{\epsilon}VV @VV{\epsilon}V @VV{\epsilon}V \\
& & (\C^{*}(\C\C^{*}))\mathcal{O} @>>{\epsilon}> (\mathcal{O}\C^{*})\mathcal{O} @>>{_{\mathcal{O}}\mu}> \C^{*}\mathcal{O} \\
& & @V{\mu_{\mathcal{O}}}VV @VV{\mu_{\mathcal{O}}}V  \\
& & (\C^{*}\C)\C^{*} @>>{\epsilon}> \mathcal{O}\C^{*} 
\end{CD}
\end{equation}
commutes.  In order to prove (\ref{eqn.mult3}) commutes, we must show the top route of (\ref{eqn.mult35}) equals the identity.  Since (\ref{eqn.mult35}) commutes, it suffices to show its bottom route equals the identity, i.e. that the diagram
\begin{equation} \label{eqn.mult4}
\begin{CD}
\mathcal{O}\C^{*} @>{\mu_{\mathcal{O}}^{-1}}>> \mathcal{O}(\C^{*}\mathcal{O}) @>{_{\mathcal{O}}\mu}>> \C^{*}\mathcal{O} @>{\mu_{\mathcal{O}}^{-1}\mathcal{O}}>> (\C^{*}\mathcal{O})\mathcal{O} \\
@A{\epsilon}AA & & & & @VV{\beta}V \\
\C^{*}(\C\C^{*}) @<<{\mu_{\mathcal{O}}}< (\C^{*}(\C\C^{*}))\mathcal{O} @<<{\epsilon}< (\C^{*}((\C\C^{*})\mathcal{O}^{*}))\mathcal{O} @<<{_{\mathcal{O}}\mu}< (\C^{*}((\mathcal{O}(\C\C^{*}))\mathcal{O}^{*}))\mathcal{O}
\end{CD}
\end{equation}
commutes.  
To prove this, it suffices to show
$$
\begin{CD}
\mathcal{O}\mathcal{O} @>{\mathcal{O}\epsilon^{-1}}>> \mathcal{O}(\mathcal{O}^{*}\mathcal{O}) @>{_{\mathcal{O}}\mu^{-1}\mathcal{O}^{*}\mathcal{O}}>> (\mathcal{O}\mathcal{O})(\mathcal{O}^{*}\mathcal{O}) \\
@V{\eta}VV @AA{\mu_{\mathcal{O}}\mathcal{O}^{*}\mathcal{O}}A @VV{\eta}V \\
(\mathcal{O}\mathcal{O}^{*})\mathcal{O} @>>{_{\mathcal{O}}\mu^{-1}\mathcal{O}^{*}\mathcal{O}}> ((\mathcal{O}\mathcal{O})\mathcal{O}^{*})\mathcal{O} @>>{\eta}> (\mathcal{O}(\C\C^{*}))(\mathcal{O}^{*}\mathcal{O}) 
\end{CD}
$$
commutes.  The right square commutes trivially while the left square commutes since $(\mathcal{O} \epsilon) \circ (\eta \mathcal{O}) = \operatorname{id}_{\mathcal{O}}$.  Let 
$$
\gamma:\C^{*} \otimes (\C \otimes \C^{*}) \rightarrow (\C^{*} \otimes \C) \otimes \C^{*}
$$ 
and 
$$
\delta:(\C^{*}\otimes \mathcal{O}) \otimes \C  \rightarrow \C^{*} \otimes (\mathcal{O} \otimes \C)
$$
denote the canonical associativity isomorphisms.  Since
$$
\begin{CD}
((\C^{*}\mathcal{O})(\C\C^{*}))(\mathcal{O}^{*}\mathcal{O}) @>{\delta}>> (\C^{*}((\mathcal{O}\C)\C^{*}))(\mathcal{O}^{*}\mathcal{O}) \\
@V{\mu_{\mathcal{O}}}VV @VV{_{\mathcal{O}}\mu}V \\
(\C^{*}(\C\C^{*}))(\mathcal{O}^{*}\mathcal{O}) @>>{=}> (\C^{*}(\C\C^{*}))(\mathcal{O}^{*}\mathcal{O})
\end{CD}
$$
commutes, in order to show the top of (\ref{eqn.mult4}) commutes, it suffices, by the naturality of associativity, to show
\begin{equation} \label{eqn.mult5}
\begin{CD}
\C^{*}\mathcal{O} @>{\mu_{\mathcal{O}}^{-1}\mathcal{O}}>> (\C^{*}\mathcal{O})\mathcal{O} @>{\epsilon^{-1}}>> (\C^{*}\mathcal{O})(\mathcal{O}^{*}\mathcal{O}) \\
@A{_{\mathcal{O}}\mu}AA & & @VV{\mu_{\mathcal{O}}^{-1}}V \\
(\mathcal{O} \C^{*}) \mathcal{O} & & & & ((\C^{*}\mathcal{O})\mathcal{O})(\mathcal{O}^{*}\mathcal{O}) \\
@A{\epsilon}AA & & @VV{\eta}V \\
((\C^{*}\C)\C^{*})(\mathcal{O}^{*}\mathcal{O}) @<<{\gamma}< (\C^{*}(\C\C^{*}))(\mathcal{O}^{*}\mathcal{O}) @<<{\mu_{\mathcal{O}}}< ((\C^{*}\mathcal{O})(\C\C^{*}))(\mathcal{O}^{*}\mathcal{O})
\end{CD}
\end{equation}
commutes.

The commutivity of (\ref{eqn.mult5}) is equivalent to the route in
\begin{equation} \label{eqn.mult6}
\begin{CD}
(\C^{*}\mathcal{O})\mathcal{O} @>{\epsilon^{-1}}>> (\C^{*}\mathcal{O})(\mathcal{O}^{*}\mathcal{O}) @>{\mu_{\mathcal{O}}^{-1}}>> ((\C^{*}\mathcal{O})\mathcal{O})(\mathcal{O}^{*}\mathcal{O}) \\
@A{\mu_{\mathcal{O}}^{-1}\mathcal{O}}AA @AA{\mu_{\mathcal{O}}^{-1}\mathcal{O}}A @VV{\eta}V \\
\C^{*}\mathcal{O} @>>{\epsilon^{-1}}> \C^{*}(\mathcal{O}^{*}\mathcal{O}) & & ((\C^{*}\mathcal{O})(\C\C^{*}))(\mathcal{O}^{*}\mathcal{O}) \\
@V{\mu_{\mathcal{O}}^{-1}}VV @VV{\mu_{\mathcal{O}}^{-1}\mathcal{O}}V @VV{\mu_{\mathcal{O}}}V \\
(\C^{*}\mathcal{O})\mathcal{O} @>>{\epsilon^{-1}}> (\C^{*}\mathcal{O})(\mathcal{O}^{*}\mathcal{O}) @>>{\mu_{\mathcal{O}}}> (\C^{*}(\C\C^{*}))\mathcal{O}^{*}\mathcal{O} \\
& & @V{\epsilon}VV @VV{\gamma}V \\
& & (\C^{*}\mathcal{O})\mathcal{O} & & ((\C^{*}\C)\C^{*})(\mathcal{O}^{*}\mathcal{O}) \\
& & @V{=}VV @VV{\epsilon}V \\
& & (\C^{*}\mathcal{O})\mathcal{O} @>>{(_{\mathcal{O}}\mu^{-1} \C^{*}\mathcal{O})\circ (\C^{*}\mathcal{O} \mu_{\mathcal{O}})}> \mathcal{O}(\C^{*}\mathcal{O}) @>>{_{\mathcal{O}}\mu}> \C^{*}\mathcal{O}
\end{CD}
\end{equation}
beginning at the second row and first column and proceeding up and around the outside of the diagram equaling the identity.  Thus, to complete the proof that (\ref{eqn.mult5}) commutes, we show every square in (\ref{eqn.mult6}) commutes, and we show the bottom route of (\ref{eqn.mult6}) is the identity.

We first prove that every square in (\ref{eqn.mult6}) commutes.  The upper and lower left square commutes by the functoriality of $\otimes$.  The commutivity of the upper right square follows, in light of the functoriality $\otimes$, from the commutivity of 
$$
\begin{CD}
\mathcal{O} @<{\mu_{\mathcal{O}}}<< \mathcal{O} \otimes \mathcal{O} \\
@V{\eta}VV @VV{\eta}V \\
\C \otimes \C^{*} @<<{_{\mathcal{O}}\mu}< \mathcal{O} \otimes (\C \otimes \C^{*}).
\end{CD}
$$
Since the maps $\mathcal{O} \otimes \mathcal{O} \overset{\mu_{\mathcal{O}}}{\rightarrow} \mathcal{O}$ and $\mathcal{O} \otimes \mathcal{O} \overset{_{\mathcal{O}}\mu}{\rightarrow} \mathcal{O}$ are equal, this diagram commutes by the naturality of the counit.

The commutivity of the lower right square in (\ref{eqn.mult6}) follows from the commutivity of 
$$
\begin{CD}
(\C^{*}\mathcal{O}) \mathcal{O} @>{\eta \otimes \epsilon^{-1}}>> (\C^{*}(\C \C^{*}))( \mathcal{O}^{*} \mathcal{O}) \\
@V{\C^{*}\mu_{\mathcal{O}}}VV @VV{\gamma}V \\
\C^{*}\mathcal{O} & & ((\C^{*} \C) \C^{*})( \mathcal{O}^{*}\mathcal{O}) \\
@V{=}VV @VV{\epsilon}V \\
\C^{*}\mathcal{O} @>>{_{\mathcal{O}}\mu^{-1}}> \mathcal{O}(\C^{*}\mathcal{O})
\end{CD}
$$
and, hence, from the commutivity of
$$
\begin{CD}
(\C^{*}\mathcal{O}) \mathcal{O} @>{\eta}>> (\C^{*}(\C \C^{*})) \mathcal{O} \\
@V{\C^{*}\mu_{\mathcal{O}}}VV @VV{\gamma}V \\
\C^{*}\mathcal{O} & & ((\C^{*} \C) \C^{*}) \mathcal{O} \\
@V{=}VV @VV{\epsilon}V \\
\C^{*}\mathcal{O} @>>{_{\mathcal{O}}\mu^{-1}}> \mathcal{O}(\C^{*}\mathcal{O}).
\end{CD}
$$
Since the composition
$$
\begin{CD}
\C^{*} @>{\mu_{\mathcal{O}}^{-1}}>> \C^{*}\mathcal{O} @>{\eta}>> \C^{*}(\C \C^{*}) @>{\gamma}>> (\C^{*}\C)\C^{*} @>{\epsilon}>> \mathcal{O}\C^{*} @>{_{\mathcal{O}}\mu}>> \C^{*}
\end{CD}
$$
is the identity, the upper route in this diagram is
$$
\begin{CD}
(\C^{*}\mathcal{O})\mathcal{O} @>{\mu_{\mathcal{O}}}>> \C^{*}\mathcal{O} @>{_{\mathcal{O}}\mu^{-1}}>> \mathcal{O}(\C^{*}\mathcal{O})
\end{CD}
$$
and so the diagram commutes.

Finally, we must show the route in (\ref{eqn.mult6}) starting in the first column and second row and proceeding downwards is the identity.  This is clear, and the proof follows.

\subsection{Proposition \ref{prop.tensor}(1), $\eta$ is a morphism in $\sf{Gr }\A$}
We show
$$
\begin{CD}
\M_{k} \otimes \A_{kl} @>{\eta_{k} \otimes \A_{kl}}>> \HU(\C,\M \Ten \C)_{k} \otimes \A_{kl} \\
@V{\mu}VV @VV{\mu}V \\
\M_{l} @>>{\eta_{l}}> \HU(\C,\M \Ten \C)_{l}
\end{CD}
$$
commutes for all $k$ and $l$.  The above diagram will commute provided the diagram
\begin{equation} \label{eqn.star6}
\begin{CD}
\M_{k} \otimes \A_{kl} @>>> ((\M \Ten \C)_{i} \otimes \C_{ki}^{*}) \otimes \A_{kl} \\
@V{\mu}VV @VVV \\
\M_{l} @>>> (\M \Ten \C)_{i} \otimes \C_{li}^{*}
\end{CD}
\end{equation}
whose right vertical is induced by the composition
\begin{equation} \label{eqn.fix}
\begin{CD}
\C_{ki}^{*} \otimes \A_{kl} @>{\mu^{*}}>> (\A_{kl} \otimes \C_{li})^{*} \otimes \A_{kl} @>{\cong}>> \C_{li}^{*} \otimes \A_{kl}^{*} \otimes \A_{kl} @>{\epsilon}>> \C_{li}^{*}
\end{CD}
\end{equation}
and whose horizontals are induced by the composition
$$
\M_{j} \overset{\eta_{j}}{\rightarrow} \HU(\C,\M \Ten \C)_{l} \rightarrow \underset{i}{\Pi}(\M \Ten \C)_{i} \otimes \C_{li}^{*} \overset{\pi_{j}}{\rightarrow} (\M \Ten \C)_{j} \otimes \C_{lj}^{*} 
$$
(whose second arrow is the canonical inclusion and whose last arrow is projection) commutes for all $i$.  To show (\ref{eqn.star6}) commutes, it suffices to show the diagram
\begin{equation} \label{eqn.eta1}
\begin{CD}
\M_{k}  \A_{kl}  @>{\eta}>> (\M_{k}(\C_{ki}\C_{ki}^{*}))\A_{kl} \\
@V{\mu}VV @VVV \\
\M_{l} & & (\M_{k}\C_{ki})\C_{li}^{*} \\
@V{\eta}VV @VVV \\
(\M_{l}\C_{li})\C_{li}^{*} @>>> (\underset{m}{\oplus}\M_{m}  \C_{mi})\C_{li}^{*} \\
& & @VVV \\
& & (\M \Ten \C)_{i} \C_{li}^{*} 
\end{CD}
\end{equation}
whose upper right vertical is induced by (\ref{eqn.fix}), whose bottom horizontal and middle right vertical are canonical inclusions, and whose bottom vertical is the canonical quotient, commutes.  To show (\ref{eqn.eta1}) commutes, we claim it suffices to prove that the diagram
\begin{equation} \label{eqn.eta2}
\begin{CD}
(\M_{k}  \A_{kl})( \C_{li} \C_{li}^{*})  @>{\eta}>> ((\M_{k}(\C_{ki}\C_{ki}^{*}))\A_{kl})(\C_{li}\C_{li}^{*}) \\
@V{\mu}VV @VVV \\
\M_{l}(\C_{li}\C_{li}^{*}) & & ((\M_{k}\C_{ki})\C_{li}^{*})(\C_{li}\C_{li}^{*}) \\
@V{\eta}VV @VVV \\
((\M_{l}\C_{li})\C_{li}^{*})(\C_{li}\C_{li}^{*}) @>>> ((\underset{m}{\oplus}\M_{m}  \C_{mi})\C_{li}^{*})(\C_{li}\C_{li}^{*}) \\
& & @VVV \\
& & ((\M \Ten \C)_{i} \C_{li}^{*})( \C_{li} \C_{li}^{*}) \\
& & @V{\epsilon}VV \\
& & (\M \Ten \C)_{i}\C_{li}^{*}
\end{CD}
\end{equation}
whose upper four rows are just (\ref{eqn.eta1}) tensored on the right by $\C_{li}\otimes \C_{li}^{*}$, commutes.  For, the vertices of (\ref{eqn.eta1}) may be mapped to the upper seven vertices of (\ref{eqn.eta2}) using maps induced by the unit $\mathcal{O} \overset{\eta}{\rightarrow} \C_{li} \otimes \C_{li}^{*}$.  Since the tensor product is functorial and since the composition
$$
\C_{li}^{*} \overset{\eta}{\rightarrow} \C_{li}^{*} \otimes \C_{li} \otimes \C_{li}^{*} \overset{\epsilon}{\rightarrow} \C_{li}^{*}
$$
is the identity, the claim follows.

Now we prove (\ref{eqn.eta2}) commutes.  To this end, we claim that each square in the diagram
\begin{equation} \label{eqn.eta3}
\begin{CD}
((\M_{k}(\C_{ki}\C_{ki}^{*})\A_{kl})(\C_{li}\C_{li}^{*}) @>{\mu}>> ((\M_{k}(\C_{ki}\C_{ki}^{*}))\C_{ki})\C_{li}^{*} \\
@VVV @VV{\epsilon}V \\
((\M_{k}\C_{ki})\C_{li}^{*})(\C_{li}\C_{li}^{*}) @>>{\epsilon}> (\M_{k}\C_{ki})\C_{li}^{*} \\
@VVV @VVV \\
((\underset{m}{\oplus}\M_{m}\C_{mi})\C_{li}^{*})(\C_{li}\C_{li}^{*}) @>>{\epsilon}> \underset{m}{\oplus}(\M_{m} \C_{mi})\C_{li}^{*} \\
@VVV @VVV \\
((\M \Ten \C)_{i} \C_{li}^{*})(\C_{li}\C_{li}^{*}) @>>{\epsilon}> (\M \Ten \C)_{i}\C_{li}^{*}
\end{CD}
\end{equation}
whose left verticals are the right verticals in (\ref{eqn.eta2}), whose center right vertical is the canonical inclusion, and whose lower right vertical is induced by the canonical quotient, commutes.  For, the bottom two squares commute by the functoriality of $\otimes$.  We leave it as an exercise in the naturality of associativity to prove the commutivity of the top square follows from Corollary \ref{corollary.tweed}.  We also claim that the diagram 
\begin{equation} \label{eqn.eta4}
\begin{CD}
\underset{k}{\oplus}\underset{l}{\oplus} ((\M_{k}(\C_{ki}\C_{ki}^{*}))\A_{kl})\C_{li} @>{\mu}>> \underset{k}{\oplus} (\M_{k} (\C_{ki} \C_{ki}^{*})) \C_{ki} \\
@A{\eta}AA @VV{\epsilon}V \\
\underset{k}{\oplus}\underset{l}{\oplus}(\M_{k}\A_{kl})\C_{li} @>>{\mu}> \underset{k}{\oplus}\M_{k}\C_{ki} \\
@V{\mu}VV @VVV \\
\underset{l}{\oplus}\M_{l}\C_{li} @>>> (\M \Ten \C)_{i}
\end{CD}
\end{equation}
whose unlabeled morphisms are canonical quotients, commutes.  To prove the claim, we note that the bottom square commutes by definition of $\M \Ten \C$, while the top commutes by the functoriality of $\otimes$.  Thus, (\ref{eqn.eta4}) commutes.  

To complete the proof that (\ref{eqn.eta2}) commutes, and hence that (\ref{eqn.eta1}) commutes, we note that the outer circuit of the diagram resulting by putting (\ref{eqn.eta2}) to the left of (\ref{eqn.eta3}) commutes since (\ref{eqn.eta4}) commutes.  Since (\ref{eqn.eta3}) commutes, we may conclude that (\ref{eqn.eta2}) commutes, as desired.

\subsection{Proposition \ref{prop.tensor}(1), $\epsilon$ is a morphism in $\sf{Gr }\A$}

We prove
$$
\begin{CD}
(\HU(\C,\M)\Ten \C)_{k} \otimes \A_{kl} @>{\epsilon_{k} \otimes \A_{kl}}>> \M_{k} \otimes \A_{kl} \\
@V{\mu}VV @VV{\mu}V \\
(\HU(\C,\M)\Ten \C)_{l} @>>{\epsilon_{l}}> \M_{l} 
\end{CD}
$$
commutes.  That is, we show the diagram
\begin{equation} \label{eqn.yetagain}
\begin{CD}
(\underset{m}{\oplus}(\HU(\C,\M))_{m}  \C_{mk})  \A_{kl} @>>> (\underset{m}{\oplus}(\underset{i}{\Pi}\M_{i}  \C_{mi}^{*})  \C_{mk})  \A_{kl} @>>> \M_{k} \A_{kl} \\
@V{\mu}VV & & @VV{\mu}V \\
\underset{m}{\oplus}(\HU(\C,\M))_{m}  \C_{ml} @>>> \underset{m}{\oplus}(\underset{i}{\Pi}\M_{i} \C_{mi}^{*})  \C_{ml} @>>> \M_{l}
\end{CD}
\end{equation}
whose left horizontals are the canonical inclusions and whose right horizontals are induced by the composition
$$
\begin{CD}
(\underset{i}{\Pi}\M_{i} \otimes \C_{mi}^{*}) \otimes \C_{mi} @>{\pi_{j}}>> (\M_{j} \otimes \C_{mj}^{*}) \otimes \C_{mj} @>{\epsilon}>> \M_{j}
\end{CD}
$$
commutes.  Let 
$$
\alpha: (\HU(\C,\M)_{m}\otimes \C_{mk})\otimes \A_{kl} \rightarrow ((\underset{i}{\Pi}\M_{i} \otimes \C_{mi}^{*}) \otimes \C_{mk}) \otimes \A_{kl}
$$
be the canonical inclusion.  We first reduce the commutivity of (\ref{eqn.yetagain}) to the commutivity of the diagram
\begin{equation} \label{eqn.yetagain0}
\begin{CD}
(\HU(\C,\M)_{m}\C_{mk})\A_{kl} \\
@V{\alpha}VV \\
((\underset{i}{\Pi}\M_{i}\C_{mi}^{*})\C_{mk})\A_{kl} @>{\pi_{l}}>> ((\M_{l}\C_{ml}^{*})\C_{mk})\A_{kl} \\
@V{\pi_{l}}VV @VVV \\
((\M_{l}\C_{ml}^{*})\C_{mk})\A_{kl} @>>> (((\M_{l}\A_{kl}^{*})\C_{mk}^{*})\C_{mk})\A_{kl}\\
@V{\mu}VV @VV{\epsilon}V \\
(\M_{l}\C_{ml}^{*})\C_{ml} @>>{\epsilon}> \M_{l} 
\end{CD}
\end{equation}
whose upper right vertical and middle horizontal are induced by the composition
\begin{equation} \label{eqn.yetagainmorph}
\begin{CD}
\C_{ml}^{*} @>{\mu^{*}}>> (\C_{mk} \otimes \A_{kl})^{*} @>{\cong}>> \A_{kl}^{*} \otimes \C_{mk}^{*}.
\end{CD}
\end{equation}
To this end, we note that, by definition of $\HU(\C,\M)$, the diagram
\begin{equation} \label{eqn.yetagain1}
\begin{CD}
(\HU(\C,\M)_{m}\C_{mk})\A_{kl} @>{\alpha}>> ((\underset{i}{\Pi}\M_{i}\C_{mi}^{*})\C_{mk})\A_{kl} \\
@V{\alpha}VV @VV{\pi_{k}}V \\
((\underset{i}{\Pi} \M_{i} \C_{mi}^{*})\C_{mk})\A_{kl} & & ((\M_{k} \C_{mk}^{*})\C_{mk})\A_{kl} \\
@V{\pi_{l}}VV @VVV \\
((\M_{l}\C_{ml}^{*})\C_{mk})\A_{kl} @>>> (((\M_{l}\A_{kl}^{*})\C_{mk}^{*})\C_{mk})\A_{kl}
\end{CD}
\end{equation}
whose bottom horizontal is induced by (\ref{eqn.yetagainmorph}) and whose bottom right vertical is induced by
\begin{equation} \label{eqn.epsilon5}
\begin{CD}
\M_{i} @>{\eta}>> \M_{i} \otimes (\A_{ij} \otimes \A_{ij}^{*}) @>{\mu}>> \M_{j} \otimes \A_{ij}^{*}
\end{CD}
\end{equation}
commutes.  In addition, the diagram
\begin{equation} \label{eqn.yetagain2}
\begin{CD}
((\M_{k}\C_{mk}^{*})\C_{mk})\A_{kl} @>{\epsilon}>> \M_{k}\A_{kl} \\
@VVV @VVV \\
(((\M_{l}\A_{kl}^{*})\C_{mk}^{*})\C_{mk})\A_{kl} @>>{\epsilon}> (\M_{l}\A_{kl}^{*})\A_{kl}
\end{CD}
\end{equation}
whose verticals are induced by (\ref{eqn.epsilon5}), commutes by the functoriality of $\otimes$.  Finally, the diagram
\begin{equation} \label{eqn.yetagain3}
\begin{CD}
\M_{k}\A_{kl} @>{\mu}>> \M_{l} \\
@VVV @VV{=}V \\
(\M_{k}\A_{kl}^{*})\A_{kl} @>>{\epsilon}> \M_{l}
\end{CD}
\end{equation}
whose left vertical is induced by (\ref{eqn.epsilon5}), commutes since the diagram
$$
\begin{CD}
\M_{k}\A_{kl} @>{\mu}>> \M_{l} \\
@A{\epsilon}AA @AA{\epsilon}A \\
(\M_{k}\A_{kl})(\A_{kl}^{*}\A_{kl}) @>>{\mu}> \M_{l}(\A_{kl}^{*}\A_{kl})
\end{CD}
$$
commutes by the functoriality of $\otimes$.  Since (\ref{eqn.yetagain1}), (\ref{eqn.yetagain2}) and (\ref{eqn.yetagain3}) commute, the outer circuit of the diagram we get by placing (\ref{eqn.yetagain3}) to the right of (\ref{eqn.yetagain2}) and (\ref{eqn.yetagain2}) to the right of (\ref{eqn.yetagain1}) commutes, i.e. the diagram
$$
\begin{CD}
(\HU(\C,\M)_{m}\C_{mk})\A_{kl} @>{\alpha}>> ((\underset{i}{\Pi}\M_{i}\C_{mi}^{*})\C_{mk})\A_{kl} \\
@V{\alpha}VV @VV{\pi_{k}}V \\
((\underset{i}{\Pi}\M_{i}\C_{mi}^{*})\C_{mk})\A_{kl} & & ((\M_{k}\C_{mk}^{*})\C_{mk})\A_{kl} \\
@V{\pi_{l}}VV @VV{\epsilon}V \\
((\M_{l}\C_{ml}^{*})\C_{mk})\A_{kl} & & \M_{k}\A_{kl} \\
@VVV @VV{\mu}V \\
((\M_{l}(\A_{kl}^{*}\C_{mk}^{*})\C_{mk})\A_{kl} & & \M_{l} \\
@V{\epsilon}VV @VV{=}V \\
(\M_{l}\A_{kl}^{*})\A_{kl} @>>{\epsilon}> \M_{l}
\end{CD}
$$
commutes.  Thus, to prove (\ref{eqn.yetagain}) commutes, it suffices to prove that the outer circuit of (\ref{eqn.yetagain0}) commutes.  We note that the top square in (\ref{eqn.yetagain0}) trivially commutes.  We leave it as an exercise in the naturality of associativity to show the commutivity of the bottom square in (\ref{eqn.yetagain0}) follows from Corollary \ref{corollary.tweed}.

\subsection{Theorem \ref{theorem.gor}, $d_{1}$ is an epimorphism}
Let $\E_{j}=\A_{j,j+1}$ so that $\E_{j+1}=\E_{j}^{*}$.  Let $\alpha:\Q_{i} \rightarrow \E_{i} \otimes \E_{i+1}$ denote the canonical inclusion, and let $\beta$ denote the isomorphism
$$
\begin{CD}
\E_{i+1} @>{\eta}>> \E_{i+1}(\Q_{i}\Q_{i}^{*}) @>{\alpha}>> \E_{i+1}((\E_{i}\E_{i+1})\Q_{i}^{*}) @>{\epsilon}>> \E_{i+1}\Q_{i}^{*}.
\end{CD}
$$
Finally, let $\psi:\Q_{i} \otimes e_{i+2}\A \rightarrow \E_{i} \otimes e_{i+1}\A$ denote the morphism of right $\A$-modules whose $k$th component is the composition
$$
\begin{CD}
\Q_{i} \otimes \A_{i+2,k} @>{\alpha}>> (\E_{i} \otimes \E_{i+1}) \otimes \A_{i+2,k} @>{\mu}>> \E_{i} \otimes \A_{i+1,k}.
\end{CD}
$$
We prove the morphism
$$
\psi':\HomA(\E_{i} \otimes e_{i+1}\A,e_{l}\A) \rightarrow \HomA(\Q_{i} \otimes e_{i+2}\A,e_{l}\A)
$$
induced by $\psi$ is an epimorphism.  Let $\phi$ be the epimorphism
$$
\begin{CD}
\A_{l,i+1}\E_{i}^{*} @>{\eta}>> (\A_{l,i+1} \E_{i}^{*})  (\Q_{i}  \Q_{i}^{*}) @>{\beta^{-1}}>> \A_{l,i+1}\E_{i+1} \Q_{i}^{*} @>{\mu}>> \A_{l,i+2}\Q_{i}^{*}.
\end{CD}
$$
We show the diagram  
\begin{equation} \label{eqn.couldbe}
\begin{CD}
\A_{l,i+1} \E_{i}^{*} @>>> \HomA(e_{i+1}\A,e_{l}\A)  \E_{i}^{*} @>>> \HomA(\E_{i} e_{i+1}\A,e_{l}\A) \\
@V{\phi}VV & & @VV{\psi'}V \\
\A_{l,i+2}  \Q_{i}^{*} @>>> \HomA(e_{i+2}\A,e_{l}\A)  \Q_{i}^{*} @>>> \HomA(\Q_{i}  e_{i+2}\A,e_{l}\A)
\end{CD}
\end{equation}
whose right horizontals are the isomorphisms in Theorem \ref{theorem.hom} (4) and whose left horizontals are induced by the composition of the unit of the adjoint pair 
$$
(- \Ten \A, \HU(\A,-))
$$ 
with the isomorphism
$$
\HU(\A,e_{l}\A \Ten \A) \rightarrow \HU(\A,e_{l}\A)
$$
induced by the multiplication isomorphism
$$
e_{l}\A \Ten \A \rightarrow e_{l}\A
$$
from Proposition \ref{prop.multiso}, commutes. 

The diagram
$$
\begin{CD}
\A_{l,i+1}\E_{i}^{*} @>{\eta_{i+1}}>> \HU(\A,e_{l}\A\Ten \A)_{i+1} \E_{i}^{*} \\
@VVV \\
(\A_{l,i+1}\E_{i}^{*})(\Q_{i}\Q_{i}^{*}) \\
@VVV @VV{\beta}V \\
(\A_{l,i+1}\E_{i}^{*})((\E_{i}\E_{i+1})\Q_{i}^{*}) \\
@VVV \\
\A_{l,i+1} (\E_{i+1} \Q_{i}^{*}) @>{\eta_{i+1}}>> \HU(\A,e_{l}\A\Ten\A)_{i+1}(\E_{i+1}\Q_{i}^{*}) \\
@V{\mu}VV @VV{\mu}V \\
\A_{l,i+2}\Q_{i}^{*} @>>{\eta_{i+1}}> \HU(\A,e_{l}\A\Ten \A)_{i+2} \Q_{i}^{*} 
\end{CD}
$$
whose left column is $\phi$, commutes since the counit of $(-\Ten \A,\HU(\A,-))$ is an $\A$-module morphism.  In addition, the diagram
$$
\begin{CD}
\HU(\A,e_{l}\A\Ten \A)_{i+1} \E_{i}^{*} @>>> \HU(\A,e_{l}\A)_{i+1}\E_{i}^{*} \\
@V{\beta}VV @VV{\beta}V \\
\HU(\A,e_{l}\A\Ten\A)_{i+1}(\E_{i+1}\Q_{i}^{*}) @>>> \HU(\A,e_{l}\A)_{i+1}(\E_{i+1}\Q_{i}^{*}) \\
 @V{\mu}VV @VV{\mu}V \\
\HU(\A,e_{l}\A\Ten \A)_{i+2} \Q_{i}^{*} @>>> \HU(\A,e_{l}\A)_{i+2}\Q_{i}^{*}
\end{CD}
$$
whose horizontals are the isomorphisms from Theorem \ref{theorem.hom} (4), commutes by functoriality of $\HU(\A,-)$ and $\otimes$.  Thus, to prove (\ref{eqn.couldbe}) commutes, it suffices to prove the diagram
\begin{equation} \label{eqn.last0}
\begin{CD}
\HomA(e_{i+1}\A,e_{l}\A) \E_{i}^{*} @>>> \HomA(\E_{i} e_{i+1}\A,e_{l}\A) \\
@V{\beta}VV @VVV \\
\HomA(e_{i+1}\A,e_{l}\A) ( \E_{i+1} \Q_{i}^{*}) & & \HomA(\E_{i} (\E_{i+1}e_{i+2}\A),e_{l}\A) \\
@V{\mu}VV @VV{\alpha}V \\
\HomA(e_{i+2}\A,e_{l}\A)  Q_{i}^{*} @>>> \HomA(\Q_{i}e_{i+2}\A,e_{l}\A)
\end{CD}
\end{equation}
whose upper right vertical is induced by multiplication $\E_{i+1} \otimes e_{i+2}\A \rightarrow e_{i+1}\A$ and whose horizontals are induced by the isomorphisms in Theorem \ref{theorem.hom} (4), commutes.  Let $\E=\E_{i}$ and let $\Q=\Q_{i}$.  To prove (\ref{eqn.last0}) commutes, it suffices to show, by the definition of $\HU(-,-)$, that the diagram
\begin{equation} \label{eqn.last1}
\begin{CD}
\A_{i+1,j}^{*}  \E^{*} @>{\cong}>> (\E \A_{i+1,j})^{*} \\
@V{\eta}VV @VV{\mu^{*}}V \\
\A_{i+1,j}^{*} (\E^{*} (\Q  \Q^{*})) & & (\E (\E^{*}  \A_{i+2,j}))^{*} \\
@V{\alpha}VV @VV{\alpha^{*}}V \\
\A_{i+1,j}^{*} (\E^{*} (( \E  \E^{*})  \Q^{*})) & & (\Q \A_{i+2,j})^{*} \\
@V{\epsilon}VV @AA{\cong}A \\
\A_{i+1,j}^{*} ( \E^{*}  \Q^{*}) & & \A_{i+2,j}^{*}  \Q^{*} \\
@V{\mu^{*}}VV @AA{\epsilon}A \\
(\E^{*} \A_{i+2,j})^{*} ( \E^{*} \Q^{*}) @>>{\cong}> (\A_{i+2,j}^{*}  \E^{**}) \otimes (\E^{*}  \Q^{*}) 
\end{CD}
\end{equation} 
commutes.  In order to prove that (\ref{eqn.last1}) commutes, we notice that (\ref{eqn.last1}) is the outer circuit of the diagram consisting of the diagram
\begin{equation} \label{eqn.last2}
\begin{CD}
\A_{i+1,j}^{*}\E^{*} @>{\cong}>> (\E\A_{i+1,j})^{*} @>{\mu^{*}}>> (\E(\E^{*}\A_{i+2,j}))^{*} \\
@V{\eta}VV @VV{\eta}V @VV{\cong}V \\
(\A_{i+1,j}^{*}( \E^{*} (\Q \Q^{*})) @>>{\cong}> (\E\A_{i+1,j})^{*}(\Q\Q^{*}) & & (\E^{*}\A_{i+2,j})^{*} \E^{*} \\
@V{\alpha}VV @VV{\alpha}V @VV{\eta}V \\
(\A_{i+1,j}^{*}(\E^{*}((\E\E^{*})\Q^{*})) @>>{\cong}> (\E\A_{i+1,j})^{*}((\E\E^{*})\Q^{*}) & & (\E^{*}\A_{i+2,j})^{*}(\E^{*} (\Q\Q^{*})) \\
@V{\epsilon}VV @VV{\mu^{*}}V  @VV{=}V \\
\A_{i+1,j}^{*}(\E^{*}\Q^{*}) & & (\E(\E^{*}\A_{i+2,j}))^{*}((\E\E^{*})\Q^{*}) & & (\E^{*}\A_{i+2,j})^{*}(\E^{*}( \Q\Q^{*})) \\
@V{\mu^{*}}VV @VV{\cong}V @VV{\alpha}V \\
(\E^{*}\A_{i+2,j})^{*} (\E^{*} \Q^{*}) @>>{\eta}> (\E^{*}\A_{i+2,j})^{*}((\E^{*}(\E\E^{*}))\Q^{*}) @>>{=}> (\E^{*}\A_{i+2,j})^{*}((\E^{*}(\E\E^{*}))\Q^{*}) 
\end{CD}
\end{equation}
to the left of the diagram  
\begin{equation} \label{eqn.last3}
\begin{CD}
(\E (\E^{*}\A_{i+2,j}))^{*} @>{\alpha^{*}}>> (\Q \A_{i+2,j})^{*} \\
@V{\cong}VV @VV{\cong}V \\
(\E^{*}\A_{i+2,j})^{*}\E^{*}  & & \A_{i+2,j}^{*}\Q^{*} \\
@V{\eta}VV @AA{\epsilon}A\\
(\E^{*}\A_{i+2,j})^{*} (\E^{*}(\Q\Q^{*})) & & (\A_{i+2,j}^{*}\E^{**})(\E^{*}\Q^{*}) \\
@V{\alpha}VV @AA{\cong}A \\
(\E^{*}\A_{i+2,j})^{*} ((\E^{*}(\E\E^{*}))\Q^{*}) @>>{\epsilon}> (\E^{*}\A_{i+2,j})^{*} (\E^{*}\Q^{*}). \\
\end{CD}
\end{equation}
Thus, to prove that (\ref{eqn.last1}) commutes, we need only prove (\ref{eqn.last2}) and (\ref{eqn.last3}) commute.  We first prove (\ref{eqn.last2}) commutes.

The upper left and middle left squares in (\ref{eqn.last2}) commute by the functoriality of $\otimes$.  To prove the  right rectangle in (\ref{eqn.last2}) commutes, we decompose it into the diagram
$$
\begin{CD}
(\E\A_{i+1,j})^{*} @>{\mu^{*}}>> (\E(\E^{*}\A_{i+2,j}))^{*} @>{\cong}>> (\E^{*}\A_{i+2,j})^{*}\E^{*} \\
@V{\eta}VV @VV{\eta}V @VV{\eta}V \\
(\E\A_{i+1,j})^{*}(\Q\Q^{*}) @>>{\mu^{*}}> (\E(\E^{*}\A_{i+2,j}))^{*}(\Q\Q^{*}) @>>{\cong}> (\E^{*}\A_{i+2,j})^{*} (\E^{*}(\Q\Q^{*})) \\
@V{\alpha}VV @VV{\alpha}V @VV{\alpha}V \\
(\E\A_{i+1,j})^{*}((\E\E^{*})\Q^{*}) @>>{\mu^{*}}> (\E(\E^{*}\A_{i+2,j}))^{*}((\E\E^{*})\Q^{*}) @>>{\cong}>  (\E^{*}\A_{i+2,j})^{*}((\E^{*}(\E\E^{*}))\Q^{*})
\end{CD}
$$
whose squares commute by the functoriality of $\otimes$.  To prove the lower left square in (\ref{eqn.last2}) commutes, it suffices, by the naturality of associativity, to show both squares in
$$
\begin{CD}
(\E\A_{i+1,j})^{*}\E @>{\mu^{*}}>> (\E(\E^{*}\A_{i+2,j}))^{*}\E \\
@V{\cong}VV @VV{\cong}V \\
(\A_{i+1,j}^{*}\E^{*})\E @>>{\mu^{*}}> (\E^{*}\A_{i+2,j})^{*}(\E^{*}\E) \\
@V{\epsilon}VV @VV{\epsilon}V \\
\A_{i+1,j}^{*} @>>{\mu^{*}}> (\E^{*}\A_{i+2,j})^{*}
\end{CD}
$$
commute.  The lower square commutes by functoriality of $\otimes$ while the upper square commutes by Corollary \ref{cor.tweed}.

Thus, in order to prove (\ref{eqn.last1}) commutes, it suffices, in light of the commutivity of (\ref{eqn.last2}), to prove that (\ref{eqn.last3}) commutes.  We leave it as an exercise in the naturality of associativity to show that Corollary \ref{cor.tweed} implies that the commutivity of (\ref{eqn.last3}) follows from the commutivity of the diagram
$$
\begin{CD}
\A_{i+2,j}^{*}(\E\E^{*})^{*} @>{\alpha^{*}}>> \A_{i+2,j}^{*}\Q^{*} \\
@V{\cong}VV \\
\A_{i+2,j}^{*}(\E^{**}\E^{*}) \\
@V{\eta}VV @AA{\epsilon}A \\
(\A_{i+2,j}^{*}(\E^{**}\E^{*}))(\Q\Q^{*}) \\
@V{\alpha}VV \\
(\A_{i+2,j}^{*}(\E^{**}\E^{*}))((\E\E^{*})\Q^{*}) @>>{\epsilon}> \A_{i+2,j}^{*}((\E^{**}\E^{*})\Q^{*}).
\end{CD}
$$
The commutivity of this diagram follows, by the naturality of associativity, from the commutivity of
$$
\begin{CD}
(\E\E^{*})^{*} @>{\alpha^{*}}>> \Q^{*} \\
@V{\cong}VV @AA{=}A \\
\E^{**}\E^{*} & & \Q^{*} \\
@V{\eta}VV @AA{\epsilon}A \\
(\E^{**}\E^{*})(\Q\Q^{*}) @>>{\alpha}> (\E^{**}\E^{*})((\E \E^{*})\Q^{*}).
\end{CD}
$$
The commutivity of this diagram follows from the commutivity of each of the squares in the diagram
$$
\begin{CD}
(\E\E^{*})^{*} @>{\alpha^{*}}>> \Q^{*} @>{=}>> \Q^{*} \\
@V{\eta}VV @AA{\epsilon}A @AA{\epsilon}A \\
(\E\E^{*})^{*}(\Q\Q^{*}) @>>{\alpha}> (\E\E^{*})^{*}((\E\E^{*})\Q^{*}) & & ((\E^{**}(\E^{*}\E))\E^{*})\Q^{*} \\
@V{\cong}VV @VV{\cong}V @AAA \\
(\E^{**}\E^{*})(\Q\Q^{*}) @>>{\alpha}> (\E^{**}\E^{*})((\E\E^{*})\Q^{*}) @>>> (((\E^{**}\E^{*})\E)\E^{*})\Q^{*} 
\end{CD}
$$
whose bottom right horizontal and bottom right vertical are associativity isomorphisms.  The upper left square commutes by definition of the dual morphism while the lower left square commutes by the functoriality of $\otimes$.  The right hand square commutes by Lemma \ref{lem.com1}, and the result follows.
\vfill
\eject

%% file: bib.tex
\address{Adam Nyman, Department of Mathematical Sciences, Mathematics Building, University of Montana, Missoula, MT 59812-0864}